\documentclass[a4paper,11pt,twoside]{article}
\usepackage[english]{babel}
\usepackage[utf8]{inputenc}

\usepackage[a4paper,inner=3cm,outer=3cm,top=4cm,bottom=4cm,pdftex]{geometry}
\usepackage{fancyhdr}
\usepackage{sectsty}
\usepackage{titlesec}
\usepackage{color}
\usepackage{bold-extra}
\usepackage{ mathrsfs }
\titleformat{\section}{\normalfont\scshape\centering}{\thesection}{1em}{}
  \titleformat{\subsection}{\bfseries}{\thesubsection}{1em}{}
  
\usepackage{comment}
\usepackage{graphics}
\usepackage{aliascnt}
\usepackage[hidelinks]{hyperref}

\usepackage{amsmath}
\usepackage{amsfonts}
\usepackage{amssymb}
\usepackage{amsthm}
\usepackage{comment}
\newcommand{\qedd}{\hfill \ensuremath{\Box}}
\newtheorem{theorem}{Theorem}[section]

\newtheorem{hypothesis}[theorem]{Hypothesis}
\newtheorem{lemma}[theorem]{Lemma}
\newtheorem{proposition}[theorem]{Proposition}
\theoremstyle{definition}
\newtheorem{definition}[theorem]{Definition}
\newtheorem{remark}[theorem]{Remark}
\numberwithin{equation}{section}

\title{\scshape \large THE GOLDBACH PROBLEM FOR PRIMES THAT ARE SUMS OF TWO SQUARES PLUS ONE}
\author{\scshape \normalsize Joni Ter\"av\"ainen}  
\date{}
\begin{document}
\maketitle

\begin{abstract} \hspace{-0.55cm}We study the Goldbach problem for primes represented by the polynomial $x^2+y^2+1$. The set of such primes is sparse in the set of all primes, but the infinitude of such primes was established by Linnik. We prove that almost all even integers $n$ satisfying certain necessary local conditions are representable as the sum of two primes of the form $x^2+y^2+1$. This improves a result of Matom\"aki, which tells that almost all even $n$ satisfying a local condition are the sum of one prime of the form $x^2+y^2+1$ and one generic prime. We also solve the analogous ternary Goldbach problem, stating that every large odd $n$ is the sum of three primes represented by our polynomial. As a byproduct of the proof, we show that the primes of the form $x^2+y^2+1$ contain infinitely many three term arithmetic progressions, and that the numbers $\alpha p \hspace{-0.1cm} \pmod 1$ with $\alpha$ irrational and $p$ running through primes of the form $x^2+y^2+1$, are distributed rather uniformly.
\end{abstract}

\section{Introduction}

Let $\mathscr{P}$ be the set of primes represented by the quadratic polynomial $x^2+y^2+1$. We consider the Goldbach problem for the set $\mathscr{P}$, our main result being the following.

\begin{theorem}\label{theo_goldbach} 
Almost all even positive integers $n\not \equiv 5,8\hspace{-0.1cm} \pmod{9}$ can be represented as $n=p+q$ with $p,q \in \mathscr{P}$. 
\end{theorem}

By ''almost all'' we mean that the number of exceptional $n\leq N$ is $o(N)$. The local condition $n\not\equiv 5,8\hspace{-0.1cm} \pmod{9}$ is necessary (unless $p$ or $q$ equals $3$ in which case we can only represent $o(N)$ integers), as is easily seen by considering primes of the form $x^2+y^2+1$ modulo $9$. An earlier result of Matom\"aki \cite{matomaki-goldbach}, using a somewhat different method, showed that one of the primes $p$ and $q$ can be taken to be from $\mathscr{P}$, the other one being a generic prime. A few years later, Tolev \cite{tolev_binarygoldbach} gave an asymptotic formula for a weighted count of the representations $n=p+q$ with $p\in \mathscr{P}$ and $q$ a generic prime for almost all even $n$. Naturally, there is a close connection between the almost all version of the binary Goldbach problem and the ternary Goldbach problem, so we can also solve the ternary problem for the primes $x^2+y^2+1$.

\begin{theorem}
\label{theo_ternary} All large enough odd positive integers $n$ can be represented as $n=p+q+r$ with $p,q, r \in \mathscr{P}$.
\end{theorem}

We remark that Tolev \cite{tolev_ternarygoldbach} established an asymptotic formula for the weighted count of the representations of $n$ as $n=p+q+r$ with $p,q\in \mathscr{P}$ but $r$ a generic prime. The proof of Theorem \ref{theo_ternary} is very similar to that of Theorem \ref{theo_goldbach}, and is remarked on in Section \ref{Sec: transference}. \\

As a byproduct of the method for proving Theorem \ref{theo_goldbach}, we will obtain an analog of Roth's theorem for the set of primes of the form $x^2+y^2+1$, so that in particular the set $\mathscr{P}$ contains infinitely many three term arithmetic progressions. 

\begin{theorem} \label{theo_roth} 
Any subset of $\mathscr{P}^{*}=\{x^2+y^2+1:\,\, x,y\,\, \textnormal{coprime}\}\cap \mathbb{P}$ having a positive upper density with respect to $\mathscr{P}^{*}$ contains infinitely many non-trivial three term arithmetic progressions.
\end{theorem}

We will also conclude from the proof of Theorem \ref{theo_goldbach} that for any irrational $\xi$, there is some uniformity in the distribution of the fractional parts of the numbers $\xi p$ with $p\in \mathscr{P}$.

\begin{theorem}\label{theo_alphap} Let  $\xi$ be irrational and $\kappa\in \mathbb{R}$. Then there are infinitely many primes $p\in \mathscr{P}$ such that $\|\xi p+\kappa\|\leq p^{-\theta}$, where $\theta=\frac{1}{80}-\varepsilon=0.0125-\varepsilon$ and $\varepsilon>0$ is arbitrary. Here $\|\cdot\|$ stands for the distance to the nearest integer.
\end{theorem}

Theorems \ref{theo_roth} and \ref{theo_alphap} are proved in Sections  \ref{Sec: restriction} and \ref{Sec: fractional parts}, respectively. In Theorem \ref{theo_alphap}, we have not pursued maximizing the value of $\theta$, and the main message is that $\theta$ can be taken to be positive.\\

It should be remarked that the distribution of $\xi p \hspace{-0.1cm}\pmod 1$ has been studied also for some other subsets of the primes, such as for Chen primes \cite{matomaki-bombieri}, \cite{shi} and very recently for Gaussian primes \cite{baier} and Piatetski-Shapiro primes \cite{guo}. In the case of Chen primes the analog of Theorem \ref{theo_alphap} with $\theta>0$ was obtained in \cite{matomaki-bombieri} (and improved in \cite{shi} to $\theta=\frac{3}{200}=0.015$).\\ 

The proof of Theorem \ref{theo_goldbach} is based on a recent paper of Matom\"aki and Shao \cite{matomaki-shao}, where a transference type theorem for additive problems of Goldbach type was established, allowing one to deduce from certain desirable properties of a set $A$ the conclusion that $A+A+A$ contains all large enough integers. One should mention that a closely related transference principle for translation invariant additive problems was famously introduced by Green \cite{green-annals} and Green-Tao \cite{green-restriction}, \cite{green-tao} to find arithmetic progressions in the primes, their principle stating that a set $A$ with certain desirable properties contains infinitely many $3$-term arithmetic progressions (or $k$-term arithmetic progressions if one assumes stronger conditions). The hypotheses of the transference type result for Goldbach type equations \cite[Theorem 2.3]{matomaki-shao} resemble the ones of the transference principle for translation invariant equations \cite[Proposition 5.1]{green-restriction}, but include an additional assumption. An additional assumption is evidently needed, since for example the primes $p$ satisfying $\|\sqrt{2}p\|<\frac{1}{100}$ contain a lot of arithmetic progressions, but most odd integers are not the sum of three such primes.\\

The first property required from a set $A$ in the transference type result of \cite{matomaki-shao} is "well-distribution" in \textit{Bohr sets}, meaning that for $\xi,\kappa\in \mathbb{R}$ and $\eta>0$ the sets $\{n:\, \|\xi n+\kappa\|\leq \eta\}$  and their intersections contain a fair proportion of the elements of $A$. The second property, which is present in \cite{green-restriction} as well, is that $A$ is "Fourier bounded", in the sense that the Fourier transform $\widehat{1_A}$ is small in $\mathcal{\ell}^r$ norm for $r>2$. The last and simplest to check condition is that there should be a lower bound of the correct order of magnitude for the number of elements in $A$ up to $N$. In \cite{matomaki-shao}, the transference type result was applied to solve the ternary Goldbach problem with three Chen primes or with three primes $p$ such that $[p,p+C]$ contains at least two primes for some large constant $C$.\\

We employ a variant of the transference type result of \cite{matomaki-shao} in this paper, the conditions for the principle being nearly identical, but with the conclusion that $A+A$ contains almost all positive integers (in the sense that there are $o(N)$ integers $n\leq N$ not representable in this form). This modification is easy to implement, so the main part of our proof is devoted to verifying the conditions involved in the transference type result in the context of the set $\mathscr{P}$. The lower bound condition follows essentially from earlier work, so we are mostly concerned with proving two requirements.\\ 

The Fourier boundedness requirement follows from the restriction theory of the primes, in the form developed by Green and Tao in \cite{green-restriction}. However, the ''enveloping sieve'' $\beta(n)$ (which is a pseudorandom majorant of a subset of the primes and enjoys certain pleasant Fourier properties) has to be modified. It turns out that the necessary modification is available in a paper of Ramaré and Ruzsa \cite{ramare}, where the enveloping sieve was developed for purposes related to additive bases, and actually the results in that paper imply that $\mathscr{P}$ is an additive basis of finite (but large and unspecified) order.\\

Proving the well-distribution of the set $\mathscr{P}$ in Bohr sets requires more work and occupies the majority of this paper. We use a strategy similar to the one that was used in \cite{matomaki-shao} to deal with Chen's primes or with primes $p$ with $[p,p+C]$ containing two primes for some large constant $C$, but we must use a different sieve to detect primes of the form $x^2+y^2+1$. The sieve suitable for this purpose is a combination of the linear sieve and the semilinear sieve (also called the half-dimensional sieve), developed by Iwaniec in \cite{iwaniec-semilinear} and used by him in \cite{iwaniec-quadraticform} to prove that the number of primes in $\mathscr{P}$ up to $N$ is $\gg N(\log N)^{-\frac{3}{2}}$ (the infinitude of the primes in $\mathscr{P}$ was established earlier by Linnik \cite{linnik}  in 1960, using his dispersion method). An upper bound for $|\mathscr{P}\cap [1,N]|$ of the same order of magnitude follows from the Selberg sieve, so $\mathscr{P}$ is a sparse set of primes.\\

When it comes to the sieve theoretic part of the argument, we proceed along the lines of \cite{matomaki-m2+n2+1} and \cite{wu} that consider the problem of finding primes from $\mathscr{P}$ in short intervals. However, unlike in these works, one cannot apply the Bombieri-Vinogradov theorem for the prime counting function, but one has to resort to a Bombieri-Vinogradov type result for exponential sums $\sum_{n\leq N}\Lambda(n)e(\alpha n)$ over primes. Such average results for exponential sums appeared for instance in  \cite{tolev_bombieri}, \cite{mikawa-bombieri}, \cite{matomaki-bombieri}, but the level of distribution achieved in these works when the weight sequence is not well-factorable (in the sense defined in \cite[Chapter 12]{friedlander}) is $\frac{1}{3}-\varepsilon$, which is not good enough for our purposes. We derive a combinatorial factorization for the semilinear sieve weights and apply \cite[Lemma 8.4]{matomaki-shao} (closely related to the estimates in \cite{mikawa-bombieri}) on Bombieri-Vinogradov type averages for $\sum_{n\leq N}\Lambda(n)e(\alpha n)$ to increase the level of distribution sufficiently and hence obtain Theorem \ref{theo_goldbach}. In particular, the results of Sections \ref{Sec: Bombieri}, \ref{Sec: sieveweight} and \ref{Sec: hypotheses} imply the following Bombieri-Vinogradov type bound.

\begin{theorem}\label{theo_sievebombieri} Let $N\geq 1$ be large and $\varepsilon>0$, $C\geq 10$ fixed, and let $\lambda_d^{+,\textnormal{SEM}}$ and $\lambda_d^{-,\textnormal{SEM}}$ be the upper and lower bound semilinear sieve weights defined by restricting the M\"obius function $\mu(d)$ to the sets
\begin{align*}
\mathcal{D}^{+,\textnormal{SEM}}&=\{p_1\cdots p_r\leq N^{\rho_{+}}:\,\, z_{+}\geq p_1> \ldots > p_r,\,\,p_1\cdots p_{2k-2}p_{2k-1}^2\leq N^{\rho_{+}}\,\, \textnormal{for all}\,\, k\geq 1\},\\
\mathcal{D}^{-,\textnormal{SEM}}&=\{p_1\cdots p_r\leq N^{\rho_{-}}:\,\, z_{-}\geq p_1> \ldots > p_r,\,\,p_1\cdots p_{2k-1}p_{2k}^2\leq N^{\rho_{-}}\,\, \textnormal{for all}\,\, k\geq 1\}, 
\end{align*}
with the choices $\rho_{+}=\frac{2}{5}-\varepsilon$, $\rho_{-}=\frac{3}{7}-\varepsilon$, $z_{+}\leq N^{\frac{1}{2}}$ and $z_{-}\leq N^{\frac{1}{3}-\varepsilon}$.
Let $\alpha$ be a real number with $|\alpha-\frac{a}{q}|\leq \frac{1}{q^2}$ for some coprime integers $a$ and $q$ with $q\in [(\log N)^{1000C},N(\log N)^{-1000C}]$. Then for any integer $b\neq 0$ we have (choosing either $+$ or $-$ sign throughout)
 \begin{align*}
\sum_{\substack{d\leq N^{\rho_{\pm}}\\(d,b)=1}}\bigg|\lambda_{d}^{\pm,\textnormal{SEM}}\sum_{\substack{n\sim N\\n\equiv b \pmod{d}}}\Lambda(n)e(\alpha n)\bigg|\ll \frac{N}{(\log N)^{C}}.
\end{align*}
\end{theorem}

We remark that the arguments of this paper would easily generalize to primes of the form $x^2+y^2+a$, where $a\neq 0$ is any integer. We also note that since for all the primes of the form $x^2+y^2+1$ appearing in the rest of the paper the only possible common prime factors of $x$ and $y$ are $2$ and $3$, Theorem \ref{theo_goldbach} could be stated in the form that almost all even $n\not \equiv 5,8\pmod 9$ are representable as $n=p+q$ with $p$ and $q$ primes and neither $p-1$ nor $q-1$ having any prime factors greater than $3$ that are $\equiv -1\pmod{4}$. One should also mention that we did not get an asymptotic formula for the number of representations of $n$ as sums of two or three primes from $\mathscr{P}$ (unlike in the work of Tolev \cite{tolev_binarygoldbach}, \cite{tolev_ternarygoldbach} on related problems), nor did we show that the number of exceptional $n$ in Theorem \ref{theo_goldbach} is $\ll \frac{N}{(\log N)^{A}}$ instead of merely $o(N)$. We can nevertheless get a lower bound of $c n(\log n)^{-3}$ for the number of representations in Theorem \ref{theo_goldbach} for almost all $n$ for some small $c>0$, and this is the correct order of magnitude.

\subsection{Structure of the proofs}
We give a brief outline of the dependencies between different theorems and propositions. The proof of Theorem \ref{theo_goldbach} is deduced from the transference type theorem (Proposition \ref{prop_transference}) in Section \ref{Sec: 3}, provided that the two key conditions in the transference type theorem are satisfied. One condition is the well-distribution of the set $\mathscr{P}$ in Bohr sets and the other one is a Fourier uniformity result for $\mathscr{P}$ (Propositions \ref{prop_bohr} and \ref{prop_restriction}, respectively). The proof of Proposition \ref{prop_restriction} is presented in Section \ref{Sec: restriction}, and in Section \ref{Sec: 3} it is shown that Propositions \ref{prop_bohr} and \ref{prop_restriction} immediately imply Theorem \ref{theo_roth}.\\

 The largest part of the paper is then devoted to proving Proposition \ref{prop_bohr} using sieve theory.  The purpose of Section \ref{Sec: reductions} is to show that Proposition \ref{prop_bohr} follows from Proposition \ref{prop2}, which involves more notation but is easier to approach. In Section \ref{Sec: weighted}, a weighted sieve for finding primes of the form $x^2+y^2+1$ is presented, in the form of Theorem \ref{t2}. Section \ref{Sec: decomposition} constructs the weighted sequence $(\omega_n)$ to which Theorem \ref{t2} is applied, as well as sets up the circle method. Section \ref{Sec: hypotheses} is then devoted to proving Hypothesis \ref{h1} for $(\omega_n)$, since this hypothesis is the requirement for applying Theorem \ref{t2}. Section \ref{Sec: hypotheses}, which finishes the proofs of Theorems \ref{theo_goldbach} and \ref{theo_sievebombieri}, involves bounding Bombieri-Vinogradov sums related to either semilinear or linear sieve coefficients and  weighted by additive characters that lie either on minor or major arcs. The type I and II input required in Section \ref{Sec: hypotheses} comes from Section \ref{Sec: Bombieri}, while the required combinatorial input comes from Section \ref{Sec: sieveweight}. As Remark \ref{rmk1} tells, the only difference in the proofs of Theorems \ref{theo_ternary} and \ref{theo_goldbach} is the form of transference type result being used. Finally, when it comes to proving Theorem \ref{theo_alphap}, one needs the sections from Section \ref{Sec: weighted} onwards, the last of which, Section \ref{Sec: fractional parts}, is required only for this purpose. We also remark that none of the sections \ref{Sec: transference}, \ref{Sec: restriction}, \ref{Sec: reductions}, \ref{Sec: weighted}, \ref{Sec: Bombieri} and \ref{Sec: sieveweight} depend on each other.

\subsection{Notation}   

The symbols $j,k,\ell,m, n$ and $q$ always denote integers, and $p$ is a prime number. We denote by $e(\alpha)=e^{2\pi i  \alpha}$ the complex exponential, by $\text{Li}(x)=\int_{2}^{x}\frac{dt}{\log t}$ the logarithmic integral, and by $\pi(x;q,a)$ the number of primes up to $x$ in the residue class $a\hspace{-0.1cm} \pmod q$. We denote by $\|\cdot\|$ the distance to the nearest integer function, by $(\cdot, \cdot)$ the greatest common divisor and by $[\cdot,\cdot]$ the least common multiple. We denote by $\mathbb{Z}_q$ the set of integers $\pmod{q}$, sometimes interpreting functions defined on this set as $q$-periodic functions on $\mathbb{Z}$ and vice versa. The expression $m^{-1} \pmod{q}$ stands for the inverse of $m$ in $\mathbb{Z}_{q}$.\\

Starting from Section \ref{Sec: 3}, there are various symbols that have been reserved a specific meaning. The integer $\mathcal{C}$ is given by \eqref{eq13}, the function $s(n)$  by \eqref{eq60}, the set $\mathcal{S}$ by \eqref{eq47}, the integer $b$ by Definition \ref{def1}, the numbers $U, J$ and $W$ by  \eqref{eq30}, the set $\mathcal{Q}$ by \eqref{eq21}, the product $\mathfrak{S}(L)$ by Definition \ref{def3}, the function $g(\ell)$ by Definition \ref{def2}, and lastly the parameter $Q$ by Lemma \ref{le11}. When it comes to sieve theoretic notation, $\lambda_d$ are sieve weights and for a set $\mathcal{A}$ of integers and $\mathcal{P}$ of primes, $S(\mathcal{A},\mathcal{P},z)$ counts the elements of $\mathcal{A}$ that are coprime to all the primes in $\mathcal{P}\cap [2,z)$, with each integer $n$ weighted by $\omega_n\geq 0$, where $(\omega_n)$ will be clear from context. The arithmetic functions $\Lambda(n)$, $\mu(n)$ and $\varphi(n)$ are the von Mangoldt, M\"obius and Euler functions, as usual, and the functions $\tau(n)$ and $\nu(n)$ count the number of divisors and distinct prime factors of $n$, respectively.\\

The parameters $\varepsilon,\eta>0$ are always assumed to be small enough, but fixed. The variables $N$ and $x$ tend to infinity, and in Sections \ref{Sec: decomposition} and \ref{Sec: hypotheses}, $A,B$ and $C$ are large enough constants (say greater than $10^{10}$). The numbers $\mathcal{C}$, $W$ and $J$ are $\ll 1$, but may be large. The expression $1_{S}$ is the indicator function of a set $S$, so that $1_{S}(n)=1$ when $n\in S$ and $1_{S}(n)=0$ otherwise. We use the usual Landau and Vinogradov asymptotic notations $o(\cdot), O(\cdot), \ll, \gg$. When we write $n\sim X$ in a summation, we mean $X\leq n<2X$. By $n\asymp X$, in turn, we mean $X\ll n\ll X$.

\subsection{Acknowledgments}

The author is grateful to his supervisor Kaisa Matom\"aki for various useful comments and discussions. The author thanks the referee for careful reading of the paper and for useful comments. While working on this project, the author was funded by UTUGS Graduate School and project number 293876 of the Academy of Finland.

\section{A transference type result}\label{Sec: transference}

We need a transference type result for binary Goldbach type problems for proving Theorem \ref{theo_goldbach}.  We begin with some definitions.\\

Let $\Omega\subset \mathbb{Z}_N$ and $\eta\in (0,\frac{1}{2})$, and write
\begin{align*}
B(\Omega,\eta)=\left\{n\in \mathbb{Z}_N:\quad \left\|\frac{\xi n}{N}\right\|\leq \eta\quad \text{for all}\quad \xi \in \Omega\right\}    
\end{align*}
for the \textit{Bohr set} associated to these parameters. We will need a function $\chi=\chi_{\Omega,\eta}:\mathbb{Z}\to \mathbb{R}_{\geq 0}$ that is a smoothed version of the characteristic function of the Bohr set $B(\Omega,\eta)$. The exact construction of $\chi$ is not necessary, and we just list the properties of $\chi$ we use, found in \cite[Lemma 3.1]{matomaki-shao}. We have 
\begin{equation}\label{eq29}\begin{split}
&0\leq \chi(n)\ll_{|\Omega|}1,\hspace{2.4cm}  \chi(n)=\chi(-n)\,\, \text{and}\,\, \chi(n+N)=\chi(n),\\
&\chi(n)\geq 1\,\, \text{for}\,\,\ n\in B(\Omega,\eta),\quad \quad \chi(n)\leq \left(\frac{\eta^2}{8}\right)^{|\Omega|},\,\, \text{for}\,\,  n \not \in B(\Omega,2\eta)\\
&\frac{1}{N}\sum_{n\in \mathbb{Z}_N}\chi(n):=\|\chi\|_1\geq \left(\frac{\eta}{2}\right)^{|\Omega|}.
\end{split}
\end{equation}
Also from \cite{matomaki-shao}, we know that $\chi$ has \textit{Fourier complexity} $\mathcal{C}\ll_{|\Omega|,\eta} 1$, where the Fourier complexity is defined as the smallest integer $\mathcal{C}$ for which we have a Fourier representation
\begin{align}\label{eq13}
\chi(n)=\sum_{k=1}^{\mathcal{C}}c_k e(\alpha_k n),\,\, |c_k|\leq \mathcal{C}\,\, \text{and}\,\, \alpha_k\in \mathbb{R}/\mathbb{Z}.    
\end{align}

The formulation of the transference type result requires harmonic analysis, so we should state which normalization of the Fourier transform we use. For functions $f,g:\mathbb{Z}_N\to \mathbb{C}$ we define the Fourier transform and the convolution as
\begin{align*}
\hat{f}(\xi)&=\frac{1}{N}\sum_{n\in \mathbb{Z}_N}f(n)e\left(-\frac{\xi n}{N}\right)\quad \text{and}\quad f*g(n)=\frac{1}{N}\sum_{k\in \mathbb{Z}_N}f(k)g(n-k),
\end{align*}
so that Parseval's identity and the convolution formula of the Fourier transform take the forms
\begin{align*}
\sum_{n\in \mathbb{Z}_N}|f(n)|^2&=N\sum_{\xi \in \mathbb{Z}_N}|\hat{f}(\xi)|^2 \quad \text{and}\quad \widehat{f*g}(\xi)=\hat{f}(\xi)\hat{g}(\xi).
\end{align*}

\begin{proposition}\label{prop_transference} Let functions $f_1, f_2:\mathbb{Z}_N\to \mathbb{R}_{\geq 0}$ and parameters $K_0\geq 1$, $\delta>0$, $\varepsilon>0$ be given. Then there exist $\eta=\eta(K_0,\delta,\varepsilon)>0$ and $\Omega\subset\mathbb{Z}_N$, $|\Omega|\ll_{K_0,\delta,\varepsilon} 1$ with $1\in \Omega$ such that the following holds. Assume that, for a function $\chi=\chi_{\Omega,\eta}:\mathbb{Z}\to \mathbb{R}_{\geq 0}$ obeying \eqref{eq29}, we have\\
(i) $f_2*\chi(t)\geq \delta \|\chi\|_1$ for all $t\in (\frac{N}{3},\frac{2N}{3})$,\\
(ii) $\displaystyle \sum_{\frac{N}{3}<n<\frac{N}{2}}f_1(n)\geq \delta N$,\\
(iii) $\displaystyle \sum_{\xi \in \mathbb{Z}_N}|\widehat{f_j}(\xi)|^r\leq K_0$ for $j\in \{1,2\}$ and $r\in \{3,4\}$.\\
Then\\
(iv) $f_1*f_2(n)\geq \frac{\delta^2}{3}$ for all but $\leq \varepsilon N$ values of $n\in [0.9N,N]$.
\end{proposition}

\textbf{Proof.} This is inspired by and similar to \cite[Theorem 2.3]{matomaki-shao} of Matom\"aki and Shao. See also \cite[Proposition 5.1]{green-restriction}, where similar ideas were applied for Roth type problems. Take $\Omega=\{\xi \in \mathbb{Z}_N:\,\, |\widehat{f_1}(\xi)|\geq \varepsilon_0\}\cup\{1\}$, where $\varepsilon_0$ will be chosen small enough in terms of $\delta$, $\varepsilon$ and $K_0$. Condition (iii) tells that $|\Omega|\leq K_0\varepsilon_0^{-3}+1$. Let $\chi=\chi_{\Omega,\eta}:\mathbb{Z}\to \mathbb{R}_{\geq 0}$ be as in the proposition (so that $\chi$ fulfills \eqref{eq29}). We will later choose $\eta$ to be small enough in terms of $\delta$, $\varepsilon$ and $K_0$. Introduce the functions
\begin{align*}
g_2=\frac{1}{\|\chi\|_1}f_2*\chi\quad \text{and}\quad h_2=f_2-g_2.   
\end{align*}
We have
\begin{align*}
\widehat{g_2}= \frac{1}{\|\chi\|_1}\widehat{f_2}\widehat{\chi}\quad \text{and}\quad \widehat{h_2}=\widehat{f_2}\left(1-\frac{\widehat{\chi}}{\|\chi\|_1}\right),   
\end{align*}
so that in particular $|\widehat{h_2}(\xi)|\leq 2|\widehat{f_2}(\xi)|$.\\

Next we estimate from above and below the average  $\frac{1}{N}\sum_{n\in \mathbb{Z}_N}|f_1*h_2(n)|^2$, starting with the lower bound. Owing to conditions (i) and (ii), for $n\in [0.9N,N]$ we have
\begin{align}\label{eq33}
f_1*g_2(n)=\frac{1}{\|\chi\|_1}f_2*\chi*f_1(n)\geq \frac{\delta }{N}\sum_{\substack{n-\frac{2N}{3}<k<n-\frac{N}{3}\\ k\in \mathbb{Z}_N}}f_1(k)\geq \delta^2  
\end{align}
since $(\frac{N}{3},\frac{N}{2})\subset (n-\frac{2N}{3},n-\frac{N}{3})$ for $n\in [0.9N,N]$. Denoting $T=\{n\in [0.9N,N]:\,\, f_1*f_2(n)<\frac{\delta^2}{3}\}$ and using the simple inequality $|a-b|^2\geq \frac{a^2}{2}-b^2$ and \eqref{eq33}, we infer that
\begin{align}\begin{split}\label{eq31}
\frac{1}{N}\sum_{n\in \mathbb{Z}_N}|f_1*h_2(n)|^2&\geq \frac{1}{N}\sum_{n\in T}\left(\frac{1}{2}|f_1*g_2(n)|^2-|f_1*f_2(n)|^2\right)\\
&\geq \left(\frac{\delta^4}{2}-\left(\frac{\delta^2}{3}\right)^2\right)\frac{|T|}{N}\geq \frac{\delta^4 }{10}\frac{|T|}{N}.   
\end{split} 
\end{align}
When it comes to an upper bound, Parseval's identity gives 
\begin{align*}
\frac{1}{N}\sum_{n\in \mathbb{Z}_N}|f_1*h_2(n)|^2&=\sum_{\xi\in \mathbb{Z}_N}|\widehat{f_1*h_2}(\xi)|^2\\
&=\sum_{\xi\in \mathbb{Z}_N}|\widehat{f_1}(\xi) \widehat{h_2}(\xi)|^2\\
&\leq \varepsilon_0^{\frac{1}{2}}\sum_{\xi\not \in \Omega}|\widehat{f_1}(\xi)|^{\frac{3}{2}} |\widehat{h_2}(\xi)|^2+\sum_{\xi\in \Omega}|\widehat{f_1}(\xi)|^2 |\widehat{h_2}(\xi)|^2.
\end{align*}
Here the first sum can be bounded with the Cauchy-Schwarz inequality and (iii), implying
\begin{align*}
\varepsilon_0^{\frac{1}{2}}\sum_{\xi\not \in \Omega}|\widehat{f_1}(\xi)|^{\frac{3}{2}} |\widehat{h_2}(\xi)|^2 \leq \varepsilon_0^{\frac{1}{2}}\left(\sum_{\xi \in \mathbb{Z}_N}|\widehat{f_1}(\xi)|^3\right)^{\frac{1}{2}}\left(\sum_{\xi \in \mathbb{Z}_N}|\widehat{h_2}(\xi)|^4\right)^{\frac{1}{2}}\leq 8\varepsilon_0^{\frac{1}{2}}K_0.    
\end{align*}
The sum over $\xi \in\Omega$ in turn can be bounded by using the fact that
\begin{align*}
\left|1-\frac{\widehat{\chi}(\xi)}{\|\chi\|_1}\right|\leq 30\eta \quad \text{for every}\quad  \xi \in \Omega,    
\end{align*}
the proof of which is contained in the proof of Theorem 2.3 in \cite[Section 4]{matomaki-shao}. After this, we may again use the Cauchy-Schwarz inequality and (iii) to get
\begin{align*}
\sum_{\xi\in \Omega}|\widehat{f_1}(\xi)|^2 |\widehat{h_2}(\xi)|^2&\leq (30\eta)^{2}\sum_{\xi\in \Omega}|\widehat{f_1}(\xi)|^2 |\widehat{f_2}(\xi)|^{2}\\
&\leq 1000\eta^{2}K_0.
\end{align*}
At this stage, we fix the choices $\varepsilon_0=\eta=\frac{\delta^8\varepsilon^2}{10^4K_0^{2}}$, so that 
\begin{align}\label{eq32}
\frac{1}{N}\sum_{n\in \mathbb{Z}_N}|f_1*h_2(n)|^2\leq 8\varepsilon_0^{\frac{1}{2}}K_0+1000\eta^{2}K_0\leq \frac{1}{10}\delta^4 \varepsilon.    
\end{align}
Combining \eqref{eq31} and \eqref{eq32} above, we discover that $|T|\leq 10\delta^{-4}\cdot \frac{1}{10}\delta^4\varepsilon N=\varepsilon N$, which concludes the proof.\qedd

\section{Deducing Theorem \ref{theo_goldbach} from the transference type result}\label{Sec: 3}

We will apply the transference type result (Proposition \ref{prop_transference}) to prove Theorem \ref{theo_goldbach}. This deduction is done in this section assuming the conditions (i)-(iii) of the transference type result, and the rest of the paper is focused on verifying these conditions. Naturally, the functions $f_1$ and $f_2$ in the transference type result are taken to be the characteristic functions of the primes of the form $x^2+y^2+1$ (restricted to a residue class), normalized in such a way that they have mean comparable to $1$. First, we introduce some notation.\\

Define the function
\begin{align}\label{eq60}
s(n)=\prod_{\substack{p\mid n\\ p\equiv -1 \hspace{-0.1cm} \pmod 4\\p\neq 3}}p,    
\end{align}
which excludes from the prime factorization of $n$ the primes $2$, $3$ and those primes that are $\equiv 1 \hspace{-0.1cm} \pmod 4$. Denote
\begin{align}\label{eq47}
\mathcal{S}=\{a^2+b^2:\quad a,b\in \mathbb{Z},\quad (a,b)\mid 6^{\infty}\}.
\end{align}
 We also define a property that we require from the linear functions we work with in what follows.

\begin{definition}\label{def1}
We say that a linear polynomial $L$ with integer coefficients is \textit{amenable} if $L(n)=Kn+b$ for some integers $K\geq 1$ and $b$, and\\
(i) $6^3\mid K$,\\
(ii) $(b,K)=(b-1,s(K))=1$,\\
(iii) $b-1=2^j 3^{2t}(4h+1)$ for some $h\in \mathbb{Z}$, $3\nmid 4h+1$ and $j,t\geq 0$ with $2^{j+2}3^{2t+1}\mid K$.
\end{definition}

What these conditions imply is that there are no local obstructions (modulo divisors of $K$) to $L(n)$ being prime and $L(n)-1$ belonging to $\mathcal{S}$ (in particular, $L(n)-1$ crucially has an even number of prime factors $p\equiv -1\hspace{-0.1cm} \pmod 4$ with multiplicities by (iii)). We note that it is essential that $b-1$ is allowed to be divisible by a power of $3$. Indeed, if $L_i(n)=Kn+b_i$ are two amenable linear functions with $3\mid K$ and $3\nmid b_1-1$, $3\nmid b_2-1$, then $L_1(m)+L_2(n)$ can only represent numbers that are $\equiv 1 \hspace{-0.1cm}\mod 3$. We also note that in our application we must allow $K$ to be divisible by arbitrarily high powers of $2$. This is due to the fact that if $L_i(n)=2^sn+b_i$ are amenable, then $L_i(n)-1\equiv 2^{a_i}\hspace{-0.1cm} \pmod{2^{a_i+2}}$ for some integers $0\leq a_i\leq s-2$, which implies that $L_1(m)+L_2(n)$ is never $\equiv 2 \pmod{2^s}$.\\

The majority of this paper is devoted to proving for functions $f_i$ related to the characteristic function of $\mathscr{P}$ the following versions of the conditions (i) and (iii) of the transference type result. Throughout the rest of the paper, we denote
\begin{align}\label{eq30}\begin{split}
U&=2^{J}\cdot 3^3 \quad \text{with} \quad 5\leq J\ll 1,\\
W&=U\cdot \prod_{5\leq p\leq w} p\quad \text{with}\quad 10^{{10}^{10}}\leq w\ll 1.    \end{split}
\end{align}

\begin{proposition}\label{prop_bohr} Let $\chi:\mathbb{Z}\to \mathbb{R}_{\geq 0}$ have Fourier complexity $\mathcal{C}\ll 1$. Let $W$ be as in \eqref{eq30} with $w\geq \mathcal{ C}^{20}$, and suppose that the linear function $Wn+b$ is amenable. For an integer $N\geq 1$, set 
\begin{align}\label{eq76}
f(n)=(\log N)^{\frac{3}{2}}\left(\frac{\varphi(W)}{W}\right)^{\frac{3}{2}}1_{Wn+b\in \mathbb{P},\,\,Wn+b-1\in \mathcal{S}} \quad \textnormal{for}\quad n\in \left(\frac{N}{3},\frac{2N}{3}\right),  
\end{align}
and $f(n)=0$ for other values of $n\in [0,N)$. Then for $N\geq N_0(w,\mathcal{C})$ we have
\begin{align*}
\sum_{n\sim  \frac{N}{3}}f(n)\chi(t-n)&\geq \delta_0 \bigg(\sum_{n\sim \frac{N}{3}}\chi(t-n)-\frac{CN}{w^{\frac{1}{3}}}\bigg)
\end{align*}
for $t\in (\frac{N}{3},\frac{2N}{3})$ and some absolute constants $\delta_0>0, C>0$.
\end{proposition}

\begin{proposition}\label{prop_restriction} Suppose that the linear function $Wn+b$ is amenable with $W$ as in \eqref{eq30}. Let $N\geq 1$ be an integer and $g:\mathbb{Z}_N\to \mathbb{R}_{\geq 0}$ with $0\leq g(n)\leq f(n)$ for $n\in [0,N)$ and $f$ as in \eqref{eq76}. Then for all $r>2$,
\begin{align*}
\sum_{\xi \in \mathbb{Z}_N}|\widehat{g}(\xi)|^r\leq K_r 
\end{align*}
for some positive constant $K_r$ depending only on $r$.
\end{proposition}

In this section, we show that Propositions \ref{prop_bohr} and \ref{prop_restriction} indeed imply Theorem \ref{theo_goldbach}. First we prove some lemmas about local representations of integers modulo powers of $2$ and $3$.

\begin{lemma} \label{le12} Let $J\geq 5$ and $n\not \equiv 0\hspace{-0.1cm} \pmod{2^{J-1}}$ be integers. Then we may write $n=a+b$ for some integers $a$ and $b$ with $a\equiv 2^{i} \hspace{-0.1cm} \pmod{2^{i+2}}$ and $b\equiv 2^j \hspace{-0.1cm} \pmod{2^{j+2}}$ for some integers $0\leq i,j\leq J-3$. 
\end{lemma}

\textbf{Proof.} Since $2^{J-1}\nmid n$, we may write $n=2^{g} s$ where $0\leq g\leq J-5$ and $s\not \equiv 0\pmod{16}$. It is easy to check that every such $s$ may be written as $s=a'+b'$ with  $a'\equiv 2^{i} \hspace{-0.1cm} \pmod{2^{i+2}}$, $b'\equiv 2^j \hspace{-0.1cm} \pmod{2^{j+2}}$ for some $0\leq i,j\leq 3$. Then $n=a+b$ with $a=2^{g}a'$, $b=2^g b'$ is a representation of the desired form. \qedd

\begin{lemma}\label{le13} Let $m'$ be any integer such that $m'\not \equiv 3,6 \hspace{-0.1cm} \pmod 9$. Then there exist integers $x_1$, $x_2$, $x_3$ and $x_4$ such that
\begin{align*}
&m'\equiv x_1^2+x_2^2+x_3^2+x_4^2 \hspace{-0.1cm} \pmod{3^{3}}\\ 
&x_1^2+x_2^2,\quad x_3^2+x_4^2\not \equiv 1\hspace{-0.1cm} \pmod 3\\
&x_1^2+x_2^2,\quad x_3^2+x_4^2\not \equiv 0\hspace{-0.1cm} \pmod{3^3}
\end{align*}
\end{lemma}

\textbf{Proof.} One easily sees that $x^2 +y^2 \hspace{-0.1cm} \pmod{27}$ attains all residue classes except those that are $\equiv 3 \hspace{-0.1cm} \pmod 9$ or $\equiv 6 \hspace{-0.1cm} \pmod 9$ as $x$ and $y$ vary. Now the lemma only states that every $m'\not \equiv 3,6 \hspace{-0.1cm} \pmod 9$ is the sum of two numbers, each of which is $0,2,5$ or $8 \hspace{-0.1cm} \pmod 9$ and neither of which is $0\hspace{-0.1cm} \pmod{27}$. This can quickly be verified by hand. \qedd\\

\textbf{Proof of Theorem \ref{theo_goldbach} assuming Propositions \ref{prop_bohr} and \ref{prop_restriction}.} Given any small $\varepsilon>0$, we must show that once $N$ is large enough, the interval $[0.9N,N]$ contains at most $\varepsilon N$ integers $m\equiv 0\hspace{-0.1cm} \pmod 2$, $m\not \equiv 5,8 \hspace{-0.1cm} \pmod{9}$ that cannot be written as $m=p+q$ with $p$ and $q$ primes of the form $x^2+y^2+1$.\\

Let $U$ and $W$ be given by \eqref{eq30} with $J=\lfloor \frac{10}{\varepsilon}\rfloor$ and $w \ll 1$ large enough. We start by showing that for any $m\in [0.9N,N]$, $m\equiv 0\hspace{-0.1cm} \pmod 2$, $m\not \equiv 5,8\hspace{-0.1cm} \pmod 9$, $m\not \equiv 2\hspace{-0.1cm} \pmod{2^{J}}$,  we may find integers $0\leq B_1, B_2\leq W-1$ such that $m= B_1+B_2$ and the linear functions $Wn+B_1$ and $Wn+B_2$ are amenable. The integers $m\equiv 2\hspace{-0.1cm} \pmod{2^{J}}$ can be disposed of since there are $\leq \frac{\varepsilon^2}{10} N$ such integers up to $N$.\\

To see that $B_1$ and $B_2$ exist, write $m=2m'+2$, so that $m'\not \equiv 3,6 \hspace{-0.1cm} \pmod 9$. Then $2^{J-1}\nmid m'$, so using Lemma \ref{le12} we may write $m'\equiv a_1+a_2\hspace{-0.1cm} \pmod{2^{J}}$ with $a_1\equiv 2^{i}\hspace{-0.1cm} \pmod{2^{i+2}}$, $a_2\equiv 2^{j}\hspace{-0.1cm} \pmod{2^{j+2}}$ for some $0\leq i,j\leq J-3$. Moreover, using Lemma \ref{le13}, we may write $m'\equiv a_1'+a_2' \hspace{-0.1cm}\pmod{3^3}$ with $a_1'$ and $a_2'$ numbers such that $3^3\nmid a_1'$, $3^3\nmid a_2'$, $2a_1'+1,2a_2+1'\not \equiv 0 \pmod 3$, and the largest powers of $3$ dividing $a_1'$ and $a_2'$ have even exponents (take $a_1'=x_1^2+x_2^2$ and $a_2'=x_3^2+x_4^2$ in that lemma and notice that the largest power of $3$ dividing $x^2+y^2$ has an even exponent).\\

Now pick numbers $b_p$ for $5\leq p\leq w$ such that $b_p\not \equiv 0,1,m,m-1\hspace{-0.1cm} \pmod{p}$. By the Chinese remainder theorem, we can find an integer $B$ such that $B\equiv 2a_1+1\hspace{-0.1cm} \pmod{2^{J}}$, $B\equiv 2a_1'+1\hspace{-0.1cm} \pmod{3^3}$, and $B\equiv b_p\hspace{-0.1cm} \pmod{p}$ for all $5\leq p<w$. Therefore, we have found some integers $B_1:=B$ and $B_2:=m-B$ such that $m=B_1+B_2$ $p\nmid B_i$, $p\nmid B_i-1$ for $5\leq p<w$, and $B_1-1$ and $B_2-1$ satisfy condition (iii) in the definition of amenability.\\

Therefore, we have a representation of any $m$ of the form above as 
\begin{align*}
m\equiv B_1(m)+B_2(m)\hspace{-0.1cm} \pmod W    
\end{align*}
with $Wn+B_1(m)$, $Wn+B_2(m)$ amenable linear functions and $0\leq B_i(m)\leq W-1$ (we use the notation $B_i(m)$ to emphasize that the $B_i$ depend on $m \hspace{-0.1cm}\pmod W$). For each $0\leq a\leq W-1$ we denote
\begin{align*}
\mathcal{B}_{a}=\{m\in [0.9N,N]:\,\, m\equiv a \pmod W\}.    
\end{align*}
We will show that each $\mathcal{B}_{a}$ with $a\equiv 0\hspace{-0.1 cm} \pmod{2}$, $a\not \equiv 5,8\hspace{-0.1cm} \pmod 9$, $a\not \equiv 2\pmod{2^{J}}$ contains at most $\varepsilon \frac{N}{2W}$ values of $m\in [0.9N,N]$,  that are not of the form $p+q$ with $p$ and $q$ primes of the form $x^2+y^2+1$, and afterwards we sum this result over $a$.\\

If $a$ satisfies the congruence conditions above, the polynomials $Wn+B_1(a)$ and $Wn+B_2(a)$ are amenable linear polynomials. Set $M'=\lfloor\frac{N}{W}\rfloor$, and for $\ell\in \{1,2\}$ set
\begin{align*}
f_{\ell}(n)=(\log N)^{\frac{3}{2}}\left(\frac{\varphi(W)}{W}\right)^{\frac{3}{2}}1_{Wn+B_{\ell}(a)\in \mathbb{P},\,\,Wn+B_{\ell}(a)-1\in \mathcal{S}}\quad \text{for}\quad n\in \left(\frac{M'}{3},\frac{2M'}{3}\right),    
\end{align*}
with $\mathcal{S}$ as in \eqref{eq47} and let $f_{\ell}(n)=0$ for $n\in [0,M')\setminus (\frac{M'}{3},\frac{2M'}{3})$.\\ 

Concerning condition (ii) of the transference type result, applying Proposition \ref{prop_bohr} to the function $\chi\equiv 1$, we see that
\begin{align*}
\sum_{\frac{M'}{3}<n<\frac{2M'}{3}}f_1(n)\geq \frac{\delta_0}{10} M',    
\end{align*}
but we evidently get the same outcome with summation over $\frac{M'}{3}<n<\frac{M'}{2}$ (since one could clearly replace $n\sim \frac{N}{3}$ with $\frac{N}{3}<n<\frac{N}{2}$ in Proposition \ref{prop_bohr}). This takes care of condition (ii).\\

Next, by Proposition \ref{prop_restriction}, 
\begin{align*}
\sum_{\xi \in \mathbb{Z}_{M'}}|\widehat{f_{\ell}}(\xi)|^{r}\leq K_0 \end{align*}
for some absolute constant $K_0$ when $r\in \{3,4\}$, so also condition (iii) holds.

Let then $\chi=\chi_{\Omega,\eta}:\mathbb{Z}_{M'}\to \mathbb{R}_{\geq 0}$ be as in Proposition \ref{prop_transference} (with $\chi$ depending on $K_0$ and $\delta_0$ that appeared above), where $\Omega\subset \mathbb{Z}_{M'}$ satisfies $1\in \Omega$, $|\Omega|\ll_{\varepsilon} 1, and 1\ll_{\varepsilon} \eta\leq 0.05$. According to \eqref{eq29}, $\chi$ is symmetric around the origin and
\begin{align*}
\sum_{\substack{n\in [-\frac{M'}{2},\frac{M'}{2}]\\|n|\geq 0.1M'}}\chi(n)\leq \left(\frac{\eta^2}{8}\right)^{|\Omega|}M'\leq \eta\left(\frac{\eta}{2}\right)^{|\Omega|}M'\leq 0.05\|\chi_1\|M'.
\end{align*}
Keeping this in mind and using Proposition \ref{prop_bohr}, for $t\in (\frac{M'}{3},\frac{2M'}{3})$ we obtain
\begin{align*}
\sum_{n\sim \frac{M'}{3}}f_2(n)\chi(t-n)&\geq \delta_0 \bigg(\sum_{n\sim \frac{M'}{3}}\chi(t-n)-\frac{CM'}{w^{\frac{1}{3}}}\bigg)\\
&\geq \frac{\delta_0}{10}\bigg(\sum_{n\in \mathbb{Z}_{M'}}\chi(t-n)-\frac{CM'}{w^{\frac{1}{3}}}\bigg)\\
&\geq \frac{\delta_0}{20} M'\|\chi\|_1
\end{align*}
for $w$ large enough, the final step coming from \eqref{eq29}, since
\begin{align*}
\|\chi\|_{1}\geq \left(\frac{\eta}{2}\right)^{|\Omega|}\geq \frac{1}{w^{0.1}}    
\end{align*}
for $w$ large enough. This means that condition (i) of the transference type result holds with $\delta=\frac{\delta_0}{20}$.\\

 From the transference type result (Proposition \ref{prop_transference}), we conclude that $f_1*f_2(n)>0$ for all $n\in [0.9M',M']$, $n\not \in T_{a}$  where $T_{a}$ is some set of integers with $|T_{a}|\leq \frac{\varepsilon}{2} M'=\varepsilon\frac{N}{2W}$. This leads to $n\equiv n_1+n_2\hspace{-0.1cm} \pmod{M'}$ with 
\begin{align}\label{eq46}\begin{split}
&Wn_i+B_{i}(a)\in \mathbb{P}, \quad Wn_i+B_{i}(a)-1\in \mathcal{S}   
\end{split}
\end{align} 
for $n\in [0.9M',M']$, $n\not \in T_{a}$. Since $n_1,n_2\in (\frac{M'}{3},\frac{2M'}{3})$, we can actually say that $n=n_1+n_2$. What we showed at the beginning of the proof is that any $m\in \mathcal{B}_{a}$, $m\in [0.9N+2W,N]$ with $m\equiv 0\hspace{-0.1cm} \pmod 2$, $m\not\equiv 5,8\hspace{-0.1cm} \pmod 9$ and $m\not \equiv 2 \hspace{-0.1cm} \pmod{2^{J}}$ can be written as $m=Wn+B_1(a)+B_2(a)$ with $n\in[0.9M',M']$ and $Wn+B_1(a)$ and $Wn+B_2(a)$ amenable (the interval $[0.9N,0.9N+2W]$ contains $\leq \frac{\varepsilon^2}{10}N$ numbers and can hence be ignored). Then 
\begin{align*}
m=(Wn_1+B_1(a))+(Wn_2+B_2(a))
\end{align*}
for some $n_1$ and $n_2$ satisfying \eqref{eq46} whenever $m\in \mathcal{B}_{a}\setminus T'_{a}$,  $m\in [0.9N+2W,N]$, $m\equiv 0\hspace{-0.1cm} \pmod 2$, $m\not \equiv 5,8\hspace{-0.1cm} \pmod 9$ and $m\not \equiv 2\hspace{-0.1cm} \pmod{2^{J}}$, where $T'_{a}=\{a+W\tau:\,\, \tau\in T_a\}$ satisfies $|T'_{a}|\leq \varepsilon \frac{N}{2W}$. Since 
\begin{align*}
\sum_{\substack{0\leq a\leq W-1\\a\equiv 0\hspace{-0.1cm}\pmod 2\\a\not \equiv 5, 8 \hspace{-0.1cm}\pmod 9\\a\not \equiv 2 \hspace{-0.1cm}\pmod{2^{J}}}}|T_{a}|\leq W\cdot \varepsilon\frac{N}{2W}=\frac{\varepsilon}{2} N,    
\end{align*}
we conclude that all but $\leq (\frac{\varepsilon}{2}+\varepsilon^2) N\leq \varepsilon N$ even integers $m\in [0.9N,N]$ satisfying $m\not\equiv 5,8\hspace{-0.1cm} \pmod{9}$ can be written as $m=p+q$ with $p,q$ primes of the form $x^2+y^2+1$.\qedd

\begin{remark}\label{rmk1}
The proof of the ternary result, Theorem \ref{theo_ternary}, goes along very similar lines. One would replace Proposition \ref{prop_transference} with the analogous ternary transference type result, namely \cite[Theorem 2.3]{matomaki-shao}. The premises in both transference type results are essentially the same (except that \cite[Theorem 2.3]{matomaki-shao} has one additional function $f_3$), and therefore the differences in the proofs can only arise when showing that the transference type theorem implies the additive result. In fact, these proofs are also very similar, and one would simply replace Lemma \ref{le12} with a version where we want to represent an arbitrary integer $n$ as a sum of three numbers of the form $2^{i}\hspace{-0.1cm} \pmod{2^{i+2}}$, and one would replace Lemma \ref{le13} with a version where there is no restriction on $m'$ and there are six variables $x_i$ (and one would define $f_3$ analogously to $f_1$ and $f_2$).
\end{remark}

\section[Restriction theory for primes]{Restriction theory for primes of the form $x^2+y^2+1$}\label{Sec: restriction}

The objective of the current section  is proving Proposition \ref{prop_restriction}, after which proving Theorem \ref{theo_goldbach} has been reduced to demonstrating Proposition \ref{prop_bohr}. As a byproduct of the arguments, we will obtain Theorem \ref{theo_roth}. The proof of Proposition \ref{prop_restriction} is based on the Green-Tao approach \cite{green-restriction} that offers a way to estimate the Fourier norms of prime-related functions and therefore to detect translation invariant constellations within the primes.  The Green-Tao approach is based on proving a restriction theorem for the Fourier transform from $\ell^r(\mathbb{Z}_N)$ to $\ell^2(\mathbb{Z}_N)$ weighted by a certain "enveloping sieve" that acts as a pseudorandom majorant for the characteristic function of the primes of the desired form. Therefore, we start by asserting that there is a suitable enveloping sieve $\beta(\cdot)$ for the primes of the form $x^2+y^2+1$.  

\begin{proposition}\label{prop3}  Let $W$ and $w$ be as in \eqref{eq30}, and suppose that $B$ is an integer for which $Wn+B$ is an amenable linear function. Then, for any large $N$, there exists a function $\beta:\mathbb{N}\to \mathbb{R}_{\geq 0}$ with the following properties (for some absolute constants $\kappa_1,\kappa_2>0$):\\
(i) $\beta(n)\geq\kappa_1 (\log N)^{\frac{3}{2}}(\log w)^{-\frac{3}{2}}$ for $n\sim \frac{N}{3}$ when $Wn+B\in \mathbb{P}\cap (\mathcal{S}+1)$,\\
(ii) $\sum_{n\leq N}\beta(n)\leq\kappa_2 N$,\\
(iii) For every fixed $\varepsilon>0$, we have $\beta(n)\ll N^{\varepsilon}$,\\
(iv) We may write, for $z=N^{0.1}$,
\begin{align}\label{eq81}
\beta(n)=\sum_{q\leq z^2}\sum_{a\in \mathbb{Z}_q^{\times}} v\left(\frac{a}{q}\right)e\left(-\frac{an}{q}\right),   
\end{align}
where $v\left(\frac{a}{q}\right)\ll q^{\varepsilon-1}$ (and $\mathbb{Z}_q^{\times}$ is the set of primitive residue classes $\hspace{-0.1cm} \pmod q$),\\
(v) We have $v(1)=1$ and $v\left(\frac{a}{q}\right)=0$ in \eqref{eq81} whenever $q$ is not square-free or $q\mid W,q\neq 1$. 
\end{proposition}

The message of the previous proposition, which we will soon prove, is that $\beta(\cdot)$ is an upper bound for the normalized characteristic function of the primes $x^2+y^2+1$ in a residue class, $\beta(\cdot)$ has average comparable to $1$, and $\beta(\cdot)$ has a Fourier expansion with small coefficients. The above result implies the following restriction theorem, which is identical to \cite[Proposition 4.2]{green-restriction}, except that $\beta(\cdot)$ has a different definition. 

\begin{proposition}\label{prop_extension} Let $\beta:\mathbb{N}\to \mathbb{R}_{\geq 0}$ be as in Proposition \ref{prop3}. Let $N\geq 1$ be large, and let $(a_n)_{n\leq N}$ be any sequence of complex numbers. Given a real number $r>2$, for some $C_r>0$ we have
\begin{align*}
\left(\sum_{\xi\in \mathbb{Z}_N}\left|\frac{1}{N}\sum_{n\leq N}a_n\beta(n)e\left(\frac{-\xi n}{N}\right)\right|^r\right)^{\frac{1}{r}}\leq C_r\left(\frac{1}{N}\sum_{n\leq N}|a_n|^2\beta(n)\right)^{\frac{1}{2}}.    
\end{align*}
\end{proposition}

\textbf{Proof of Proposition \ref{prop_extension} assuming Proposition \ref{prop3}}: Our function $\beta(\cdot)$ fulfills the same axioms as in the paper of Green-Tao (except the pointwise lower bound, which is not used for the proof of  \cite[Proposition 4.2]{green-restriction}). Therefore, the proof of  \cite[Proposition 4.2]{green-restriction} goes through in this setting.\qedd\\

At this point, we show that Proposition \ref{prop_extension} easily implies Proposition \ref{prop_restriction}, which corresponds to condition (iii) in the transference type result.\\

\textbf{Proof of Proposition \ref{prop_restriction} assuming Proposition \ref{prop3}}: We already know that if Proposition \ref{prop3} is true, so is Proposition \ref{prop_extension}. We choose $a_n=\frac{g(n)}{\beta(n)}$ whenever $\beta(n)\neq 0$, and $a_n=0$ otherwise. Since $0\leq g(n)\leq f(n)\leq \kappa_1^{-1} \beta(n)$ in the notation of Proposition \ref{prop_restriction}, from Proposition \ref{prop_extension} we immediately derive
\begin{align*}
\left(\sum_{\xi \in \mathbb{Z}_N}|\widehat{g}(\xi)|^r\right)^{\frac{1}{r}}&\leq C_r \left(\frac{1}{N}\sum_{\substack{n\leq N\\ \beta(n)\neq 0}}\frac{g(n)^2}{\beta(n)}\right)^{\frac{1}{2}}\leq C_r \left(\frac{\kappa_1^{-2}}{N}\sum_{n\leq N}\beta(n)\right)^{\frac{1}{2}}\leq C_r\kappa_1^{-1}\kappa_2^{\frac{1}{2}}
\end{align*}
by part (ii) of Proposition \ref{prop3}.\qedd\\

What remains to be shown is that the enveloping sieve promised by Proposition \ref{prop3} exists. This is based on an argument of Ramaré and Ruzsa \cite{ramare} (which incidentally developed the enveloping sieve for purposes unrelated to restriction theory). The enveloping sieve $\beta(n)$ turns out to be a normalized Selberg sieve corresponding to sifting primes of the form $p=x^2+y^2+1$, $p\equiv B\pmod W$.\\

\textbf{Proof of Proposition \ref{prop3}}: We first introduce some notation. For a prime $p$, let $\mathcal{A}_p\subset \mathbb{Z}_p$ denote the residue classes $\hspace{-0.1cm} \pmod p$ that are sifted away when looking for primes of the form $x^2+y^2+1\equiv B\hspace{-0.1cm} \pmod W$. In other words,
\begin{align*}
\mathcal{A}_p=\begin{cases}\emptyset\,\, \text{for}\,\, p\leq w,\\
\{0\}\,\, \text{for}\,\,p\equiv 1 \hspace{-0.1cm} \pmod 4, \,\,p>w\\
\{0,1\}\,\, \text{for}\,\,p\equiv -1 \hspace{-0.1cm} \pmod 4,\,\,p>w.
\end{cases}
\end{align*}
Further, for square-free $d$ let 
\begin{align*}
\mathcal{A}_d=\bigcap_{p\mid d}\mathcal{A}_p,
\end{align*} 
where $\mathcal{A}_d$ is interpreted as a subset of $\mathbb{Z}_d$. Set also $\mathcal{A}_1=\mathbb{Z}_1$ and $\mathcal{A}_d=\emptyset$ when $d$ is not square-free. For $d\geq 2$, we have $|\mathcal{A}_d|=\omega(d)$, where $\omega(\cdot)$ is a multiplicative function supported on the square-free integers and having the values
\begin{align*}
\omega(p)=\begin{cases}0\,\, \text{for}\,\, p\leq w,\\
1\,\, \text{for}\,\,p\equiv 1 \hspace{-0.1cm} \pmod 4,\,\,p>w,\\
2 \,\, \text{for}\,\, p\equiv -1\hspace{-0.1cm} \pmod 4,\,\,p>w.
\end{cases}
\end{align*}

For later use, we also define
\begin{align}\label{eq80}
\mathcal{ K}_1=\mathbb{Z}_1,\quad \mathcal{K}_p=\mathbb{Z}_p\setminus \mathcal{A}_p,\quad \mathcal{K}_d=\bigcap_{p\mid d}\mathcal{K}_p\quad \text{for}\quad \mu(d)^2=1
\end{align}
and let  $\mathcal{K}_d=\mathbb{Z}_d$ for $\mu(d)=0$.\\

Let the Selberg sieve coefficients $\rho_d$ (not the same as sieve weights) be given by
\begin{align*}
&\rho_d=\mu(d)\frac{G_d(z)}{G_1(z)},\quad \text{where}\quad z=N^{0.1},\quad G_d(z)=\sum_{\substack{\delta\leq z\\ [d,\delta]\leq z}}h(\delta),\\
&h(\delta)=\prod_{p\mid \delta}h(p)\quad \text{and}\quad h(p)=\frac{\omega(p)}{p-\omega(p)}.
\end{align*}
The above notations are otherwise the same as in \cite[Section 4]{ramare}, except that $\lambda_d$ there has been replaced with $\rho_d$ and $\mathcal{L}_d$ with $\mathcal{A}_d$. We define
\begin{align}\label{eq78}
\beta(n)=G_1(z)\bigg(\sum_{\substack{d\mid P(z)\\ Wn+B\in \mathcal{A}_d}}\rho_d\bigg)^2,    
\end{align}
where 
\begin{align*}
P(z)=\prod_{w<p<z}p.    
\end{align*}
In \cite{ramare} the factor $G_1(z)$ does not appear in their definition of $\beta(n)$, but this is just a normalization constant. In \eqref{eq78} the condition $m\in \mathcal{A}_d$ means $m\hspace{-0.1cm} \pmod d \in \mathcal{A}_d$. Now we can check parts (i)-(v) of Proposition \ref{prop3}.\\

For part (i), first observe that if $Wn+B=x^2+y^2+1\in \mathbb{P}\cap (\mathcal{S}+1)$ with $n\sim \frac{N}{3}$, then $x^2+y^2+1\not \equiv 0\hspace{-0.1cm} \pmod p$ for $w< p<z=N^{0.1}$ and $x^2+y^2\not \equiv 0 \pmod p$ for $p\equiv -1\hspace{-0.1cm} \pmod 4$, $w< p<z$, since $(x,y)\mid 6^{J}$. This means that if $Wn+B=x^2+y^2+1\in \mathbb{P}\cap (\mathcal{S}+1)$ with $n\sim \frac{N}{3}$, then $\beta(n)=G_1(z)$. Now the assertion follows from
\begin{align*}
G_1(z)\geq 10^{-10}\prod_{w<p<z}\left(1-\frac{\omega(p)}{p}\right)^{-1}\geq 10^{-20}(\log N)^{\frac{3}{2}}(\log w)^{-\frac{3}{2}}.
\end{align*}

Part (ii) in turn follows by applying the Selberg sieve \cite[Chapter 7]{iwaniec-kowalski} to estimate
\begin{align*}
&G_1(z)\sum_{n\leq N}\bigg(\sum_{\substack{d\mid P(z)\\Wn+B\in \mathcal{A}_d}}\rho_d\bigg)^2\\
&\leq 10^{10}(\log N)^{\frac{3}{2}}(\log w)^{-\frac{3}{2}}\cdot \left(N\prod_{w<p<z}\left(1-\frac{\omega(p)}{p}\right)+z^3\right)\\
&\leq 10^{20} (\log N)^{\frac{3}{2}}(\log w)^{-\frac{3}{2}}\cdot N\left(\frac{\log w}{\log z}\right)^{\frac{3}{2}}\leq 10^{30} N.    
\end{align*}

Part (iii) is verified as follows. From the definition of $\rho_d$ it is clear that $|\rho_d|\leq 1$, so that
\begin{align}\label{eq79}
\beta(n)\leq G_1(z)\bigg(\sum_{\substack{d\mid P(z)\\Wn+B\in \mathcal{A}_d}}1\bigg)^2.
\end{align}
Note that if $Wn+B\in \mathcal{A}_p$ for some $w<p\leq z$, then $p\mid Wn+B$ or $p\mid Wn+B-1$, so that $p$ can be chosen in at most $\nu(Wn+B)+\nu(Wn+B-1)$ ways, where $\nu(\cdot)$ is the number of distinct prime factors. Since $d$ is square-free and a product of such primes $p$,  $d$ can be chosen in at most $2^{\nu(Wn+B)+\nu(Wn+B-1)}\ll N^{\frac{\varepsilon}{3}}$ ways in \eqref{eq79}. Therefore, \eqref{eq79} is $\ll (\log N)^{\frac{3}{2}}N^{\frac{2}{3}\varepsilon}\ll N^{\varepsilon}$.\\

Part (iv), which is the most crucial part concerning pseudorandomness, was verified in \cite{ramare}. Namely, our set of primes of the form $Wn+B=x^2+y^2+1$ is "sufficiently sifted" in the sense of the definition given on pages 1 and 2 of \cite{ramare} (to see that, take in that paper $A$ to be the set of primes of the form under consideration up to $N$ and $\kappa=\frac{3}{2}$). This property is all that is needed to obtain (iv) with the bound $v\left(\frac{a}{q}\right)\ll q^{-\frac{1}{2}}$, by formula (4.1.19) of \cite{ramare}. It is clear that this can be replaced with the stronger bound $v\left(\frac{a}{q}\right)\ll q^{\varepsilon-1}$, since we have defined the sets $\mathcal{K}_d$ in \eqref{eq80} so that formula (4.1.18) of \cite{ramare} holds for $\xi=\frac{\varepsilon}{2}$, instead of just some $0<\xi<\frac{1}{2}$.\\

We are then left with part (v). Equations (4.1.13) and (4.1.21) of \cite{ramare} reveal that \eqref{eq81} holds when $v(\frac{a}{q})$ is defined for $(a,q)=1$ by 
\begin{align*}
&v\left(\frac{a}{q}\right)=G_1(z)\sum_{q\mid [d_1,d_2]}\frac{\rho_{d_1}^{*}\rho_{d_2}^{*}}{[d_1,d_2]}|\mathcal{K}_{[d_1,d_2]}|\cdot \frac{\sum_{b\in \mathcal{K}_q}e\left(\frac{ab}{q}\right)}{|\mathcal{K}_q|}\,\,\text{with}\,\,\rho_{\ell}^{*}=\sum_{d\equiv 0\hspace{-0.1cm} \pmod \ell}\mu\left(\frac{d}{\ell}\right)\mu(d)\rho_d,
\end{align*}
where the set $\mathcal{K}_d$ is given by \eqref{eq80}. As in formula (4.1.17) of \cite{ramare}, we have
\begin{align*}
\left|\sum_{b\in \mathcal{K}_q}e\left(\frac{ab}{q}\right)\right|=\left|\sum_{b\in \mathbb{Z}_q\setminus\mathcal{K}_q}e\left(\frac{ab}{q}\right)\right|\leq |\mathbb{Z}_q\setminus \mathcal{K}_q|\leq \prod_{p^{\alpha}\mid \mid q}(p^{\alpha}-|\mathcal{K}_{p^{\alpha}}|),
\end{align*}
which immediately gives $v(\frac{a}{q})=0$ unless $q$ is square-free and $(q,W)=1$. In addition, by formula (4.1.13) of the same paper (with the right-hand side multiplied by $G_1(z)$), we have 
\begin{align}\label{eq82}
v\left(\frac{a}{q}\right)=G_1(z)w_{q}^{\#} \cdot \frac{\sum_{b\in \mathcal{K}_q}e\left(\frac{ab}{q}\right)}{|\mathcal{K}_q|},
\end{align}
where by (4.1.14) we have
\begin{align*}
w_{q}^{\#}=\frac{1}{G_1(z)}\sum_{\delta\leq z}h(\delta)\rho_z(q,\delta),
\end{align*}
and $\rho_z(q,\delta)$ satisfies (4.1.15). Putting $q=1$ into (4.1.15), we clearly get $w_1^{\#}=\frac{1}{G_1(z)}$, so that $v(1)=1$ by \eqref{eq82}.\qedd\\

We have now proved Proposition \ref{prop_restriction}, which will be needed in the proof of Theorem \ref{theo_goldbach}. As a consequence of the above considerations, we can now establish Theorem \ref{theo_roth}, that is, Roth's theorem for the subset $\mathscr{P}$ of primes.\\

\textbf{Proof of Theorem \ref{theo_roth}:} This is very similar to the proof of \cite[Theorem 1.2]{green-restriction}. Let $\mathcal{A}\subset \mathscr{P}^{*}$ have positive upper density in $\mathscr{P}^{*}$. Then there is $\delta>0$ (which may be assumed small) such that $|\mathcal{A}\cap (\frac{N}{3},\frac{2N}{3})|\geq \delta |\mathscr{P}^{*}\cap (\frac{N}{3},\frac{2N}{3})|$ for $N\in \mathcal{N}$, where $\mathcal{N}$ is some infinite set of positive integers.  Let $W$, $w$ and $J$ be as in \eqref{eq30} with $J=\lfloor \frac{10}{\delta}\rfloor$.\\

Let $S_B=S\cap \{Wn+B:\,\, n\geq 1\}$ for any set $S$ and integer $B$. Note that if $n=x^2+y^2+1\in (\frac{N}{3},\frac{2N}{3})$ is a prime with $(x,y)=1$ and $N\geq 10W$, then $(n,W)=(n-1,s(W))=1$ and $(n-1,3)=1$, $4\nmid n-1$. Therefore,
\begin{align*}
\sum_{\substack{1\leq B\leq W\\Wn+B \,\, \text{amenable}}} \left|\mathcal{A}_B\cap (\frac{N}{3},\frac{2N}{3})\right|&=\left|\mathcal{A}\cap (\frac{N}{3},\frac{2N}{3})\right|\geq \delta \left|\mathscr{P}^{*}\cap (\frac{N}{3},\frac{2N}{3})\right|,
\end{align*}
for $N\geq 10W$ and $N\in \mathcal{N}$, so using the pigeonhole principle and the lower bound for $|\mathscr{P}^{*}\cap (\frac{N}{3},\frac{2N}{3})|$ coming from Proposition \ref{prop_bohr} with $\chi\equiv 1$, we can find a value of $B\in [1,W]$ such that the polynomial $Wn+B$ is amenable and 
\begin{align}\label{eq102}
\left|\mathcal{A}_B\cap  (\frac{N}{3},\frac{2N}{3})\right|\geq \delta_1\cdot \delta(\log w)^{\frac{3}{2}}\frac{N}{W(\log N)^{\frac{3}{2}}}
\end{align} 
for $N\in \mathcal{N}'$  with $\mathcal{N}'$ an infinite set of positive integers and for some small absolute constant $\delta_1>0$, since the Chinese remainder theorem shows that there are $\leq 10^{10} W(\log w)^{-\frac{3}{2}}$ amenable functions $Wn+B$ with $1\leq B\leq W$.\\

Next, set 
\begin{align*}
g(n)=\delta_2(\log N)^{\frac{3}{2}}(\log w)^{-\frac{3}{2}}1_{\mathcal{A}_B\cap (\frac{N}{3},\frac{2N}{3})}(n)\quad \text{for}\quad N\in \mathcal{N}'\quad \text{and}\quad 1\leq n\leq N
\end{align*}
with $\delta_2>0$ small, and extend $g$ periodically to $\mathbb{Z}_N$. The assertion of the theorem will follow from the Green-Tao transference principle \cite[Proposition 5.1]{green-restriction} as soon as we check formulas (5.3)-(5.6) of that paper for the functions $g(n)$ and $\nu(n)=\beta(n)1_{[1,N]}(n)$ (extended periodically to $\mathbb{Z}_N$) with $\beta(\cdot)$ given by Proposition \ref{prop3}. We know (5.3) from Proposition \ref{prop3} and (5.6) from Proposition \ref{prop_restriction}. Formula (5.5) follows from the properties (i)-(v) of $\beta(n)$ just as in \cite[Chapter 6]{green-restriction}. We are left with (5.4), which follows (for a different value of $\delta$) for $N\in \mathcal{N}'$ from \eqref{eq102}. Now, as mentioned, \cite[Proposition 5.1]{green-restriction} yields the result, since any triple of the form $(a,a+d+j_1N,a+2d+j_2N)$ is an arithmetic progression in $\mathbb{Z}$ if $a,a+2d+j_1N, a+2d+j_2N\in (\frac{N}{3},\frac{2N}{3})$. \qedd

\section{Reductions for finding primes in Bohr sets}\label{Sec: reductions}

The proof of Proposition \ref{prop_bohr} goes through an intermediate result (namely Proposition \ref{prop2} below) that resembles  it and is slightly more technical, but at the same time easier to approach. The proof of Proposition \ref{prop2} uses among other things the circle method, Bombieri-Vinogradov type estimates, and ideas similar to Iwaniec's proof \cite{iwaniec-quadraticform} of the infinitude of primes $x^2+y^2+1$, and will occupy Sections  \ref{Sec: weighted} to \ref{Sec: hypotheses}.
\begin{proposition}\label{prop2}
Let $\chi:\mathbb{Z}\to \mathbb{R}_{\geq 0}$ have Fourier complexity $\mathcal{C}\ll 1$. Let $N\geq 1$ be an integer and $W$ be as in \eqref{eq30} with $w\geq \mathcal{C}^{20}$, and suppose that $Wn+b$ is an amenable linear function. There exists  an integer $Q\leq (\log N)^{B}$, depending only on $\chi$, with $B\ll_{\mathcal{C}} 1$, such that the following holds. For $N\geq N_0(w,\mathcal{C})$, $|t|\leq 5N$ and $c_0\in \mathcal{Q}$ we have
\begin{align*}
 \sum_{\substack{n\sim N\\n\equiv c_0\hspace{-0.1cm} \pmod Q\\Wn+b\in \mathbb{P}\\Wn+b-1\in \mathcal{S}}}\chi(t-n)\geq \frac{\delta_1}{(\log N)^{\frac{3}{2}}}\left(\frac{W}{\varphi(W)}\right)^{\frac{3}{2}}\frac{Q}{|\mathcal{Q}|}\bigg(\sum_{\substack{n\sim N\\ n\equiv c_0\hspace{-0.1cm} \pmod Q}}\chi(t-n)+o\left(\frac{N}{Q}\right)\bigg),
\end{align*}
where $\delta_1>0$ is an absolute constant and 
\begin{align}\label{eq21}
\mathcal{Q}=\{c_0\hspace{-0.2cm}\hspace{-0.1cm} \pmod Q:\,\, (Wc_0+b,Q)=(Wc_0+b-1,s(Q))=1\}.    
\end{align}
\end{proposition}

We remark that, by the Chinese remainder theorem,
\begin{align}\label{eq62}
|\mathcal{Q}|=Q\prod_{\substack{p\mid Q\\p\nmid W\\p\equiv 1\hspace{-0.1cm} \pmod 4}}\left(1-\frac{1}{p}\right)\prod_{\substack{p\mid Q\\p\nmid W\\p\equiv -1\hspace{-0.1cm}\pmod 4}}\left(1-\frac{2}{p}\right) ,
\end{align}
considering that $(b,W)=(b-1,s(W))=1$ by the definition of amenability.\\

In this section, we will show that Proposition \ref{prop2} implies Proposition \ref{prop_bohr}, by appealing to the following lemma. 

\begin{lemma}\label{lemma1} Let $\chi:\mathbb{Z}\to \mathbb{R}_{\geq 0}$ have Fourier complexity at most $\mathcal{C}$. Let $N,Q\geq 1$ be such that $N\geq 2Q^2$. Let $\mathcal{Q}$ be a collection of residue classes $\hspace{-0.1cm} \pmod Q$ such that for all $q\mid Q, q\neq 1$ and for  all $(a,q)=1$ we have
\begin{align*}
\bigg|\sum_{c_0\in \mathcal{Q}}e\left(\frac{a}{q}c_0\right)\bigg|\leq \eta_0 |\mathcal{Q}|    
\end{align*}
for some $\eta_0>0$. Then, with the same notations as in Proposition \ref{prop2}, for some absolute constant $C'>0$ and for all integers $t$ we have
\begin{align*}
\frac{Q}{|\mathcal{Q}|}\sum_{\substack{c_0\in \mathcal{Q}}}\sum_{\substack{n\sim N\\ n\equiv c_0\hspace{-0.1cm} \pmod Q}}\chi(t-n)\geq \sum_{n\sim N}\chi(t-n)-C'(\eta_0 \mathcal{C}^2N+Q\mathcal{C}^2N^{\frac{1}{2}}).    
\end{align*}
\end{lemma}

\textbf{Proof.} This is \cite[Lemma 7.4]{matomaki-shao}.\qedd\\

Note that the conclusion of Proposition \ref{prop_bohr} (with $\frac{N}{3}$ replaced with $N$) can be rewritten as 
\begin{align}\label{eq77}
\sum_{\substack{n\sim N\\Wn+b\in \mathbb{P}\\Wn+b-1\in \mathcal{S}}}\chi(t-n)\geq \frac{\delta_0}{(\log N)^{\frac{3}{2}}}\left(\frac{W}{\varphi(W)}\right)^{\frac{3}{2}}\bigg(\sum_{n\sim N}\chi(t-n)-\frac{CN}{w^{\frac{1}{3}}}\bigg),    
\end{align}
for $N\geq N_0(w,\mathcal{C})$ and $t\in \left(N,3N\right)$, with $\delta_0>0$ and $C>0$ absolute constants. In view of the previous lemma, Proposition \ref{prop_bohr} follows immediately from Proposition \ref{prop2} by splitting in \eqref{eq77} the sum over $n$ on the left-hand side to a sum over $n$ in different residue classes $\hspace{-0.1cm} \pmod Q$, provided that the premise of Lemma \ref{lemma1} is true for $\eta_0=w^{-\frac{1}{2}}$. This is what we will prove in the remainder of this section.

\begin{lemma}\label{lemma2} Let $Q\geq 1$, and let $\mathcal{Q}$ be defined by \eqref{eq21} (and $W$  and $w$  in the definition of $\mathcal{Q}$ given by \eqref{eq30}). Let $a$ and $q\mid Q$  be positive integers with $(a,q)=1$, $q\neq 1$. We have 
\begin{align}\label{eqq5}
\bigg|\sum_{c_0\in \mathcal{Q}}e\left(\frac{a}{q}c_0\right)\bigg|\leq w^{-\frac{1}{2}}|\mathcal{Q}|.    
\end{align}
\end{lemma}

Before proving this, we present another lemma, which will be used to prove Lemma \ref{lemma2}.

\begin{lemma}\label{lemma3} Let $a$ and $q$ be positive integers, $q \neq 1$, $(a,q)=1$, and let $Wn+b$ be an amenable linear polynomial with $W$ and $w$ as in \eqref{eq30}. Let $V\geq 1$ be an integer with $(q,V)=1$. Then
\begin{align}\label{eqq1}
\bigg|\sum_{\substack{n\hspace{-0.1cm} \pmod{q}\\(WVn+b,q)=1\\(WVn+b-1,s(q))=1}}e\left(\frac{a}{q}n\right)\bigg|\leq  \tau(q)\cdot 1_{(q,W)=1}.  
\end{align}
\end{lemma}

\textbf{Proof.} Using M\"obius inversion, the sum in question (without absolute values) becomes
\begin{align}\label{eqq4}
\sum_{d\mid q}\mu(d)\sum_{k\mid s(q)}\mu(k)\sum_{\substack{n\hspace{-0.1cm} \pmod{q}\\WVn\equiv -b\hspace{-0.1cm} \pmod d\\WVn\equiv -(b-1) \hspace{-0.1cm} \pmod k}}e\left(\frac{a}{q}n\right).    
\end{align}
Now consider the sum
\begin{align}\label{eqq2}
\sum_{\substack{n\hspace{-0.1cm} \pmod{q}\\WVn\equiv -b\hspace{-0.1cm} \pmod d\\WVn\equiv -(b-1) \hspace{-0.1cm} \pmod k}}e\left(\frac{a}{q}n\right).
\end{align}
Note that the sum is nonempty only if $(d,k)=1$. Let $x_{1},\ldots, x_{R(d,k)} \hspace{-0.1cm} \pmod{dk}$ be the pairwise incongruent solutions to the system $WVx\equiv -b\hspace{-0.1cm} \pmod d$, $WVx\equiv -(b-1)\hspace{-0.1cm} \pmod k$ (if there are none, the sum \eqref{eqq2} is empty). Since $dk=[d,k]\mid q$, after writing $n=x_j+dk t$ for some $1\leq j\leq R(d,k)$ and $1\leq t\leq \frac{q}{dk}$, \eqref{eqq2} transforms into
\begin{align}\label{eqq3}
\sum_{j=1}^{R(d,k)}\sum_{\substack{n\hspace{-0.1cm} \pmod{q}\\ n\equiv x_j \hspace{-0.1cm} \pmod{dk}}}e\left(\frac{an}{q}\right)&=\sum_{j=1}^{R(d,k)}e\left(\frac{ax_j}{q}\right)\sum_{t \hspace{-0.1cm} \pmod{\frac{q}{dk}}}e\left(\frac{at}{\frac{q}{dk}}\right).   
\end{align}
The inner sum is nonzero only when $dk=q$, in which case it is $1$. Taking these considerations into account, \eqref{eqq4} has absolute value at most
\begin{align}\label{eq34}
\sum_{\substack{d\mid q\\ k\mid s(q)\\ dk=q}}R(d,k)|\mu(d)||\mu(k)|.    
\end{align}
We estimate this differently depending on whether $(q,W)>1$ or $(q,W)=1$. In the former case, there is some prime $p$ such that $p\mid q$, $p\mid W$, so $ dk=q$ tells that $p$ divides either $d$ or $k$. If $p\mid d$, then supposing that $R(d,k)\neq 0$, the congruence $WVx\equiv -b\hspace{-0.1cm} \pmod p$ must be solvable. It however is not solvable, since $p\nmid b$ for $p\mid W$ by the amenability of $Wn+b$. If $p\mid k$, then $k\mid s(q)$ implies that $p\equiv -1\hspace{-0.1cm} \pmod 4$, $p\neq 3$. If $R(d,k)\neq 0$, the congruence $WVx\equiv -(b-1) \hspace{-0.1cm} \pmod p$ has a solution, but $p\nmid b-1$ by amenability, so we have a contradiction. We deduce that all the summands in \eqref{eq34} vanish for $(q,W)>1$.\\

Then let  $(q,W)=1$. As $d,k\mid q$ in \eqref{eq34}, we also have $(d,W)=(k,W)=1$ and $(d,V)=(k,V)=1$. Now clearly both of the congruences $WVx\equiv -b\hspace{-0.1cm} \pmod d$, $WVx\equiv -(b-1)\hspace{-0.1cm} \pmod k$ have a unique solution, so if the two congruences are thought of as a simultaneous equation, it has at most one solution $\hspace{-0.1cm} \pmod{dk}$. Therefore $R(d,k)\leq 1$, which leads to \eqref{eq34} being at most
\begin{align*}
\sum_{dk=q}1\leq \tau(q),    
\end{align*}
as asserted.\qedd\\

\textbf{Proof of Lemma \ref{lemma2}.} This is similar to the argument on page 21 of \cite{matomaki-shao}. We can find unique $q'$ and $Q'$ such that $Q=qq'Q'$ and $(q,Q')=1$ and all the prime divisors of $q'$ divide $q$. Writing $c_0=c_1q+c_2Q'$, $c_0$ runs through each residue class $\hspace{-0.1cm} \pmod{Q}$ exactly once as $c_1$ runs through residue classes $\hspace{-0.1cm} \pmod{q'Q'}$ and $c_2$ runs independently through residue classes $\hspace{-0.1cm} \pmod{q}$. Now the left-hand side of \eqref{eqq5} (without absolute values) becomes
\begin{align}\label{eqq6}
\Sigma:=\sum_{\substack{c_1\hspace{-0.1cm} \pmod{q'Q'}\\(Wqc_1+b,Q')=1\\(Wqc_1+b-1,s(Q'))=1}}\sum_{\substack{c_2\hspace{-0.1cm} \pmod{q}\\ (WQ'c_2+b,q)=1\\(WQ'c_2+b-1,s(q))=1}}e\left(\frac{aQ'}{q}c_2\right).    
\end{align}
Since $(aQ',q)=1$, the inner sum is exactly of the form appearing in Lemma \ref{lemma3}. Therefore, 
\begin{align*}
|\Sigma|\leq \sum_{\substack{c_1\hspace{-0.1cm} \pmod{q'Q'}\\(Wqc_1+b,Q')=1\\(Wqc_1+b-1,s(Q'))=1}}\tau(q)\cdot 1_{q>w}.
\end{align*}
Since $w\geq 10^{{10}^{10}}$, estimating the divisor function crudely yields
\begin{align*}
|\Sigma|\leq 1_{q> w}\cdot q^{0.1}\sum_{\substack{c_1\hspace{-0.1cm} \pmod{q'Q'}\\(Wqc_1+b,Q')=1\\(Wqc_1+b-1,s(Q'))=1}}1&=1_{q>w}\cdot q'q^{0.1}\sum_{\substack{c_1\hspace{-0.1cm} \pmod{Q'}\\(Wqc_1+b,Q')=1\\(Wqc_1+b-1,s(Q'))=1}}1\\
&= 1_{q>w}\cdot q'q^{0.1}Q'\prod_{\substack{p\mid Q'\\ p>w}}\left(1-\frac{\omega(p)}{p}\right)
\end{align*}
where $\omega(p)\in \{1,2\}$ and $\omega(p)=2$ precisely when $p\equiv -1\hspace{-0.1cm} \pmod 4$. The previous expression is, for $q>w\geq 10^{{10}^{10}}$,
\begin{align*}
&\leq q'q^{0.2}\prod_{\substack{p\mid q\\p>w}}\left(1-\frac{\omega(p)}{p}\right)\cdot Q'\prod_{\substack{p\mid Q'\\p>w}}\left(1-\frac{\omega(p)}{p}\right)\\
&=\frac{Q}{q^{0.8}}\prod_{\substack{p\mid Q\\p>w}}\left(1-\frac{\omega(p)}{p}\right)\leq \frac{|\mathcal{Q}|}{w^{\frac{1}{2}}},
\end{align*}
where the last step comes from \eqref{eq62}.\qedd\\

From Lemma \ref{lemma2}, we conclude that proving Proposition \ref{prop2} is enough for establishing Proposition \ref{prop_bohr} (and hence Theorem \ref{theo_goldbach}).

\section[Weighted sieve for primes]{Weighted sieve for primes of the form $p=x^2+y^2+1$} \label{Sec: weighted}

Next we investigate primes of the form $x^2+y^2+1$ in Bohr sets and prove Proposition \ref{prop2} concerning these, from which Theorem \ref{theo_goldbach} will follow. We will prove in this section Theorem \ref{t2} about weighted counting of primes in the shifted set $\mathcal{S}+1=\{s+1:s\in \mathcal{S}\}$. The proof resembles Iwaniec' s proof \cite{iwaniec-quadraticform} of the infinitude of primes of the form $x^2+y^2+1$, as well as the later works \cite{wu}, \cite{matomaki-m2+n2+1} on the same problem in short intervals, but the theorem involves a weighted version of the sieve procedure and hence requires a hypothesis about the weights. We will later verify the conditions of this hypothesis for a weight function related to the function $\chi(n)$ in Proposition \ref{prop2}, and this will imply Proposition \ref{prop2} and consequently Theorem \ref{theo_goldbach}. To formulate Theorem \ref{t2}, we first introduce the hypothesis regarding our weight coefficients. To this end, we need a couple of definitions.

\begin{definition}\label{def3}
Given a linear function $L$, let $\mathfrak{S}(L)$ be the \textit{singular product}
\begin{align*}
\mathfrak{S}(L)&=\prod_{\substack{p\equiv -1\hspace{-0.1cm} \pmod 4\\p\neq 3}}\left(1-\frac{|\{n\in \mathbb{Z}_p:\,\, L(n)\equiv 0\,\, \textnormal{or}\,\, 1 \hspace{-0.1cm} \pmod p\}|}{p}\right)\left(1-\frac{2}{p}\right)^{-1}\\
&\cdot\prod_{p\not \equiv -1\hspace{-0.1cm} \pmod 4}\left(1-\frac{|\{n\in \mathbb{Z}_p:\,\, L(n)\equiv 0\hspace{-0.1cm} \pmod p\}|}{p}\right)\left(1-\frac{1}{p}\right)^{-1}.
\end{align*}
\end{definition}

\begin{definition}\label{def2}
We say that a sequence $(g({\ell}))_{\ell\geq 1}$ of complex numbers is of \textit{convolution type} (for a given large integer $N$ and constant $\sigma\in (3,4)$) if 
\begin{align*}
g(\ell)=\sum_{\substack{\ell=km\\N^{\frac{1}{\sigma}}\leq k\leq N^{1-\frac{1}{\sigma}}}}\alpha_k\beta_m
\end{align*}
for some complex numbers $|\alpha_k|, |\beta_k|\leq \tau(k)^2\log k$. 
\end{definition}

\begin{definition}\label{def4} For $\frac{1}{3}<\rho_2<\rho_1<\frac{1}{2}$ and $\sigma\in (3,4)$, let $\textnormal{H}(\rho_1,\rho_2,\sigma)$ be the proposition
\begin{align}\label{eq96}
\frac{1}{2\sqrt{\rho_2}}\int_{1}^{\rho_2 \sigma}\frac{dt}{\sqrt{t(t-1)}}>\frac{1}{2\rho_1}\int_{2}^{\sigma}\frac{\log(t-1)}{t(1-\frac{t}{\sigma})^{\frac{1}{2}}}dt+10^{-10}.
\end{align}
 \end{definition}

In the proof of Theorem \ref{theo_goldbach}, we will use the fact that 
\begin{align*}
\text{H}\left(\frac{1}{2}-\varepsilon,\frac{3}{7}-\varepsilon,3+\varepsilon\right)\quad \text{is true for small enough}\quad  \varepsilon>0.
\end{align*}
This holds for $\varepsilon=0$ by a numerical computation and by continuity in a small neighborhood of $0$. Indeed, the difference between the integrals in \eqref{eq96} is then $>10^{-3}$. We are ready to state our Bombieri-Vinogradov type hypothesis, whose validity depends on the weight sequence $(\omega_n)$, as well as on the parameters $\rho_1,\rho_2$ and $\sigma$.

\begin{hypothesis}\label{h1} Let  $L(n)=Kn+b$ be an amenable linear function with $K\ll (\log N)^{O(1)}$.  Let $(\omega_n)_{n\sim N}$ be a nonnegative sequence of real numbers, and let $\delta=(b-1,K)$. Let $\varepsilon>0$ be any small number. Let $\frac{1}{3}<\rho_2<\rho_1<\frac{1}{2}-\varepsilon$, $\sigma\in (3,4)$.  Then for any sequence $(g(\ell))_{\ell\leq N^{0.9}}$ of convolution type (with parameter $\sigma$) 
\begin{align*}
\sum_{\substack{d\leq N^{\rho_1}\\(d,K)=1}}\lambda_d^{+,\textnormal{LIN}}\sum_{\substack{\ell\leq N^{0.9}\\(\ell,K)=\delta\\(\ell,d)=1}}g(\ell)\bigg(\sum_{\substack{n\sim N\\L(n)=\ell p+1\\L(n)\equiv 0\hspace{-0.1cm} \pmod d}}\omega_n-\frac{1}{\varphi(d)}\frac{K}{\varphi(\frac{K}{\delta})}\sum_{n\sim N}\frac{\omega_n}{\ell \log \frac{Kn}{\ell}}\bigg)&\ll \frac{\sum_{n\sim N}\omega_n}{(\log N)^{100}},\\
\sum_{\substack{d\leq N^{\rho_2}\\ (d,K)=1}}\lambda_d^{-,\textnormal{SEM}}\bigg(\sum_{\substack{n\sim N\\L(n)\in \mathbb{P}\\L(n)\equiv 1\hspace{-0.1cm} \pmod d}}\omega_n-\frac{1}{\varphi(d)}\frac{ K}{\varphi(K)}\sum_{n\sim N}\frac{\omega_n}{\log(Kn)}\bigg)&\ll \frac{\sum_{n\sim N}\omega_n}{(\log N)^{100}},
\end{align*}
where $\lambda_{d}^{+,\textnormal{LIN}}$ are the upper bound linear sieve weights with sifting parameter $z_1= N^{\frac{1}{5}}$ and $\lambda_{d}^{-,\textnormal{SEM}}$ are the lower bound semilinear sieve weights with sifting parameter $z_2=N^{\frac{1}{\sigma}}$ (the weights $\lambda_d^{\pm,\textnormal{SEM}}$ were defined in Theorem \ref{theo_sievebombieri}, and the weights $\lambda_d^{\pm,\textnormal{LIN}}$ are defined analogously by replacing $\beta=1$ by $\beta=2$ in that definition). 
\end{hypothesis}

\begin{theorem}\label{t2} Assume Hypothesis \ref{h1} for a linear form $L(n)$, sequence $(\omega_n)_{n\sim N}$, and parameters $\rho_1,\rho_2,\sigma$ satisfying $\textnormal{H}(\rho_1,\rho_2,\sigma)$. Then
\begin{align*}
\sum_{\substack{n\sim N\\L(n)\in \mathbb{P}\\L(n)-1\in \mathcal{S}}}\omega_n\geq \frac{\delta_0\cdot \mathfrak{S}(L)}{(\log N)^{\frac{3}{2}}}\sum_{n\sim N}\omega_n+O(N^{\frac{1}{2}}),
\end{align*}
where $\delta_0>0$ is an absolute constant.
\end{theorem}

\begin{remark} We will be able to prove Hypothesis \ref{h1} in Section \ref{Sec: hypotheses} for $\rho_1=\frac{1}{2}-\varepsilon$, $\rho_2=\frac{3}{7}-\varepsilon$ and $\sigma=3+\varepsilon$ when $L(n)$ is suitable and $\omega_n$ is of bounded Fourier complexity. It would suffice to prove the same with $\rho_2=0.385$ instead of $\rho_2=\frac{3}{7}-\varepsilon=0.428\ldots$ (since then $\text{H}(\rho_1,\rho_2,\sigma)$ is true). On the other hand, existing Bombieri-Vinogradov estimates such as \cite[Lemma 12]{tolev_bombieri} would only give us $\rho_2=\frac{1}{3}-\varepsilon=0.333\ldots$, which falls short of what we need.
\end{remark}

\textbf{Proof.} Put
\begin{align*}
\mathcal{A}&=\{L(n)-1:\, n\sim N, L(n)\in \mathbb{P}\}\\
\mathcal{P}_{4,-1}&=\{p\in \mathbb{P}:\, p\equiv -1\hspace{-0.1cm} \pmod 4,\, p\neq 3\},\\
P(z)&=\prod_{\substack{p< z\\p\in \mathcal{P}_{4,-1}}}p,\\
\mathcal{P}_{4,1}^{*}&=\{n\geq 1:\, p\mid n\Rightarrow p\equiv 1 \hspace{-0.1cm} \pmod 4\}.
\end{align*}
If we weigh the elements of $\mathcal{A}$ by $\nu_n=\omega_{(L^{-1}(n+1))}$, where $L^{-1}$ is the inverse function of $L$, the sifting function is
\begin{align*}
S(\mathcal{A},\mathcal{P}_{4,-1},z)=\sum_{\substack{n\sim N\\L(n)\in \mathbb{P}\\(L(n)-1,P(z))=1}}\omega_n.
\end{align*}
Note that $L(n)-1\equiv 2^{\beta} \hspace{-0.1cm} \pmod{2^{\beta+2}}$ for some $\beta\geq 1$ by the definition of amenability, so that $L(n)-1$ has an even number of prime factors that are $\equiv -1\hspace{-0.1cm} \pmod 4$ (counted with multiplicity). We have
\begin{align}\label{eq85}
\sum_{\substack{n\sim N\\L(n)\in \mathbb{P}\\L(n)-1\in \mathcal{S}}}\omega_n&=S(\mathcal{A},\mathcal{P}_{4,-1},(3KN)^{\frac{1}{2}}),
\end{align}
since the right-hand side counts with weight $\omega_n$ the numbers  $L(n)-1=2^{\alpha_1}3^{\alpha_2}k\in \mathcal{A}$ with $k\in \mathcal{P}_{4,1}^{*}$, and we claim that these numbers are precisely the numbers in $\mathcal{S}\cap \mathcal{A}$. We have $2^{\alpha_1}3^{\alpha_2}k=L(n)-1$, so by amenability $\alpha_2\equiv 0\hspace{-0.1cm} \pmod 2$. It is a fact in elementary number theory that for $k\in \mathcal{P}_{4,1}^{*}$, both $k$ and $2k$ can be expressed in the form $a^2+b^2$ with $(a,b)=1$, and additionally no number of the form $2^{\alpha_1}3^{\alpha_2}k$ with $(k,6)=1$ and $\alpha_2$ odd or $k\not \in \mathcal{P}_{4,1}^{*}$ is of the form $x^2+y^2$ with $(x,y)\mid 6^{\infty}$. Hence both sides of \eqref{eq85} indeed count the same integers.\\

 Buchstab's identity reveals that
\begin{align*}
S(\mathcal{A},\mathcal{P}_{4,-1},(3KN)^{\frac{1}{2}})&=S(\mathcal{A},\mathcal{P}_{4,-1},N^{\frac{1}{\sigma}})-\sum_{\substack{n\sim N\\L(n)\in \mathbb{P}}}\sum_{\substack{p_2\mid L(n)-1\\N^{\frac{1}{\sigma}}\leq p_2<(3KN)^{\frac{1}{2}}\\(L(n)-1,P(p_2))=1\\p_2\in \mathcal{P}_{4,-1}}}\omega_n.
\end{align*}
The condition $p_2\mid L(n)-1\equiv 2^{\beta} \hspace{-0.1cm} \pmod{2^{\beta+2}}$ implies that $L(n)-1$ has either exactly $2$ prime divisors from $\mathcal{P}_{4,-1}$ or at least $4$ such prime divisors (with multiplicities). The second case is impossible, since all the prime divisors of $L(n)-1$ that are from $\mathcal{P}_{4,-1}$ are $\geq p_2$ and $p_2^{4}\geq N^{\frac{4}{\sigma}}>L(2N)-1.$ This means that we may write $L(n)-1=p_1p_2m',$ $p_1\geq p_2,$ $p_1\in \mathcal{P}_{4,-1},$ with $m'$ having no prime divisors from $\mathcal{P}_{4,-1}$. Now $\delta\mid L(n)-1=Kn+b-1$ with $\delta=(b-1,K)$, and since $p_1\geq p_2\geq N^{\frac{1}{\sigma}}>K$, we have $\delta\mid m'$. Hence we may write $m'=\delta m$, where $m\in \mathcal{P}_{4,1}^{*}$ (we have $3\nmid m$, since $K$ is divisible by a larger power of $3$ than $b-1$ is, by the definition of amenability. Similarly $2\nmid m$). We claim that $(m,\frac{K}{\delta})=1$. Indeed, if $p\mid m$ and $p\mid \frac{K}{\delta}$, we must have $p\mid \frac{b-1}{\delta}$, a contradiction to $(K,b-1)=\delta$. Now we have
\begin{align}\label{eq39}
S(\mathcal{A},\mathcal{P}_{4,-1},(3KN)^{\frac{1}{2}})= S-T.
\end{align}
Here
\begin{align*}
S=S(\mathcal{A},\mathcal{P}_{4,-1},N^{\frac{1}{\sigma}}),\quad T=\sum_{\substack{n\sim N\\L(n)\in \mathbb{P}}}\sum_{\substack{L(n)-1=\delta  p_1p_2m\\p_1,p_2\in \mathcal{P}_{4,-1}\\N^{\frac{1}{\sigma}}\leq p_2\leq p_1\\m\in \mathcal{P}_{4,1}^{*}}}\omega_n\leq \sum_{\ell \in \mathcal{L}}S(\mathcal{M}(\ell),\mathcal{P}(\ell),N^{\frac{1}{6}}),
\end{align*} 
with
\begin{align*}
\mathcal{L}&=\{\delta  p_2m:\,\, N^{\frac{1}{\sigma}}\leq p_2\leq (3KNm^{-1})^{\frac{1}{2}},\,\, p_2\in \mathcal{P}_{4,-1},\,\, m\in \mathcal{P}_{4,1}^{*},\,\, (m,\frac{K}{\delta})=1\},\\
\mathcal{M}(\ell)&=\{L(n):\, L(n)=\ell p+1: n\sim N, p\in \mathbb{P}\},\\
\mathcal{P}(\ell)&=\{p\in \mathbb{P}:(p,2\ell)=1\},\quad Q(z)=\prod_{\substack{p< z\\ p\in \mathcal{P}(\ell)}}p,
\end{align*}
and $M(\ell)$ has been assigned the weights $\nu_n=\omega_{L^{-1}(n)}$, so that
\begin{align*}
S(\mathcal{M}(\ell),\mathcal{P}(\ell),z)=\sum_{\substack{n\sim N\\L(n)=\ell p+1\\(L(n),Q(z))=1}}\omega_n.
\end{align*}
We carry out bounding $S$ from below and bounding $T$ from above  separately.\\

\textbf{Bounding $S$.} For $d\mid P(z)$, $(d,K)=1$, let
\begin{align*}
r(\mathcal{A},d)&=\sum_{\substack{n\sim N\\L(n)\in \mathbb{P}\\L(n)-1\equiv 0\hspace{-0.1cm} \pmod d}}\omega_n-\frac{1}{\varphi(d)}\frac{K}{\varphi(K)}\sum_{n\sim N}\frac{\omega_n}{\log(Kn)},
\end{align*}
and for $(d,K)>1$ we let $r(\mathcal{A},d)=0$ (since if $p\mid d$, $p\mid K$ and $p\in \mathcal{P}_{4,-1}$, then $p$ does not divide any element of $\mathcal{A}$ by the amenability of $L(n)$). Let $\sigma\in (3,4)$ be as in Hypothesis \ref{h1}. The semilinear sieve \cite[Theorem 11.13]{friedlander}, with $\beta=1$, sifting parameter $z=N^{\frac{1}{\sigma}}$, and level $D=z^s$, $1\leq s\leq 2$,
gives
\begin{align}\label{eq35}\begin{split}
&S(\mathcal{A},\mathcal{P}_{4,-1},N^{\frac{1}{\sigma}}) \\
&\geq \frac{K}{\varphi(K)}\sum_{n\sim N}\frac{\omega_n}{\log(Kn)} V_K^{\textnormal{SEM}}(N^{\frac{1}{\sigma}}) \left(f(s)+O((\log N)^{-0.1})\right)+\sum_{d\leq N^{\frac{s}{\sigma}}}\lambda_d^{-,\textnormal{SEM}}r(\mathcal{A},d),  \end{split}  
\end{align}
where $\lambda_d^{-,\textnormal{SEM}}$ are the lower bound semilinear weights with sifting parameter $z=N^{\frac{1}{\sigma}}$ and we have introduced the quantities
\begin{align*}
f(s)=\sqrt{\frac{e^{\gamma}}{\pi s}}\int_{1}^s \frac{dt}{\sqrt{t(t-1)}}\quad \text{and}\quad V_K^{\textnormal{SEM}}(z)=\prod_{\substack{p<z\\p\equiv -1 \hspace{-0.1cm} \pmod 4\\p\nmid K}}\left(1-\frac{1}{\varphi(p)}\right)
\end{align*}
We take $s=\rho_2 \sigma \in [1,2]$, where $\rho_2$ is as in Hypothesis \ref{h1}. Now Hypothesis \ref{h1} permits replacing the last sum in \eqref{eq35} with an error of $\ll \frac{\sum_{n\sim N}\omega_n}{(\log N)^{100}}$ (since the terms of that sum in \eqref{eq35} vanish unless $(d,K)=1$). Moreover, the term $V_K^{\textnormal{SEM}}(N^{\frac{1}{\sigma}})$ can be computed asymptotically using  \cite[Proposition 1]{wu}, which implies that
\begin{align*}
V_K^{\textnormal{SEM}}(z)=(1+o(1))\prod_{\substack{p\mid K\\p\equiv -1 \hspace{-0.1cm} \pmod 4}}\left(1-\frac{1}{p-1}\right)^{-1}\cdot 2AC_{4,-1}\cdot \left(\frac{\pi e^{-\gamma}}{\log z}\right)^{\frac{1}{2}},
\end{align*}
where
\begin{align*}
A=\frac{1}{2\sqrt{2}}\prod_{p\equiv -1 \hspace{-0.1cm} \pmod 4}\left(1-\frac{1}{p^2}\right)^{\frac{1}{2}}\quad \text{and}\quad C_{4,i}=\prod_{p\equiv i \hspace{-0.1cm} \pmod 4}\left(1-\frac{1}{(p-1)^2}\right)
\end{align*}
for $i\in \{-1,1\}$. Therefore, we end up with the bound
\begin{align}\label{eq44}
S&\geq \frac{4AC_{4,-1}+o(1)}{(\log N)^{\frac{1}{2}}} \cdot I_1(\rho_2,\sigma)\frac{K}{\varphi(K)}\prod_{\substack{p\mid K\\ p\equiv -1 \hspace{-0.1cm} \pmod 4}}\left(1-\frac{1}{p-1}\right)^{-1}\cdot \sum_{n\sim N}\frac{\omega_n}{\log(Kn)} \nonumber  \\
&=\frac{4AC_{4,-1}+o(1)}{(\log N)^{\frac{3}{2}}} \cdot I_1(\rho_2,\sigma)\frac{K}{\varphi(K)}\prod_{\substack{p\mid K\\ p\equiv -1 \hspace{-0.1cm} \pmod 4}}\left(1-\frac{1}{p-1}\right)^{-1}\cdot \sum_{n\sim N}\omega_n,
\end{align}
where
\begin{align*}
I_1(\rho_2,\sigma)&=\frac{1}{2\sqrt{\rho_2}}\int_{1}^{\rho_2 \sigma}\frac{dt}{\sqrt{t(t-1)}}.
\end{align*}
\textbf{Bounding $T$.} Write, for $d\mid Q(z)$, $(d,K)=1$, $(\ell,d)=1$ and $(\ell,K)=\delta$, 
\begin{align*}
r(\mathcal{M}(\ell),d)=\sum_{\substack{n\sim N\\L(n)-1= \ell p\\L(n)\equiv 0\hspace{-0.1cm} \pmod d}}\omega_n-\frac{1}{\varphi(d)}\frac{K}{\varphi(\frac{K}{\delta})}\sum_{n\sim N}\frac{\omega_n}{\ell \log \frac{Kn}{\ell}}.
\end{align*}
For all other $d$ such that $d\mid Q(z)$, let $r(\mathcal{M}(\ell),d)=0$ (since if $(d,K)>1$, then $L(n)-1=\ell p,$ $L(n)\equiv 0\hspace{-0.1cm} \pmod d$ is impossible). With these notations, for $1\leq s \leq 3$ the linear sieve \cite[Theorem 11.13]{friedlander} with $\beta=2$ provides the bound
\begin{align}\label{eq1}
S(\mathcal{M}(\ell),\mathcal{P}(\ell),N^{\frac{1}{6}})&\leq \frac{(1+o(1))K}{\varphi(\frac{K}{\delta})}\sum_{n\sim N}\frac{\omega_n}{\ell \log \frac{Kn}{\ell}} V_K^{\textnormal{LIN}}(N^{\frac{1}{5}},\ell)F(s)+\sum_{d\leq N^{\frac{s}{5}}}\lambda_d^{+,\textnormal{LIN}}r(\mathcal{M}(\ell),d),
\end{align}
where $\lambda_d^{+,\textnormal{LIN}}$ are the upper bound linear sieve coefficients with sifting parameter $z=N^{\frac{1}{5}}$, $F(s)=\frac{2e^{\gamma}}{s}$, and
\begin{align*}
V_K^{\textnormal{LIN}}(z,\ell)&=\prod_{\substack{p\in \mathcal{P}(\ell)\\ p< z\\p\nmid K}}\left(1-\frac{1}{\varphi(p)}\right)=\prod_{2< p<z}\left(1- \frac{1}{p-1}\right)\prod_{\substack{p\mid K\ell\\2<p<z}}\left(1-\frac{1}{p-1}\right)^{-1}.\nonumber
\end{align*}
Applying formula (4.6) of \cite{wu}, we get the asymptotic
\begin {align}\label{eq48}
&V_K^{\textnormal{LIN}}(z,\ell)=(1+o(1))\frac{2C_{4,1}C_{4,-1}e^{-\gamma}\mathfrak{f}(K \ell)}{\log z},\textnormal{where}\quad \mathfrak{f}(d)=\prod_{\substack{p\mid d\\p>2}}\left(1-\frac{1}{p-1}\right)^{-1}.
\end{align}
We take $s=5\rho_1\in [1,3]$ in the linear sieve. Then we have
\begin{align}\label{eq75}
\sum_{\ell \in \mathcal{L}}\sum_{d\leq N^{\rho_1}}\lambda_d^{+,\textnormal{LIN}}r(\mathcal{M}(\ell),d)&=\sum_{\substack{d\leq N^{\rho_1}\\(d,K)=1}}\lambda_d^{+,\textnormal{LIN}}\sum_{\substack{\ell\leq N^{\frac{3}{4}+\varepsilon}\\(\ell,d)=1\\(\ell,K)=\delta}}1_{\mathcal{L}}(\ell)r(\mathcal{M}(\ell),d),
\end{align}
since $1_{\mathcal{L}}(\ell)$ is supported on $\ell\leq 3K^2N^{1-\frac{1}{\sigma}}\leq N^{\frac{3}{4}+\varepsilon}$. Concerning the error sum in \eqref{eq1}, observe that
\begin{align*}
1_{\mathcal{L}}(\ell)=\sum_{\substack{\ell=k\cdot \delta m\\N^{\frac{1}{\sigma}}\leq k\leq (3KN)^{\frac{1}{2}}\\k\leq \left(\frac{3KN}{m}\right)^{\frac{1}{2}}}}1_{\mathcal{P}_{4,-1}}(k)1_{\mathcal{P}_{4,1}^{*}}(m)1_{(m,\frac{K}{\delta})=1},
\end{align*}
so $1_{\mathcal{L}}(\ell)$ is of convolution type (for the value of $\sigma$ we are considering), except for the cross condition $k\leq \left(\frac{3KN}{m}\right)^{\frac{1}{2}}$. We use Perron's formula in the form 
\begin{align*}
1_{(1,\infty)}(y)&=\frac{1}{\pi}\int_{-N^4}^{N^4}\frac{\sin(t \log y)}{t}dt+O\left(\frac{1}{N^4|\log y|}\right)\\
&=\frac{2}{\pi}\int_{N^{-5}}^{N^4}\frac{\sin(t \log y)}{t}dt+O\left(\frac{1}{N^4|\log y|}+\frac{|\log y|}{N^5}\right)
\end{align*}
for $N^{-3}<y\leq N^3, y\neq 1$ to dispose of the cross condition. We choose $y=\frac{3KN}{k^2m}$, which satisfies $|y-1|\geq \frac{1}{3KN^2}$ after altering $N$ by $\leq 1$ if necessary, so that the error term in Perron's formula becomes $O(\frac{K}{N^2})$. According to the addition formula for sine, we have
\begin{align*}
\sin (t\log y)=\sin(t\log (3KN)-t\log k^2)\cos(t\log m)-\cos(t\log (3KN)-t\log k^2) \sin(t\log m)
\end{align*}
which permits us to separate the variables $k$ and $m$. Then we have
\begin{align*}
1_{\mathcal{L}}(\ell)=\frac{2}{\pi}\int_{N^{-4}}^{N^3}\frac{1}{t}\sum_{\substack{\ell=k\cdot \delta m\\N^{\frac{1}{\sigma}}\leq k\leq (3KN)^{\frac{1}{2}}}}(\alpha_k^{(1)}(t)\beta_m^{(1)}(t)-\alpha_k^{(2)}(t)\beta_m^{(2)}(t))\, dt+O\left(\frac{1}{N^{2-\varepsilon}}\right),
\end{align*}
where $|\alpha_k^{(j)}(t)|, |\beta_m^{(j)}(t)|\leq 1$ and $t\mapsto \alpha_k^{(j)}(t)$ and $t\mapsto \beta_m^{(j)}(t)$ are continuous and $\alpha_k^{(j)}(t)$ is supported on $N^{\frac{1}{\sigma}}\leq k\leq (3KN)^{\frac{1}{2}}$. Substituting this to \eqref{eq75}, Hypothesis \ref{h1} tells that
\begin{align*}
\sum_{\ell \in \mathcal{L}}\sum_{d\leq N^{\rho_1}}\lambda_d^{+,\textnormal{LIN}}r(\mathcal{M}(\ell),d)\ll \frac{\sum_{n\sim N}\omega_n}{(\log N)^{99}}+O(N^{\frac{1}{2}-\varepsilon}).
\end{align*}
We sum \eqref{eq1} over $\ell \in \mathcal{L}$ and make use of \eqref{eq48}, after which we have obtained 
\begin{align*}
&\sum_{\ell \in \mathcal{L}}S(\mathcal{M}(\ell),\mathcal{P}(\ell),N^{\frac{1}{5}})\\
& \leq (F(s)+o(1))\cdot \frac{K}{\varphi(\frac{K}{\delta})}\sum_{n\sim N}\sum_{\ell \in \mathcal{L}}\frac{\omega_n}{\ell \log \frac{Kn}{\ell}} V_K^{\textnormal{LIN}}(N^{\frac{1}{6}},\ell)+O\left(\frac{\sum_{n\sim N}\omega_n}{(\log N)^{99}}\right)\\
&=\left(\frac{2e^{\gamma}}{5\rho_1}+o(1)\right)\cdot \frac{K}{\varphi(\frac{K}{\delta})}\sum_{\ell\in \mathcal{L}}\frac{\mathfrak{f}(K \ell)}{\ell \log \frac{KN}{\ell}}\cdot \sum_{n\sim N}\omega_n\cdot\frac{2C_{4,1}C_{4,-1}e^{-\gamma}}{\frac{1}{5}\log N}+O\left(\frac{\sum_{n\sim N}\omega_n}{(\log N)^{99}}\right).  
\end{align*}

We analyze the sum over $\mathcal{L}$ in the above formula. Denoting $\mathcal{L}'=\{\frac{\ell}{\delta}:\, \ell\in \mathcal{L}\}$, it is 
\begin{align*}
\sum_{\ell \in \mathcal{L}}\frac{\mathfrak{f}(K\ell)}{\ell \log \frac{KN}{\ell}}&=\left(\frac{1}{\delta}+o(1)\right)\sum_{\ell' \in \mathcal{L}'}\frac{\mathfrak{f}(K\ell')}{\ell'\log \frac{KN}{\ell'}}1_{(\ell',\frac{K}{\delta})=1},
\end{align*}
since $\delta\mid K$. The previous sum can be written as 
\begin{align}\label{eq83}
(1+o(1))\sum_{m\leq N^{1-\frac{2}{\sigma}+\varepsilon}}\frac{u(m)\mathfrak{f}(Km)1_{(m,\frac{K}{\delta})=1}}{m}\sum_{\substack{N^{\frac{1}{\sigma}}\leq p\leq (\frac{3KN}{m})^{\frac{1}{2}}\\p\equiv -1\hspace{-0.1cm} \pmod 4}}\frac{1}{p\log \frac{N}{pm}},
\end{align}
where $u(m)$ is the characteristic function of $\mathcal{P}_{4,1}^{*}$. To evaluate this sum, we study the sum
\begin{align}\label{eq83a}
\sum_{m\leq x}u(m)\mathfrak{f}(Km)1_{(m,\frac{K}{\delta})=1}.
\end{align}
The sum can be written as
\begin{align*}
\mathfrak{f}(K)\sum_{m\leq x}u(m)\mathfrak{f}(\psi_K(m))1_{(m,\frac{K}{\delta})=1}, \quad \text{where}\quad \psi_K(m)=\prod_{\substack{p\mid m\\p\nmid K}}p, 
\end{align*}
and the advantage is that $\mathfrak{f}(\psi_K(m))$ is a multiplicative function. By Wirsing's theorem \cite[Satz 1]{wirsing} applied to the nonnegative multiplicative function $h(m)=u(m)\mathfrak{f}(\psi_K(m))1_{(m,\frac{K}{\delta})=1}$ (which is bounded by $2$ at prime powers and fulfills $\sum_{p\leq x}h(p)\log p=(\frac{1}{2}+o(1))x$), we see that \eqref{eq83a} equals
\begin{align*}
&(\mathfrak{f}(K)+o(1))\frac{e^{-\frac{\gamma}{2}}}{\sqrt{\pi}}\frac{x}{\log x}\prod_{\substack{p\leq x\\p\nmid \frac{K}{\delta}\\p\equiv 1 \hspace{-0.1cm} \pmod 4}}\left(1+\frac{h(p)}{p}+\frac{h(p^2)}{p^2}+\cdots\right)\\
&=(\mathfrak{f}(K)+o(1))\frac{e^{-\frac{\gamma}{2}}}{\sqrt{\pi}}\frac{x}{\log x}\prod_{\substack{p\leq x\\p\nmid K\\p\equiv 1 \hspace{-0.1cm} \pmod 4}}\left(1+\frac{1}{p-2}\right)\prod_{\substack{p\mid K\\p\nmid \frac{K}{\delta}\\p\equiv 1 \hspace{-0.1cm} \pmod 4}}\left(1-\frac{1}{p}\right)^{-1}.
\end{align*}
Applying Wirsing's theorem reversely, this is 
\begin{align*}
&(\mathfrak{f}(K)+o(1))\sum_{m\leq x}u(m)\mathfrak{f}(m)\cdot \prod_{\substack{p\mid K\\p\equiv 1 \hspace{-0.1cm} \pmod 4}}\left(1+\frac{1}{p-2}\right)^{-1}\prod_{\substack{p\mid K\\p\nmid \frac{K}{\delta}\\p\equiv 1 \hspace{-0.1cm} \pmod 4}}\left(1-\frac{1}{p}\right)^{-1}.
\end{align*}
By  \cite[Lemma 3]{wu}, we have
\begin{align*}
\sum_{m\leq x}u(m)\mathfrak{f}(m)=(1+o(1))\frac{A}{C_{4,1}}\frac{x}{(\log x)^{\frac{1}{2}}}.
\end{align*}
Now, using the same argument  as in the proof of \cite[Lemma 5]{matomaki-m2+n2+1}, we compute that \eqref{eq83} equals 
\begin{align*}
\frac{A+o(1)}{C_{4,1}(\log N)^{\frac{1}{2}}}\cdot \frac{1}{2}\int_{2}^{\sigma}\frac{\log(t-1)}{t(1-\frac{t}{\sigma})^{\frac{1}{2}}}dt\cdot \frac{\mathfrak{f}(K)}{\delta}\prod_{\substack{p\mid K\\p\equiv 1 \hspace{-0.1cm} \pmod 4}}\left(1+\frac{1}{p-2}\right)^{-1}\prod_{\substack{p\mid K\\p\nmid \frac{K}{\delta}\\p\equiv 1 \hspace{-0.1cm} \pmod 4}}\left(1-\frac{1}{p}\right)^{-1}.
\end{align*}

\textbf{Concluding the proof.} Now we have
\begin{align}\label{eq22}
T\leq \frac{4AC_{4,-1}+o(1)}{(\log N)^{\frac{3}{2}}} \frac{ I_2(\rho_1,\sigma)K\mathfrak{f}(K)}{\delta\varphi(\frac{K}{\delta})}\prod_{\substack{p\mid K\\p\equiv 1 \hspace{-0.1cm} \pmod 4}}\hspace{-0.1cm}\left(1+\frac{1}{p-2}\right)^{-1}\hspace{-0.2cm}\prod_{\substack{p\mid K\\p\nmid \frac{K}{\delta}\\p\equiv 1 \hspace{-0.1cm} \pmod 4}}\hspace{-0.1cm}\left(1-\frac{1}{p}\right)^{-1}\sum_{n\sim N}\omega_n,    
\end{align}
where
\begin{align*}
I_2(\rho_1, \sigma)=\frac{1}{2\rho_1}\int_{2}^{\sigma}\frac{\log(t-1)}{t(1-\frac{t}{\sigma})^{\frac{1}{2}}}dt.    
\end{align*}
We claim that the local factors in \eqref{eq44} and \eqref{eq22} are identical, or in other words that
\begin{align}\label{eq84}\begin{split}
&\prod_{p\mid K}\left(1-\frac{1}{p}\right)^{-1}\prod_{\substack{p\mid K\\ p\equiv -1\hspace{-0.1cm} \pmod 4}}\left(1-\frac{1}{p-1}\right)^{-1}\\
&=\prod_{p\mid \frac{K}{\delta}}\left(1-\frac{1}{p}\right)^{-1}\prod_{\substack{p\mid K\\p>2}}\left(1-\frac{1}{p-1}\right)^{-1}\prod_{\substack{p\mid K\\p\equiv 1 \hspace{-0.1cm} \pmod 4}}\left(1+\frac{1}{p-2}\right)^{-1}\prod_{\substack{p\mid K\\p\nmid \frac{K}{\delta}\\p\equiv 1 \hspace{-0.1cm} \pmod 4}}\left(1-\frac{1}{p}\right)^{-1}.\end{split}
\end{align}
By the identity $(1+\frac{1}{p-2})^{-1}=1-\frac{1}{p-1}$,  \eqref{eq84} is equivalent to
\begin{align*}
\prod_{p\mid K}\left(1-\frac{1}{p}\right)^{-1}=\prod_{p\mid \frac{K}{\delta}}\left(1-\frac{1}{p}\right)^{-1}\prod_{\substack{p\mid K\\p\nmid \frac{K}{\delta}\\p\equiv 1 \hspace{-0.1cm} \pmod 4}}\left(1-\frac{1}{p}\right)^{-1},
\end{align*}
which in turn is equivalent to the nonexistence of a prime $p\not\equiv 1\hspace{-0.1cm} \pmod 4$ for which $p\mid K$, $p\nmid \frac{K}{\delta}$. If $p\geq 5$ were such a prime, we would have $p\mid  \delta$, so $p\mid b-1$, which contradicts the definition of amenability. We also cannot have $p=2$ or $p=3$, since $2\mid \frac{K}{\delta}$ and $3\mid \frac{K}{\delta}$ for $\delta=(b-1,K)$ by amenability.\\ 

Thus no such $p$ exists and \eqref{eq84} holds. Furthermore, it is clear that \eqref{eq84} is at least $0.01\mathfrak{S}(L)$. Consequently, 
\begin{align*}
 S-T\geq (0.01+o(1))4AC_{4,-1}\mathfrak{S}(L)(I_1(\rho_2,\sigma)-I_2(\rho_1,\sigma))\frac{\sum_{n\sim N}\omega_n}{(\log N)^{\frac{3}{2}}}+O(N^{\frac{1}{2}}).
\end{align*}
Owing to the fact that $\text{H}(\rho_1,\rho_2,\sigma)$ is assumed to be true, we have $I_1(\rho_2,\sigma)-I_2(\rho_1,\sigma)\geq 10^{-10}$, and  this completes the proof of Theorem \ref{t2} in view of \eqref{eq85} and \eqref{eq39}.\qedd

\section{Preparation for the verifying the hypothesis}\label{Sec: decomposition}

The sequence $(\omega_n)$ to which we will apply Theorem \ref{t2} will be determined by a function $\chi(n)$ having a Fourier series of the form \eqref{eq13}. In \eqref{eq13} it is natural to separate the phases $\alpha_i$ into major and minor arc parameters. This partition arises from the following lemma. 

\begin{lemma}\label{le11}
Let $\alpha_1,\ldots, \alpha_{\mathcal{C}}$ be real numbers with $\mathcal{C}\ll 1$, and let $W\ll 1$ be as in \eqref{eq30}. Also let the constants $A,B\geq 1$ be related by $B=A(3\mathcal{C})^{\mathcal{C}}$. Then for any large $N$ there exists a positive integer $Q\leq (\log N)^{B}$ such that each $\alpha_k$ may be written as
\begin{align*}
\alpha_k&=W\frac{a_k}{q_k}+\varepsilon_k,\,\, (a_k,q_k)=1,\,\, 1\leq q_k\leq \frac{N}{(\log N)^{100B}},\,\,|\varepsilon_k|\leq \frac{(\log N)^{100B}W}{q_k N},     
\end{align*}
and either $q_k\mid Q$ or $q_k\geq \frac{q_k}{(q_k,Q^2)}\geq (\log N)^A$.
\end{lemma}

\textbf{Proof.} This is Lemma 3.2 in  \cite{matomaki-shao}.\qedd\\

From now on, $A$ (and therefore also $B$) will be large enough quantities (say $A,B\geq 10^{10}$). Let us define the sequence $(\omega_n)$ to which we will apply Theorem \ref{t2} in order to prove Proposition \ref{prop2}.  Let $\chi:\mathbb{Z}\to \mathbb{R}_{\geq 0}$ be any function with Fourier complexity $\leq \mathcal{C}$ (i.e., $\chi$ satisfies \eqref{eq13}). Given an integer $t$ with $|t|\leq 5N$, we choose 
\begin{align*}
(\omega_n)_{n\sim \frac{N}{Q}}=(\chi(t-(Qn+c_0)))_{n\sim \frac{N}{Q}},    
\end{align*}
where $Q$ is determined by the $\alpha_i$ in \eqref{eq13} with the help of Lemma \ref{le11} and $c_0\in \mathcal{Q}$ with
\begin{align*}
\mathcal{Q}=\{c_0\hspace{-0.2cm}\hspace{-0.1cm} \pmod Q:\,\, (Wc_0+b,Q)=(Wc_0+b-1,s(Q))=1\}.   
\end{align*} 
Recall that $|\mathcal{Q}|$ is given by \eqref{eq62}.\\

From now on, let
\begin{align*}
x=\frac{N}{Q},\quad L(n)=QWn+Wc_0+b,\quad c_0\in \mathcal{Q}.
\end{align*}
To prove Proposition \ref{prop2} and hence Proposition \ref{prop_bohr} and Theorem \ref{theo_goldbach}, it suffices to show that for $W$ as in \eqref{eq30} and $\mathfrak{S}(L)$ as in Definition \ref{def3} we have 
\begin{align}\label{eq12}
\sum_{\substack{n\sim x\\L(n)\in \mathbb{P}\\L(n)-1\in \mathcal{S}}}\chi(t-(Qn+c_0))\geq \frac{\delta_0\cdot \mathfrak{S}(L)}{(\log x)^{\frac{3}{2}}}\sum_{n\sim x}\chi(t-(Qn+c_0))+o\left(\frac{x}{(\log x)^{\frac{3}{2}}}\right), 
\end{align}
since $L(n)$ is amenable and since by \eqref{eq62}
\begin{align*}
\mathfrak{S}(L)&\asymp \prod_{\substack{p\equiv -1\hspace{-0.1cm} \pmod 4\\ p\mid QW\\p\nmid W}}\left(1-\frac{1}{p}\right)^{-2}\prod_{\substack{p\not \equiv -1 \hspace{-0.1cm} \pmod 4\\ p\mid QW\\p\nmid W}}\left(1-\frac{1}{p}\right)^{-1}\\
&\cdot  \prod_{\substack{p\equiv -1\hspace{-0.1cm} \pmod 4\\p\mid W}}\left(1-\frac{1}{p}\right)^{-2}\prod_{\substack{p\not \equiv -1\hspace{-0.1cm} \pmod 4\\p\mid W}}\left(1-\frac{1}{p}\right)^{-1}\\
&\asymp \left(\frac{W}{\varphi(W)}\right)^{\frac{3}{2}}\frac{Q}{|\mathcal{Q}|}.
\end{align*}
By Theorem  \ref{t2} and the remark after it, formula \eqref{eq12} will follow once we have verified Hypothesis \ref{h1} for our sequence $(\chi(t-(Qn+c_0)))_{n\sim x}$ and linear function $L(n)$ and parameters
\begin{align}\label{eq87}
\rho_1=\frac{1}{2}-10\varepsilon,\quad \rho_2=\frac{3}{7}-10\varepsilon,\quad \text{and}\quad \sigma=3+\varepsilon.    
\end{align}
By formula \eqref{eq13} for $\chi(n)$ and Lemma \ref{le11}, it suffices to inspect Hypothesis \ref{h1} with the choices \eqref{eq87} for $(e(\xi n))_{n\sim x}$, where $\xi$ is an arbitrary real number satisfying, for some $Q\leq 2(\log x)^{B}$,
\begin{align}\label{eq88}
\left|\xi-\frac{QWa}{q}\right|\leq \frac{2(\log x)^{102B}}{qx}\,\, \text{for}\,\, (a,q)=1,\,\, q\leq  \frac{x}{(\log x)^{99B}},\,\, \text{and}\,\, q\mid Q\,\, \text{or}\,\, \frac{q}{(q,Q^2)}\geq (\log x)^{A}.
\end{align}
Moreover, we may assume in \eqref{eq12} that
\begin{align*}
\sum_{n\sim x}\chi(t-(Qn+c_0))\gg \frac{x}{(\log x)\mathfrak{S}(L)},  
\end{align*}
since otherwise we have nothing to prove, and consequently it suffices to prove Hypothesis \ref{h1} for  $(e(\xi n))_{n\sim x}$ with $(\sum_{n\sim x}\omega_n) (\log x)^{-100}$ replaced by $x(\log x)^{-200}$ in that hypothesis.

\section{Bombieri-Vinogradov sums weighted by additive characters}\label{Sec: Bombieri}

We will establish Hypothesis \ref{h1} in the setting of Section \ref{Sec: decomposition} subsequently in Section \ref{Sec: hypotheses}.  For that purpose as well as for proving Theorem \ref{theo_alphap} in Section \ref{Sec: fractional parts}, we need the following Bombieri-Vinogradov type estimates for type I and II exponential sums. We employ for positive integers $q$ and $v$ the notation
\begin{align*}
q_{v}=\frac{q}{(q,v^2)}.
\end{align*}

\begin{lemma}\label{le8}  Let $M\leq N^{0.4}$, $R\leq N^{0.1}$, and $\rho\leq \frac{1}{2}-\varepsilon$ for some $\varepsilon\in (0,\frac{1}{6})$. Let $\xi$ be a real number with $|\xi-\frac{a}{q}|\leq \frac{1}{(qv)^2}$ for some coprime $a$ and $q\in [1,N]$ and some positive integer $v\leq N^{0.1}$. Then for any complex numbers $|\alpha_m|\leq \tau(m)^2\log m$ and any $t\in [N,2N]$ we have
\begin{align*}
&\sum_{0<|r|\leq R}\,\,\sum_{d\leq N^{\rho}}\max_{(c,dv)=1}\bigg|\sum_{\substack{N\leq mn\leq t\\mn\equiv c \hspace{-0.1cm} \pmod{dv}\\m\leq M}}\alpha_m e(\xi r mn)\bigg|\\
&\ll \left(\frac{RN}{v}\right)^{\frac{1}{2}}\left(RMN^{\rho}+\frac{RN}{vq_v}+q_v\right)^{\frac{1}{2}}(\log N)^{1000}.
\end{align*}
\end{lemma}

\textbf{Proof.} We follow the proof of  \cite[Lemma 8.3]{matomaki-shao}. It suffices to consider the sum over $0<r\leq R$. Our task is to estimate 
\begin{align*}
S_r=\sum_{d\leq N^{\rho}}\max_{(c,dv)=1}\bigg|\sum_{\substack{N\leq mn\leq t\\mn\equiv c \pmod{dv}\\m\leq M}}\alpha_m e(\xi rmn)\bigg|
\end{align*}
for $r\leq R$. The inner sum in the definition of $S_r$ is a geometric sum in the variable $n$, so evaluating it provides the bound
\begin{align*}
S_r\ll \sum_{d\leq N^{\rho}}\sum_{m\leq M}|\alpha_m|\min\left\{\frac{RN}{rmdv},\frac{1}{\|r \xi md v\|}\right\}.
\end{align*}
Observe that $\left|v\xi-\frac{av}{q}\right|\leq \frac{1}{q^2}$. Based on this, writing $d'=rmd$ and using a standard bound for sums over fractional parts \cite[Lemma B.3]{matomaki-shao} (taking $x=\frac{RN}{v}$ in that lemma), we get
\begin{align*}
\sum_{r\leq R}S_r&\ll \sum_{d'\leq RMN^{\rho}}\tau(d')^5\min\left\{\frac{RN}{d' v},\frac{1}{\|v\xi d'\|}\right\}(\log N)\\
&\ll \left(\frac{RN}{vq_{v}^{\frac{1}{2}}}+\left(\frac{RN\cdot RMN^{\rho}}{v}\right)^{\frac{1}{2}}+\left(\frac{RN}{v}q_v\right)^{\frac{1}{2}}\right)(\log N)^{1000}\\
&\ll \left(\frac{RN}{v}\right)^{\frac{1}{2}}\left(RMN^{\rho}+\frac{RN}{vq_v}+q_v\right)^{\frac{1}{2}}(\log N)^{1000},
\end{align*}
as wanted.\qedd

\begin{lemma}\label{le9}  Let $M\in [N^{\frac{1}{2}},N^{\frac{3}{4}}]$ and $\Delta_1,\Delta_2\geq 1$, $\Delta_1\Delta_2\leq N^{\frac{1}{2}}$, $\Delta_1 \Delta_2^2\leq \frac{M}{v}$ for some positive integer $v\leq N^{0.1}$.  Let $\xi$ be a real number with $|\xi-\frac{a}{q}|\leq \frac{1}{(qv)^2}$ for some coprime $a$ and $q\in [1,N]$. 
Then  for any complex numbers $|\alpha_m|,|\beta_m|\leq \tau(m)^2\log m$ and any integer $c'\neq 0$ and number $t\in [N,2N]$ we have
\begin{align*}
&\sum_{0<|r|\leq R}\,\,\sum_{\substack{d_1\sim \Delta_1}}\sum_{\substack{d_2\sim \Delta_2\\(d_2,c'd_1v)=1}}\max_{(c,d_1v)=1}\bigg|\sum_{\substack{N\leq mn\leq t\\mn\equiv c \hspace{-0.1cm}\pmod{d_1v}\\mn\equiv c'\hspace{-0.1cm} \pmod{d_2}\\m\sim M}}\alpha_m\beta_ne(\xi r mn)\bigg|\\
&\ll \frac{RN}{v}\min\{F_1,F_2\}(\log N)^{1000},
\end{align*}
with
\begin{align*}
F_1&=\left(\frac{\Delta_1 Mv}{N}+\Delta_1\Delta_2^2\frac{v}{M}\right)^{\frac{1}{2}}+\left(\frac{1}{\Delta_1}+\frac{1}{q_{v}}+\frac{q_{v}v^2}{RN}\right)^{\frac{1}{8}},\\
F_2&=\Delta_1 \Delta_2\left(\frac{1}{q_v^{\frac{1}{2}}}+\frac{v}{M^{\frac{1}{2}}}+\frac{v^2M}{N}+\frac{q_{v}^{\frac{1}{2}}v}{(RN)^{\frac{1}{2}}}\right)^{\frac{1}{2}}.
\end{align*}
\end{lemma}

\begin{remark} In Section \ref{Sec: hypotheses}, we will only need the case $R=1$, while the dependence on $v$ will be crucial. In Section \ref{Sec: fractional parts}, on the other hand, $v=1$ but the dependence on $R$ will be crucial.
\end{remark}

\textbf{Proof.} We follow the proof of  \cite[Lemma 8.4]{matomaki-shao}, which in turn is based on an argument of Mikawa \cite{mikawa-bombieri}. It suffices to consider the case $r>0$. We will first prove the lemma in the case $F_1=\min\{F_1,F_2\}$. Let us write 
\begin{align*}
I_r=\sum_{d_1\sim \Delta_1}\sum_{\substack{d_2\sim \Delta_2\\(d_2,c'd_1v)=1}}\max_{(c,d_1v)=1}\bigg|\sum_{\substack{N\leq mn\leq t\\mn\equiv c \hspace{-0.1cm} \pmod{d_1v}\\mn\equiv c'\pmod{d_2}\\m\sim M}}\alpha_m\beta_ne(\xi r mn)\bigg|,
\end{align*}
so that $\sum_{r\leq R}I_r$ is what we are interested in. Since $\Delta_1\Delta_2^2\leq \frac{M}{v}$, a formula on page 37 of \cite{matomaki-shao} tells (with $x=N,$ $D=\Delta_1$, $\alpha=r\xi$) that
\begin{align*}
I_r^2\ll N(\log N)^{100}\left(\Delta_1\sum_{d_1\sim \Delta_1}\sum_{0<|j|\leq \frac{8\Delta_2^2 N}{\Delta_1 M v}}\tau_3(j)\min\left\{\frac{RN}{r(d_1v)^2|j|},\frac{1}{\|r \xi(d_1v)^2 |j|\|}\right\}+\frac{\Delta_1 M}{v}\right)
\end{align*}
(since the term $\frac{x^2}{Q^2}(\log x)^{-C+10}$ present in that formula of \cite{matomaki-shao} can be replaced with $\frac{DMx}{Q}(\log x)^{100}$ without changing anything in the proof). Using the Cauchy-Schwarz inequality, we obtain 
\begin{align}\label{eq65}
&\frac{1}{(\log N)^{200}}\sum_{r\leq R}I_r\leq \frac{1}{(\log N)^{200}}R^{\frac{1}{2}}\left(\sum_{r\leq R}I_r^2\right)^{\frac{1}{2}}\nonumber\\
&\leq (RN)^{\frac{1}{2}}\left(\Delta_1 \sum_{d_1\sim \Delta_1}\sum_{0<|j|\leq \frac{8\Delta_2^2 N}{\Delta_1Mv}}\sum_{r\leq R}\tau_3(j)\min\left\{\frac{RN}{r(d_1v)^2 |j|},\frac{1}{\|r\xi (d_1v)^2 |j|\|}\right\}+\frac{\Delta_1 RM}{v}\right)^{\frac{1}{2}}\nonumber\\
&\ll (RN)^{\frac{1}{2}}\left(\Delta_1 \sum_{d_1\sim \Delta_1}\sum_{1\leq \ell\leq \frac{8\Delta_2^2 R N}{\Delta_1 Mv}}\tau_4(\ell)\min\left\{\frac{RN}{(d_1v)^2\ell},\frac{1}{\|v^2\xi d_1^2\ell\|}\right\}+\frac{\Delta_1 RM}{v}\right)^{\frac{1}{2}},
\end{align}
after writing $\ell=rj$. When it comes to the sum above, we can estimate it using the lemma on page 6 of \cite{mikawa-bombieri} (with $\tau_3(\cdot)$ replaced by $\tau_4(\cdot)$), stating that
\begin{align}\label{eq64}
&\Delta_1\sum_{d_1\sim \Delta_1}\sum_{\ell\sim J}\tau_4(\ell)\min\left\{\frac{x}{d_1^2\ell},\frac{1}{\|\xi' d_1^2 \ell\|}\right\}\ll (\Delta_1^2 J+x^{\frac{3}{4}}(q'+\frac{x}{q'}+\frac{x}{\Delta_1})^{\frac{1}{4}})(\log x)^{100}
\end{align}
for $1\leq J\leq 10x$ and any real number $\xi'$ satisfying $|\xi'-\frac{a'}{q'}|\leq \frac{1}{q'^2}$ for some coprime $a'$ and $q'\leq x$. In the case $q'>x$, \eqref{eq64} continues to hold, by trivial estimates. We substitute \eqref{eq64} with $x=\frac{RN}{v^2}$, $\xi'=v^2\xi$ and $J\leq \frac{8\Delta_2^2RN}{\Delta_1 Mv}$ into \eqref{eq65} (we have $J\leq 10\frac{RN}{v^2}$ since $\Delta_1\Delta_2^2\leq \frac{M}{v}$), making use of our assumption on $\xi$, which implies that $\left|v^2\xi-\frac{\frac{av^2}{(q,v^2)}}{q_{v}}\right|\leq \frac{1}{q_v^2}$. This results in the claimed bound.\\

Then let $F_2=\min\{F_1,F_2\}$.  In this situation, we use the orthogonality of characters to bound the sum in Lemma \ref{le9} with
\begin{align}\label{eq85a}
&\sum_{r\leq R}\sum_{d_1\sim \Delta_1}\sum_{d_2\sim \Delta_2} \max_{\psi \hspace{-0.1cm}\hspace{-0.1cm} \pmod{d_1d_2}}\bigg|\sum_{\substack{N\leq mn\leq t\\mn\equiv c_{v}(d_1,d_2) \pmod v\\m\sim M}}\alpha_m \psi(m)\beta_n \psi(n)e(\xi r mn)\bigg|,
\end{align}
where $c_v(d_1,d_2)$ is a suitably chosen integer coprime to $v$. Estimating the sums over $d_1$ and $d_2$ trivially and using the Cauchy-Schwarz inequality and expanding a square, we find that \eqref{eq85a} is, for some $|\beta_n'|\leq \tau(n)^2 \log n$ and some $c_v$ coprime to $v$,
\begin{align}\label{eq106}
&\leq  \Delta_1 \Delta_2 (RM)^{\frac{1}{2}}\bigg(\sum_{r\leq R}\sum_{m\leq M}\bigg|\sum_{\substack{\frac{N}{m}\leq n\leq \frac{t}{m}\\n\equiv c_vm^{-1}\hspace{-0.1cm}\pmod v}}\beta_n' e(\xi r mn)\bigg|^2 \bigg)^{\frac{1}{2}}(\log M)^{100}\nonumber\\
&= \Delta_1 \Delta_2 (RM)^{\frac{1}{2}}\bigg(\sum_{r\leq R}\sum_{\substack{\frac{N}{2M}\leq n_i\leq \frac{2N}{M}\\n_1\equiv n_2\pmod v\\\text{for}\,\, i\in \{1,2\}}}\beta_{n_1}'\overline{\beta_{n_2}'}\sum_{\substack{m\leq M\\\frac{N}{n_i}\leq m\leq \frac{t}{n_i}\\m\equiv c_vn_i^{-1}\hspace{-0.1cm}\pmod v\\\text{for}\,\,i\in \{1,2\}}}e(\xi r m(n_1-n_2))\bigg)^{\frac{1}{2}}(\log M)^{100}\nonumber\\
&\ll \Delta_1 \Delta_2 (RN)^{\frac{1}{2}}\left(RM+\sum_{r\leq R}\sum_{\substack{1\leq n\leq \frac{2N}{M}\\n\equiv 0\pmod v}}T(n)\min\left\{\frac{RN}{rnv}+1,\frac{1}{\|v\xi r n\|}\right\}\right)^{\frac{1}{2}}(\log M)^{101},
\end{align}
where
\begin{align*}
T(n)=\frac{M}{N}\sum_{\substack{n=n_1-n_2\\n_1,n_2\leq \frac{2N}{M}}}\tau(n_1)^2\tau(n_2)^2.    
\end{align*}
We can write $n=kv$ and $\ell=kr$ to bound \eqref{eq106} with
\begin{align}\label{eq103}
\ll \Delta_1 \Delta_2 (RN)^{\frac{1}{2}}\left(RM+\sum_{\ell\leq \frac{2RN}{Mv}}U(\ell)\min\left\{\frac{RN}{\ell v^2}+1,\frac{1}{\|v^2\xi \ell\|}\right\}\right)^{\frac{1}{2}}(\log N)^{101},
\end{align}
where 
\begin{align*}
U(\ell)=\sum_{\substack{\ell=\ell_1 \ell_2\\\ell_1\leq \frac{2N}{Mv}}}T(\ell_1 v).    
\end{align*}
We apply \cite[Lemma B.3]{matomaki-shao} (with $k=20$) to \eqref{eq103}. The weight function $U(\ell)$ is not a divisor function, but the only property of the weight function needed in that lemma is a second moment bound. Therefore, \eqref{eq103} can be bounded with
\begin{align}\label{eq105}
&\ll  \Delta_1 \Delta_2 (RN)^{\frac{1}{2}}\bigg(\frac{RN}{q_{v}^{\frac{1}{2}}v^2}+\frac{RN}{(v^2M)^{\frac{1}{2}}}+RM+\bigg(\frac{RNq_v}{v^2}\bigg)^{\frac{1}{2}}\bigg)^{\frac{1}{2}}(\log N)^{1000},
\end{align}
once we prove that
\begin{align}\label{eq107}
\sum_{\ell\leq \frac{2RN}{Mv}}U(\ell)^2\ll \frac{RN}{Mv}(\log N)^{100}.    
\end{align}
We calculate
\begin{align}\label{eq104}
&\sum_{\ell\leq \frac{2RN}{Mv}}\left(\sum_{\substack{\ell=\ell_1\ell_2\\\ell_1\leq \frac{2N}{Mv}}}T(\ell_1v)\right)^2\ll\frac{RN}{Mv}\sum_{\substack{\ell_1\leq \frac{2N}{Mv}\\\ell_1'\leq \frac{2N}{Mv}}}\frac{T(\ell_1v)T(\ell_1'v)}{[\ell_1,\ell_1']}\nonumber\\
&\ll \frac{RN}{Mv}\sum_{d\leq \frac{2N}{Mv}}\frac{1}{d}\sum_{\substack{\ell_1\leq \frac{2N}{dMv}\\\ell_1'\leq \frac{2N}{dMv}}}\frac{T(\ell_1 dv)T(\ell_1' dv)}{\ell_1\ell_1'}=\frac{RN}{Mv}\sum_{d\leq \frac{2N}{Mv}}\frac{1}{d}\left(\sum_{\ell\leq \frac{2N}{dMv}}\frac{T(\ell dv)}{\ell}\right)^2.
\end{align}
We can estimate the sum inside the square using
\begin{align*}
&\sum_{\substack{n\leq \frac{2N}{M}\\n\equiv 0\pmod{c}}}\frac{T(n)}{n}\ll \frac{M}{N}\sum_{\substack{n_1\leq \frac{2N}{M}\\n_2\leq \frac{2N}{M}\\n_1\equiv n_2\pmod c\\n_1>n_2}}\frac{\tau(n_1)^2\tau(n_2)^2}{n_1-n_2}\\
&\ll\frac{M}{Nc}\sum_{1\leq a\leq c}\sum_{\substack{n_1'\leq \frac{2N}{Mc}\\n_2'\leq \frac{2N}{Mc}\\n_1'>n_2'}}\frac{\tau(cn_1'+a)^2\tau(cn_2'+a)^2}{n_1'-n_2'}\ll \frac{M}{Nc}\sum_{1\leq a\leq c}\sum_{n\leq \frac{2N}{Mc}}\tau(cn+a)^4\\
&\ll \frac{M}{Nc}\sum_{m\leq \frac{2N}{M}+c}\tau(m)^4\ll \frac{1}{c}(\log N)^{15},
\end{align*}
for $c\leq \frac{2N}{M}$, where we used Hilbert's inequality \cite[Chapter 7]{montgomery} in the third last step. Taking $c=dv$, and substituting to \eqref{eq104}, we see that \eqref{eq107} holds, as claimed. Therefore, we indeed have the bound \eqref{eq105} for \eqref{eq103}, and that bound can be rewritten as the desired bound $F_2$.\qedd

\section{Factorizing sieve weights}\label{Sec: sieveweight}

The linear and semilinear sieve weights will play a crucial role in verifying Hypothesis \ref{h1}, since we aim to split the summation over $d\leq x^{\rho}$ in that hypothesis to summations over $d_1\sim \Delta_1$, $d_2\sim \Delta_2$ for various values of $\Delta_1$ and $\Delta_2$. If such a factorization can be done, it provides more flexibility in our Bombieri-Vinogradov sums, and hence gives better bounds. This advantage can be seen from Lemma \ref{le9}, which often produces better bounds when $\Delta_1$ and $\Delta_2$ are of somewhat similar size, as opposed to the choice $\Delta_1=x^{\rho}$, $\Delta_2=1$. The following lemmas about the  combinatorial structure of sieve weights have been tailored so that the estimate given by Lemma \ref{le9} will be $\ll Nv^{-1}(\log N)^{-1000}$ if $\Delta_1$ and $\Delta_2$ satisfy the conditions for $d_1$ and $d_2$ in Lemma \ref{le10} or \ref{le1} with $D=\frac{x^{1-\varepsilon^2}}{M}$, $\theta=0$, $R=1$ and $q$ suitably large, and additionally $\rho=\frac{3}{7}(1-4\theta)-\varepsilon$ in the case of Lemma \ref{le10} or $\rho=\frac{1}{2}(1-4\theta)-\varepsilon$ in the case of Lemma \ref{le1}.  It should be remarked that in Section \ref{Sec: hypotheses} we will only need the case $\theta=0$ of the following lemmas, but for the proof of Theorem \ref{theo_alphap} we will choose $\theta=\frac{1}{80}-\varepsilon$.

\subsection{Linear sieve weights}\label{sub: linear}

\begin{lemma}\label{le10} Let $\varepsilon>0$ be small,  $0\leq \theta\leq \frac{1}{30}$, and $\rho=\frac{1}{2}(1-4\theta)-\varepsilon$. Let
\begin{align*}
\mathcal{D}^{+,\textnormal{LIN}}=\{p_1\cdots p_r\leq x^{\rho}:\,\, z_1\geq p_1> \ldots > p_r,\,\,p_1\cdots p_{2k-2}p_{2k-1}^3\leq x^{\rho}\,\, \textnormal{for all}\,\, k\geq 1\}   
\end{align*}
be the support of the upper bound linear sieve weights with level $x^{\rho}$ and sifting parameter $z_1\leq x^{\frac{1}{2}}$. Then, for any $D\in [x^{\frac{1}{5}},x^{\rho}]$, every $d\in \mathcal{D}^{+,\textnormal{LIN}}$ can be written as $d=d_1d_2$, where the positive integers $d_1$ and $d_2$ satisfy $d_1\leq D$, $d_1d_2^2\leq \frac{x^{1-4\theta-2\varepsilon^2}}{D}$. Moreover, we can take either $d_1\geq x^{0.1}$ or $d_2=1$.
\end{lemma}

\textbf{Proof.} The proof is similar to the proof of \cite[Lemma 12.16]{friedlander} (which essentially says that the linear sieve weights $\lambda^{+,\textnormal{LIN}}_d$ are well-factorable for any sifting parameter $z\leq x^{\frac{1}{2}-\varepsilon}$). We will actually show that any $d=p_1\cdots p_r\in \mathcal{D}^{+,\textnormal{LIN}}$ can be written as $d=d_1d_2$ with $d_1\leq D$, $d_2\leq \frac{x^{\rho}}{D}$ and either $d_1\geq x^{0.1}$ or $d_2=1$. After that statement has been proved, we have proved the lemma, because then $d_1d_2^2\leq \frac{x^{2\rho}}{D}\leq \frac{x^{1-4\theta-2\varepsilon^2}}{D}$. We use induction on $r$ to prove the existence of such $d_1$ and $d_2$. For $r=1$, we can simply take $d_1=p_1$ and $d_2=1$, since $p_1\leq x^{\frac{\rho}{3}}\leq x^{\frac{1}{6}}$. If $r=2$, we can take $d_1=p_1p_2$ , $d_2=1$, unless $p_1p_2>D$. In the case $p_1p_2>D$, in turn, the choice $d_1=p_1$, $d_2=p_2$ works, since $p_1\leq x^{\frac{1}{6}}$ and $p_2\leq \frac{x^{\rho}}{p_1p_2}\leq \frac{x^{\rho}}{D}$. Suppose then that $r\geq 3$ and that case $r-1$ has been proved and consider the case $r$. We  have $p_1\cdots p_{r-1}\in \mathcal{D}^{+,\textnormal{LIN}}$, so by the induction assumption $p_1\cdots p_{r-1}=d_1'd_2'$ with $d_1'\leq D$, $d_2'\leq \frac{x^{\rho}}{D}$ and either $d_1'\geq x^{0.1}$ or $d_2'=1$. We claim that we can take either $d_1=d_1'p_r$, $d_2=d_2'$ or $d_1'=d_1$, $d_2=d_2'p_r$. Firstly, if $d_1'< x^{0.1}$, then $d_2'=1$ and $d_1'=p_1\cdots p_{r-1}$. Since $r\geq 3$, this yields $p_1p_2<x^{0.1}$, so $p_2<x^{0.05}$. Now the choice $d_1=d_1'p_r$, $d_2=d_2'=1$ works because $d_1< x^{0.1}p_r\leq x^{0.15}\leq D$. Secondly, if in the opposite case $d_1'\geq x^{0.1}$ neither of the choices for $(d_1,d_2)$ works, then $d_1'd_2'p_r^2>x^{\rho}$. However, $d_1'd_2'p_r^2=p_1\cdots p_{r-1}p_r^2\leq x^{\rho}$ by the definition of $\mathcal{D}^{+,\textnormal{LIN}}$, so we have a contradiction and the induction works.\qedd

\subsection{Semilinear sieve weights}      \label{sub: semilinear}

\begin{lemma}\label{le1} Let $\varepsilon>0$ be small, $0\leq \theta\leq \frac{1}{30}$, and $\rho=\frac{3}{7}(1-4\theta)-\varepsilon$. Let
\begin{align*}
\mathcal{D}^{-,\textnormal{SEM}}=\{p_1\cdots p_r\leq x^{\rho}:\,\, z_2\geq p_1> \ldots > p_r,\,\,p_1\cdots p_{2k-1}p_{2k}^2\leq x^{\rho}\,\, \textnormal{for all}\,\, k\geq 1\}.    
\end{align*}
be the support of the lower bound semilinear sieve weights with level $x^{\rho}$ and sifting parameter $z_2\leq x^{\frac{1}{3}-2\theta-2\varepsilon^2}$. Then, for any $D\in [x^{\frac{1}{3}-2\theta-2\varepsilon^2},x^{\rho}]$, every $d\in \mathcal{D}^{-,\textnormal{SEM}}$ can be written as $d=d_1d_2$, where the positive integers $d_1$ and $d_2$ satisfy $d_1\leq D$, $d_1d_2^2\leq \frac{x^{1-4\theta-2\varepsilon^2}}{D}$. Moreover, we can take either $d_1\geq x^{0.1}$ or $d_2=1$.
\end{lemma}

\begin{remark}The exponent $\rho=\frac{3}{7}(1-4\theta)-\varepsilon$ is optimal in Lemma \ref{le1}. Namely, if $\rho=\frac{3}{7}(1-4\theta)+3\varepsilon$, then the lemma is false for $D=x^{\frac{3}{7}(1-4\theta)}$ and $p_1p_2p_3\in \mathcal{D}^{-,\textnormal{SEM}}$, $p_1,p_2,p_3\sim \frac{1}{2}x^{\frac{1}{7}(1-4\theta)+\varepsilon}$. 
\end{remark}

\begin{remark} We remark that an argument almost identical to the proof of Lemma \ref{le1} below shows that the lemma holds also for the set
\begin{align*}
\mathcal{D}^{+,\textnormal{SEM}}=\{p_1\cdots p_r\leq x^{\rho}:\,\, x^{\frac{1}{2}}\geq p_1\geq \ldots \geq p_r,\,\,p_1\cdots p_{2k-2}p_{2k-1}^2\leq x^{\rho}\,\, \textnormal{for all}\,\, k\geq 1\},    
\end{align*}
which is the support of the upper bound semilinear weights, when $\rho=\frac{2}{5}(1-4\theta)-\varepsilon$, $\theta\leq \frac{1}{40}$, and all the other parameters are as before. This observation will be used in the proof of Theorem \ref{theo_sievebombieri}. This exponent is also optimal, as is seen by taking $\rho=\frac{2}{5}(1-4\theta)+2\varepsilon$ and $D=x^{\frac{2}{5}(1-4\theta)}$, $p_1p_2\in \mathcal{D}^{+, \textnormal{SEM}}$, $p_1,p_2\sim \frac{1}{2}x^{\frac{1}{5}(1-4\theta)+\varepsilon}$.
\end{remark}

\textbf{Proof of Lemma \ref{le1}.} The proof resembles some arguments related to Harman's sieve \cite[Chapter 3]{harman-sieves}. Let $d=p_1\cdots p_r \in \mathcal{D}^{-,\textnormal{SEM}}$. The claim is that the set $\{p_1,\ldots,p_r\}$ can be partitioned into two subsets $S_1$ and $S_2$ in such a way that the products $P_1$ and $P_2$ of the elements of $S_1$ and $S_2$ satisfy $P_1\leq D$, $P_1P_2^2\leq \frac{x^{1-4\theta-2\varepsilon^2}}{D}$, and additionally $P_1\geq x^{0.1}$ or $P_2=1$. Note that for $r=1$ one can take $S_1=\{p_1\}$ and $S_2=\emptyset$. Assume then that $r\geq 2$. If $p_1\cdots p_r\leq D$, we may take $S_1=\{p_1,\ldots,p_r\}$, $S_2=\emptyset$. Indeed, then $P_1\leq D$, $P_2=1$ and $P_1P_2^2\leq D\leq \frac{x^{1-4\theta-2\varepsilon^2}}{D}$. Now we may assume that $p_1\cdots p_r>D$. Since $p_1\leq D$, we can select the largest $j$ for which $p_1\cdots p_j\leq D$. We have $j\leq r-1$ and $p_{j+1}\leq p_2\leq x^{\frac{\rho}{3}}$, so
\begin{align*}
p_1\cdots p_j=\frac{p_1\cdots p_{j+1}}{p_{j+1}}\geq \frac{D}{x^{\frac{\rho}{3}}}.
\end{align*}
We claim that the choice $S_1=\{p_1,\ldots p_j\}$, $S_2=\{p_{j+1},\ldots,p_r\}$ works. First of all, we have $P_1\geq \frac{D}{x^{\frac{\rho}{3}}}\geq x^{0.1}$. Supposing that the claim does not hold for $S_1$ and $S_2$, we have $(P_1P_2)^2>P_1\frac{x^{1-4\theta-2\varepsilon^2}}{D}$. Using $P_1P_2\leq x^{\rho}$ and $P_1\geq \frac{D}{x^{\frac{\rho}{3}}}$, this yields $x^{2\rho}> x^{1-4\theta-\frac{\rho}{3}-2\varepsilon^2}$, from which we solve $\rho>\frac{3}{7}(1-4\theta)-\frac{6}{7}\varepsilon^2$, a contradiction to our choice of $\rho$.\qedd

\section{Verifying the Hypothesis}\label{Sec: hypotheses}

\subsection{Splitting variables}\label{sub: splitting}

Based on Section \ref{Sec: decomposition}, the proof of Hypothesis \ref{h1} for the sequence $(\omega_n)_{n\sim x}$ and linear function $L(n)$ defined in that section has been reduced to showing that
\begin{align}
\sum_{\substack{d\leq x^{\rho_2}\\(d,QW)=1}}\hspace{-0.1cm}\lambda_{d}^{-,\textnormal{SEM}}\bigg(\hspace{-0.1cm}\sum_{\substack{n\sim x\\L(n)\in \mathbb{P}\\L(n)\equiv 1\hspace{-0.1cm} \pmod d}}e(\xi n)-\frac{1}{\varphi(d)}\frac{QW}{\varphi(QW)}\sum_{n\sim x}\frac{e(\xi n)}{\log(QWn)}\bigg)\quad \text{and}\label{eq41}
\end{align}
\begin{align}
\sum_{\substack{d\leq x^{\rho_1}\\(d,QW)=1}}\hspace{-0.1cm}\lambda_d^{+,\textnormal{LIN}}\hspace{-0.1cm}\sum_{\substack{\ell \leq x^{1-\varepsilon}\\(\ell,QW)=\delta\\(\ell,d)=1}}g(\ell)\bigg(\hspace{-0.1cm}\sum_{\substack{n\sim x\\L(n)=\ell p+1\\L(n)\equiv 0\hspace{-0.1cm} \pmod d}}e(\xi n)-\frac{1}{\varphi(d)}\frac{QW}{\varphi(\frac{QW}{\delta})}\sum_{n \sim x}\frac{e(\xi n)}{\ell \log \frac{QWn}{\ell}}\bigg)\label{eq93}
\end{align}
are $\ll x(\log x)^{-200}$, where $\delta=(Wc_0+b-1,QW)$, $(g(\ell))_{\ell\geq 1}$ is a sequence of convolution type (with parameter $\sigma$), the sieve weights $\lambda_d^{+,\textnormal{LIN}}, \lambda_d^{-,\textnormal{SEM}}$ have respective sifting parameters $z_1\leq x^{\frac{1}{5}+\varepsilon}$, $z_2\leq x^{\frac{1}{3+\frac{\varepsilon}{2}}}$, and $\rho_1,\rho_2$, $\sigma$ are as in \eqref{eq87}, and $\xi$ is subject to \eqref{eq88}. It would actually suffice to replace $\ell\leq x^{1-\varepsilon}$ by $\ell\leq x^{0.9+\varepsilon}$ above, but this would not simplify the argument.\\

As mentioned in Section \ref{Sec: sieveweight}, we wish to split the sum over $d$ into a double sum. This is enabled by Lemmas \ref{le10} and \ref{le1}. If $D$ is as in Lemma \ref{le1} with $0\leq \theta\leq \frac{1}{30}$, we may write
\begin{align}\label{eq89}
|\lambda_d^{-,\textnormal{SEM}}|&\leq \min_{D}\sum_{\substack{d=d_1d_2\\d_1\leq D\\d_1d_2^2\leq \frac{x^{1-4\theta-2\varepsilon^2}}{D}\\(d_1,d_2)=1\\d_1\geq x^{0.1}\,\textnormal{or}\,d_2=1}}1\leq \left(\frac{\log x}{\log 2}\right)^2\min_{D}\max_{\Delta_1,\Delta_2}\sum_{\substack{d=d_1d_2\\d_1\sim \Delta_1\\d_2\sim \Delta_2\\(d_1,d_2)=1}}1,
\end{align}
where the maximum and minimum are over those $\Delta_1, \Delta_2\geq 1$ and $D\geq 1$ that satisfy
\begin{align}\label{eq90}\begin{split}
&D\in [x^{\frac{1}{3}-2\theta-2\varepsilon^2},x^{\rho_2}],\,\,\Delta_1\leq D,\,\, \Delta_1\Delta_2^2\leq \frac{x^{1-4\theta-2\varepsilon^2}}{D},\,\, \Delta_1\Delta_2\leq x^{\rho_2},\\
&\text{and either}\,\, \Delta_1\geq x^{0.1}\,\,\text{or}\,\, \Delta_2=1.
\end{split}
\end{align}
By Lemma \ref{le10}, formula \eqref{eq89} continues to hold with $\lambda_{d}^{-,\textnormal{SEM}}$ replaced with $\lambda_{d}^{+,\textnormal{LIN}}$ and \eqref{eq90} replaced with
\begin{align}\label{eq91}\begin{split}
&D\in [x^{\frac{1}{5}},x^{\rho_1}],\,\,\Delta_1\leq D,\,\, \Delta_1\Delta_2^2\leq \frac{x^{1-4\theta-2\varepsilon^2}}{D},\,\, \Delta_1\Delta_2\leq x^{\rho_1},\\ 
&\text{and either}\,\, \Delta_1\geq x^{0.1}\,\,\text{or}\,\, \Delta_2=1.
\end{split}
\end{align}
We take $\theta=0$ in this section, but in Section \ref{Sec: fractional parts} we will employ the same formulas with $\theta>0$.
As a conclusion,  we see that \eqref{eq41} and \eqref{eq93} are bounded by $(\frac{\log x}{\log 2})^2$ times 
\begin{align}
&\sum_{\substack{d_1\sim \Delta_1\\(d_1,QW)=1}}\sum_{\substack{d_2\sim \Delta_2\\(d_2,QW)=1\\(d_1,d_2)=1}}\hspace{-0.1cm}\bigg|\hspace{-0.1cm}\sum_{\substack{n\sim x\\L(n)\in \mathbb{P}\\L(n)\equiv 1\hspace{-0.1cm} \pmod{d_1d_2}}}\hspace{-0.1cm}e(\xi n)-\frac{QW}{\varphi(d_1d_2)\varphi(QW)}\sum_{n\sim x}\frac{e(\xi n)}{\log(QWn)}\bigg|\quad \text{and}\label{eq15}\\
&\sum_{\substack{d_1\sim \Delta_1\\(d_1,QW)=1}}\sum_{\substack{d_2\sim \Delta_2\\(d_2,QW)=1\\(d_1,d_2)=1}}\bigg|\sum_{\substack{\ell \leq x^{1-\varepsilon}\\(\ell,QW)=\delta\\(\ell,d_1d_2)=1}}g(\ell)\bigg(\hspace{-0.1cm}\sum_{\substack{n\sim x\\L(n)=\ell p+1\\L(n)\equiv 0\hspace{-0.1cm} \pmod{d_1d_2}}}\hspace{-0.1cm}e(\xi n)-\frac{QW}{\varphi(d_1d_2)\varphi(\frac{QW}{\delta})}\sum_{n \sim x}\frac{e(\xi n)}{\ell \log \frac{QWn}{\ell}}\bigg)\bigg|\label{eq92},
\end{align}
respectively, where $\Delta_1$ and $\Delta_2$ are any numbers constrained by \eqref{eq90} or \eqref{eq91}, depending on whether we consider \eqref{eq15} or \eqref{eq92}. At this point, it is also natural to split into two cases depending on whether $\xi$ lies on a major arc or minor arc (that is, whether $q\mid Q$ or $\frac{q}{(q,Q^2)}\geq (\log x)^{A}$ holds in \eqref{eq88}).

\subsection{Major arcs for the semilinear sieve} \label{sub: maj sem}

We first assume the major arc  condition $q\mid Q$ in the definition of $\xi$ in \eqref{eq88}. By partial summation, \eqref{eq41} becomes
\begin{align*}
=\int_{x}^{2x}e(\pm\|\xi\| t)\,d\bigg\{\sum_{\substack{d\leq x^{\rho_2}\\(d,QW)=1}}\lambda_d^{-,\textnormal{SEM}} \bigg(\sum_{\substack{x\leq n\leq t\\L(n)\in \mathbb{P}\\L(n)\equiv 1\hspace{-0.1cm} \pmod d}}1-\frac{QW}{\varphi(QW)}\frac{1}{\varphi(d)}\sum_{\substack{x\leq n\leq t}}\frac{1}{\log(QW n)}\bigg)\bigg\}.
\end{align*}
Naming the function inside $d\{\ldots\}$ as $G(t)$, partial integration tells that the previous expression is
\begin{align}\label{eq38}
=G(2x)e(\pm2\|\xi\|x)\mp2\pi i \|\xi\|\int_{x}^{2x}e(\pm\|\xi\| t)G(t)dt\ll (1+\|\xi\|x)\max_{x\leq t\leq 2x}|G(t)|.
\end{align}
Since $\frac{1}{\log(QWn)}=\frac{1}{QW}\int_{QWn}^{QW(n+1)}\frac{du}{\log u}+O(\frac{1}{n})$, putting $c_1=Wc_0+b$ we have
\begin{align*}
G(t)&\leq\sum_{\substack{d\leq x^{\rho_2}\\(d,QW)=1}}|\lambda_d^{-,\textnormal{SEM}}|\bigg|\sum_{\substack{QWx\leq p\leq QWt\\p\equiv c_1\hspace{-0.1cm} \pmod{QW}\\p\equiv 1 \hspace{-0.1cm} \pmod d}}1-\frac{1}{\varphi(QWd)}\int_{QWx}^{QWt}\frac{du}{\log u}\bigg|+O(x^{\frac{1}{2}})\\
&\leq \sum_{\substack{d\leq x^{\rho_2}\\(d,QW)=1}}\max_{(r,QWd)=1}\left|\pi(QWt;QWd,r)-\frac{1}{\varphi(QWd)}\textnormal{Li}(QWt)\right|\\
&\quad +\sum_{\substack{d\leq x^{\rho_2}\\(d,QW)=1}}\max_{(r,QWd)=1}\left|\pi(QWx;QWd,r)-\frac{1}{\varphi(QWd)}\textnormal{Li}(QWx)\right|+O(x^{\frac{1}{2}})\\
&\ll \frac{x}{(\log x)^{1000B}}
\end{align*}
by the Bombieri-Vinogradov theorem \cite[Theorem 17.1]{iwaniec-kowalski}. As $\xi$ is on a major arc, by \eqref{eq88} we have $\|\xi\|\leq \frac{2(\log x)^{102B}}{x}$, so \eqref{eq38} is $\ll x(\log x)^{-1000}$. Therefore, the major arc case for the semilinear sieve has been dealt with.

\subsection{Major arcs for the linear sieve}

Again assume $q\mid Q$ in \eqref{eq88}. After applying partial summation, \eqref{eq93} takes the form
\begin{align*}
\int_{x}^{2x}e(\pm \|\xi\| t)\,d\bigg\{\hspace{-0.1cm}\sum_{\substack{d\leq x^{\rho_1}\\(d,QW)=1}}\lambda_d^{+,\textnormal{LIN}}\bigg(\hspace{-0.1cm}\sum_{\substack{x\leq n\leq t\\L(n)=\ell p+1\\L(n)\equiv 0\hspace{-0.1cm} \pmod d\\ \ell\leq x^{1-\varepsilon}\\(\ell,QW)=\delta\\(\ell,d)=1}}\hspace{-0.1cm}g(\ell)-\frac{QW}{\varphi(d)\varphi(\frac{QW}{\delta})}\sum_{\substack{x\leq n\leq t\\ \ell\leq x^{1-\varepsilon}\\(\ell,QW)=\delta\\(\ell,d)=1}}\frac{g(\ell)}{\ell \log \frac{QWn}{\ell}}\bigg)\bigg\},
\end{align*}
so we want this to be $\ll x(\log x)^{-202}$. Proceeding as in Subsection \ref{sub: maj sem}, it suffices to prove for that $t\in [x,2x]$
\begin{align*}
\sum_{\substack{d\leq x^{\rho_1}\\(d,QW)=1}}\bigg|\sum_{\substack{x\leq n\leq t\\L(n)=\ell p+1\\L(n)\equiv 0\hspace{-0.1cm} \pmod d\\ \ell\leq x^{1-\varepsilon}}}g(\ell)1_{(\ell,QW)=\delta,\,\,(\ell,d)=1}-\frac{QW}{\varphi(d)\varphi(\frac{QW}{\delta})}\sum_{\substack{x\leq n\leq t\\ \ell\leq x^{1-\varepsilon}\\(\ell,QW)=\delta\\(\ell,d)=1}}\frac{g(\ell)}{\ell \log \frac{QWn}{\ell}}\bigg|
\end{align*}
is $\ll x(\log x)^{-1000B}$.\\

We start by analyzing the second sum inside the absolute values in the previous expression. Since $QW\ll  (\log x)^{B+1}$ and $\ell\leq x^{1-\varepsilon}$,  a change of variables and the prime number theorem give
\begin{align*}
\frac{QW}{\varphi(\frac{QW}{\delta})}\sum_{x\leq n\leq t}\frac{1}{\ell \log \frac{QWn}{\ell}}&=\frac{QW}{\varphi(\frac{QW}{\delta})}\int_{x}^{t}\frac{du}{\ell \log \frac{QWu}{\ell}}+O(QW)\\
&=\frac{1}{\varphi(\frac{QW}{\delta})}\int_{\frac{QWx}{\ell}}^{\frac{QWt}{\ell}}\frac{du}{\log u}+O(QW)\\
&=\frac{1}{\varphi(\frac{QW}{\delta})}\sum_{\substack{QWx\leq \ell p\leq QWt}}1+O\left(\frac{x}{\ell}(\log x)^{-3000B}\right).
\end{align*}
The error term remains still $\ll x(\log x)^{-2000B}$ after multiplying it by $\frac{|g(\ell)|}{\varphi(d)}$ and summing over $d\leq x^{\rho_1}$,$\ell\leq x^{1-\varepsilon}$. Hence, what we wish to show is that
\begin{align}\label{eq69}
\sum_{\substack{d\leq x^{\rho_1}\\(d,QW)=1}}\bigg|\sum_{\substack{QWx\leq \ell p\leq QWt\\\ell p\equiv -1\hspace{-0.1cm} \pmod{d}\\\ell p\equiv c_1-1\hspace{-0.1cm} \pmod{QW}\\\ell\leq x^{1-\varepsilon}\\(\ell,QW)=\delta\\(\ell,d)=1}}g(\ell)-\frac{1}{\varphi(\frac{QWd}{\delta})}\sum_{\substack{QWx\leq \ell p\leq QWt\\\ell\leq x^{1-\varepsilon}\\(\ell,QW)=\delta\\(\ell,d)=1}}g(\ell)\bigg|
\end{align}
is $\ll \frac{x}{(\log x)^{1000B}}$ for $t\in [x,2x]$ and $c_1=Wc_0+b$. Since $(\ell,QW)=\delta$, $(\ell,d)=1$ and $(d,\delta)=1$, the congruences $\ell p\equiv -1\hspace{-0.1cm} \pmod d$, $\ell p\equiv c_1-1\hspace{-0.1cm} \pmod{QW}$ can be rewritten as $\ell' p\equiv -\delta^{-1}\hspace{-0.1cm} \pmod{d}$, $\ell' p\equiv \frac{c_1-1}{\delta} \hspace{-0.1cm} \pmod{\frac{QW}{\delta}}$ with $\ell'=\frac{\ell}{\delta}$. By the Chinese remainder theorem, these congruences are equivalent to $\ell' p\equiv c \hspace{-0.1cm} \pmod{\frac{QWd}{\delta}}$ for some $c$ depending on $Q,W, d$ and $\delta$ and coprime to $\frac{QWd}{\delta}$. Concerning the second sum inside absolute values in \eqref{eq69}, we wish to add the constraint $(\ell' p,\frac{QWd}{\delta})=1$ to that summation (where again $\ell'=\frac{\ell}{\delta}$). We know that $(\ell',\frac{QW}{\delta})=(\ell',d)=1$, and clearly $p\geq x^{\varepsilon}$ in \eqref{eq69}, so $(p,QW)=1$. Therefore, we have shown that we may insert the constraint $(\ell' p,QWd)=1$ if the case $p\mid d$ has a small enough contribution to the aforementioned sum. That case contributes at most
\begin{align*}
\sum_{\substack{p\mid d\\p\geq x^{\varepsilon}}}\sum_{\ell\leq \frac{2QWx}{p}}|g(\ell)|\ll_{\varepsilon}x^{1-\frac{\varepsilon}{2}}, 
\end{align*}
which is $\ll x^{1-\varepsilon^2}$ when multiplied by $\frac{1}{\varphi(\frac{QWd}{\delta})}$ and summed over $d\leq x^{\rho_1}$. Summarizing, our aim has been reduced to showing that
\begin{align}\label{eq70}
\sum_{\substack{d\leq x^{\rho_1}\\(d,QW)=1}}\max_{(c,\frac{QWd}{\delta})=1}\bigg|\sum_{\substack{\frac{QWx}{\delta}\leq \ell' p\leq \frac{QWt}{\delta}\\\ell' p\equiv c \hspace{-0.1cm} \pmod{\frac{QWd}{\delta}}\\\ell'\leq x^{1-\varepsilon}/\delta}}g(\delta\ell')-\frac{1}{\varphi(\frac{QWd}{\delta})}\sum_{\substack{\frac{QWx}{\delta}\leq \ell' p\leq \frac{QWt}{\delta}\\(\ell' p,\frac{QWd}{\delta})=1\\\ell'\leq x^{1-\varepsilon}/\delta}}g(\delta \ell')\bigg|
\end{align}
is $\ll \frac{x}{(\log x)^{1000B}}$ for $t\in [x,2x]$.\\

 To obtain this estimate, we apply \cite[Theorem 17.4]{iwaniec-kowalski} to the sequences $(\alpha_{\ell'})_{\ell' \leq x^{1-\varepsilon}/\delta}=(g(\delta \ell'))_{\ell'\leq x^{1-\varepsilon}/\delta}$ and $(\beta_k)_{k\geq 1}=(1_{\mathbb{P}}(k))_{k\geq 1}$ -- that theorem is applicable since the sequence $(1_{\mathbb{P}}(k))_{k\geq 1}$ is well-distributed in the sense of formula (17.13) of \cite{iwaniec-kowalski} (with $\Delta=(\log x)^{-20000B}$ there) by the Siegel-Walfisz theorem. Now, since in \eqref{eq70} we have $\ell'\geq x^{\frac{\varepsilon}{2}}$, $p\geq x^{\varepsilon}$, $\rho_1<\frac{1}{2}$ and $|\alpha_{\ell'}|\leq \tau(\ell')^2\log \ell'$, the claimed Bombieri-Vinogradov type estimate follows immediately from the theorem cited above.

\subsection{Minor arcs for the semilinear sieve} \label{sub: min sem}

We assume then that $\xi$ is on a minor arc, meaning that $\frac{q}{(q,Q^2)}\geq (\log x)^A$ in \eqref{eq88}. We study the sum \eqref{eq15}. Using partial summation, we see that
\begin{align*}
\sum_{n\sim x}\frac{e(\xi n)}{\log(QWn)}\ll \max_{x\leq t\leq 2x}\left|\sum_{x\leq n\leq t}e\left(\xi n\right)\right|\ll \frac{1}{\|\xi\|}.
\end{align*} 
We have $(q,QW)\leq W(q,Q)\leq \frac{Wq}{(\log x)^{A}}<q$, so $q\nmid QW$. Taking this and \eqref{eq88} into account, $\|\xi\|\geq \frac{1}{q}-\frac{2(\log x)^{102B}}{q x}\geq \frac{1}{2q}$, so the second expression inside absolute values in \eqref{eq15} is $\ll \frac{q}{\varphi(d)}\ll \frac{x}{(\log x)^{99B}\varphi(d)}$. Hence it contributes $\ll x(\log x)^{-98B}$ when summing over $d$.\\

 When it comes to the first expression inside absolute values in \eqref{eq15}, it equals
\begin{align*}
\sum_{\substack{n\sim x\\L(n)\in \mathbb{P}\\L(n)\equiv 1\hspace{-0.1cm} \pmod{d_1d_2}}}e(\xi n)=e\left(\frac{-\xi c_1}{QW}\right)\sum_{\substack{p\sim QWx\\p\equiv c_1\hspace{-0.1cm} \pmod{QW}\\p\equiv 1\hspace{-0.1cm} \pmod{d_1d_2}}}e\left(\frac{\xi}{QW}p\right)+O(QW),
\end{align*}
where the error $O(QW)$ remains $\ll x^{\frac{1}{2}}$ when summed over $d\leq x^{\rho_2}$. With partial summation, we may bound the sum on the right-hand side by
\begin{align*}
\bigg|\sum_{\substack{n\sim QWx\\n\equiv c_1\hspace{-0.1cm} \pmod{QW}\\n\equiv 1\hspace{-0.1cm} \pmod{d_1d_2}}}\hspace{-0.1cm}\Lambda(n)e\left(\frac{\xi}{QW}n\right)\bigg|+\int_{QWx}^{2QWx}\hspace{-0.1cm}\sum_{\substack{QWx\leq n\leq t\\n\equiv c_1\hspace{-0.1cm} \pmod{QW}\\n\equiv 1\hspace{-0.1cm} \pmod{d_1d_2}}}\hspace{-0.1cm}\Lambda(n)e\left(\frac{\xi}{QW}n\right)\, \frac{dt}{t\log^2 t}+O(x^{\frac{1}{2}+\varepsilon}),
\end{align*}
the error coming from the values of $n$ that are prime powers, and the error being $\ll x^{1-\varepsilon^2}$ after summing over $d\leq x^{\rho_2}$. This means that it suffices to prove
\begin{align}\label{eq94}
\sum_{\substack{d_1\sim \Delta_1\\(d_1,QW)=1}}\sum_{\substack{d_2\sim \Delta_2\\(d_2,QW)=1\\(d_1,d_2)=1}}\bigg|\sum_{\substack{QWx\leq n\leq t\\n\equiv c_1\hspace{-0.1cm}\pmod{QW}\\n\equiv 1 \hspace{-0.1cm}\pmod{d_1d_2}}}\Lambda(n)e\left(\frac{\xi}{QW}n\right)\bigg|\ll \frac{x}{(\log x)^{1000}}
\end{align}
uniformly for $t\in [QWx,2QWx]$. We may now apply Vaughan's identity (in the form of \cite[Proposition 13.4]{iwaniec-kowalski} with $y=z=(QWx)^{\frac{1}{3}}$ there), which transforms the sum inside absolute values in \eqref{eq94} (up to error $O(x^{\frac{1}{3}+\varepsilon})$) into a sum of $\ll (\log x)^{10}$ type I and type II sums of the form
\begin{align*}
\widetilde{R}_{d_1d_2}^{\text{I}}(t)=\hspace{-0.2cm}\sum_{\substack{QWx\leq mn\leq t\\mn\equiv c_1 \hspace{-0.1cm} \pmod{QW}\\ mn\equiv 1\hspace{-0.1cm} \pmod{d_1d_2}\\m\asymp M}}\alpha_m e\left(\frac{\xi m n}{QW}\right)\,\, \text{and}\,\, \widetilde{R}_{d_1d_2}^{\text{II}}(t)=\hspace{-0.2cm}\sum_{\substack{QWx\leq mn\leq t\\mn\equiv c_1 \hspace{-0.1cm} \pmod{QW}\\ mn\equiv 1\hspace{-0.1cm} \pmod{d_1d_2}\\m\asymp M}}\alpha_m \beta_n e\left(\frac{\xi m n}{QW}\right),
\end{align*}
with $|\alpha_m|,|\beta_m|\leq \tau(m)^2\log m$ some complex numbers and $M\leq (2QWx)^{\frac{1}{3}}$ in the case of $\widetilde{R}^{\text{I}}_{d_1d_2}(t)$, while $M\in [(QWx)^{\frac{1}{3}},(2QWx)^{\frac{2}{3}}]$ in the case of $\widetilde{R}^{\text{II}}_{d_1d_2}(t)$. Moreover, we may assume in the latter case that $M\in [(QWx)^{\frac{1}{2}},(2QWx)^{\frac{2}{3}}]$ by flipping the roles of the variables if necessary. We may replace the type I and type II sum with the (possibly larger) sums
\begin{align}\begin{split}\label{eq100}
R_{d_1d_2}^{\text{I}}(t)&=\max_{(c,d_1d_2QW)=1}\bigg|\sum_{\substack{QWx\leq mn\leq t\\mn\equiv c \hspace{-0.1cm} \pmod{d_1d_2QW}\\m\asymp M}}\alpha_m e\left(\frac{\xi}{QW}m n\right)\bigg|\quad \text{and}\\
\quad R_{d_1d_2}^{\text{II}}(t)&=\max_{(c,d_1QW)=1}\bigg|\sum_{\substack{QWx\leq mn\leq t\\ mn\equiv c\hspace{-0.1cm} \pmod{d_1QW}\\mn\equiv 1 \hspace{-0.1cm}\pmod{d_2}\\m\asymp M}}\alpha_m \beta_n e\left(\frac{\xi}{QW} m n\right)\bigg|.
\end{split}
\end{align}
We are now in a position to apply the Bombieri-Vinogradov lemmas \ref{le8} and \ref{le9}. Note that, by \eqref{eq88}, we either have $|\xi-\frac{QWa}{q}|\leq \frac{1}{(QWq)^2}$ or $q>\frac{x}{2(\log x)^{102B}(QW)^2}$. If the latter happens, we have $|e(\frac{\xi}{QW} mn)-e(\frac{a}{q}mn)|\leq |\frac{\xi}{QW}-\frac{a}{q}|mn\leq \frac{8(QW)^3(\log x)^{204B}}{x}$ for $mn\leq 2QWx$. This implies that $e(\frac{\xi}{QW} mn)$ can be replaced by $e(\frac{a}{q}mn)$ in the type I and II sums. In conclusion, we can assume in any case that $|\xi-\frac{QWa}{q}|\leq \frac{1}{(QWq)^2}$.\\

The type I Bombieri-Vinogradov sums cause no problems, as Lemma \ref{le8} with the choices $R=1$, $N=QWx$, $v=QW$, $M= x^{\frac{1}{3}+\varepsilon}$, $\rho\leq \frac{1}{2}-\varepsilon$ tells at once that
\begin{align*}
\sum_{\substack{d_1\sim \Delta_1\\(d_1,QW)=1}}\sum_{\substack{d_2\sim \Delta_2\\(d_2,QW)=1\\(d_1,d_2)=1}}R_{d_1d_2}^{\text{I}}(t)\ll \frac{x}{(\log x)^{\frac{A}{10}}},
\end{align*}
since $\frac{q}{(q,(QW)^2)}\geq W^{-2}(\log x)^{A}$ and $\Delta_1\Delta_2\leq x^{\rho_2}$.\\

We know that $(QWx)^{\frac{1}{2}}\leq M\leq (2QWx)^{\frac{2}{3}}$ in the sum $R_{d_1d_2}^{\text{II}}(t)$. We divide the analysis of this sum into three cases.\\

\textbf{Case 1.} Assume that $M\geq x^{1-\rho_2-\varepsilon^2}$, $\Delta_1\geq (\log x)^{\frac{A}{10}}$. Take $D=\frac{x^{1-\varepsilon^2}}{M}$. We know that $x^{\frac{1}{3}-\varepsilon^2}(\log x)^{-B}\leq D\leq x^{\rho_2}$ by the bound on $M$. In view of \eqref{eq90} with $\theta=0$, this means in particular that $\Delta_1\leq \frac{x^{1-\varepsilon^2}}{M}$ and $\Delta_1\Delta_2^2\leq \frac{x^{1-2\varepsilon^2}}{D}= Mx^{-\varepsilon^2}$. Now we apply Lemma \ref{le9} (in the case of $F_1$) with $R=1$, $N=QWx$, $v=QW$, $\rho=\rho_2\leq \frac{3}{7}-\varepsilon$ to deduce that
\begin{align*}
&\sum_{\substack{d_1\sim \Delta_1\\(d_1,QW)=1}}\sum_{\substack{d_2\sim \Delta_2\\(d_2,QW)=1\\(d_1,d_2)=1}}R_{d_1d_2}^{\text{II}}(t)\\
&\ll x\left(\left(\frac{1}{\Delta_1}+\frac{W^2}{(\log x)^{A}}+(\log x)^{-99B}(QW)^2\right)^{\frac{1}{8}}+\left(\frac{\Delta_1 M}{x}+\Delta_1\Delta_2^2\frac{QW}{M}\right)^{\frac{1}{2}}\right)(\log x)^{1000},
\end{align*}
which is $\ll \frac{x}{(\log x)^{\frac{A}{100}}}$ for $A$ large enough by the lower bound on $\Delta_1$.\\

\textbf{Case 2.} Assume then that $M\geq x^{1-\rho_2-\varepsilon^2}$, $\Delta_1< (\log x)^{\frac{A}{10}}$. Since $\Delta_1<x^{0.1}$, we know that $\Delta_2=1$, so applying Lemma \ref{le9} (in the case of $F_2$) we obtain, for $A$ large enough,
\begin{align*}
\sum_{\substack{d_1\sim \Delta_1\\(d_1,QW)=1}}\sum_{\substack{d_2\sim \Delta_2\\(d_,QW)=1\\(d_1,d_2)=1}}R_{d_1d_2}^{\text{II}}(t)\ll x(\log x)^{\frac{A}{5}}\left(\frac{W}{(\log x)^{\frac{A}{2}}}+\frac{QW}{M^{\frac{1}{2}}}+\frac{(QW)^2M}{x}+\frac{(QW)^{\frac{1}{2}}}{(\log x)^{\frac{99B}{2}}}\right)^{\frac{1}{2}},
\end{align*}
and this is again $\ll \frac{x}{(\log x)^{\frac{A}{100}}}$ for $A$ large.\\

\textbf{Case 3.} Lastly, assume that $M<x^{1-\rho_2-\varepsilon^2}$. Then we estimate \eqref{eq41} instead of \eqref{eq15}. This amounts to just replacing $d_1\sim \Delta_1$, $d_2\sim \Delta_2$ with $d_1\leq x^{\rho_2}$, $d_2=1$ throughout this subsection. We have $x^{\rho_2}\leq \frac{x^{1-\varepsilon^2}}{M}$ and $x^{\rho_2}\leq Mx^{-\varepsilon^2}$, so we can bound the type II sums in the same way as for $M\geq x^{1-\rho_2-\varepsilon^2}$ (considering again the cases $\Delta_1\geq (\log x)^{\frac{A}{10}}$ and $\Delta_1< (\log x)^{\frac{A}{10}}$ separately), so also Case 3 contributes $\ll \frac{x}{(\log x)^{\frac{A}{100}}}$.\\

 Consequently, we have shown that the contribution of the minor arcs for the semilinear sieve is small enough.

\subsection{Minor arcs for the linear sieve}

We assume again $\frac{q}{(q,Q^2)}\geq (\log x)^{A}$.  We first look at the second expression inside absolute values in \eqref{eq92}. We have by partial summation
\begin{align*}
\sum_{n\sim x}\frac{e(\xi n)}{\ell \log\frac{QWn}{\ell}}\ll \frac{1}{\ell\|\xi\|}
\end{align*}
for $\ell\leq x^{1-\varepsilon}$ just as in Subsection \ref{sub: min sem}. We showed earlier that $\frac{1}{\|\xi\|}\ll \frac{x}{(\log x)^{99B}}$ when $\frac{q}{(q,Q^2)}\geq (\log x)^{A}$, so the second expression inside absolute values in \eqref{eq92} is $\ll \frac{x}{\ell\varphi(d)(\log x)^{98B}}$, which is $\ll x(\log x)^{-97B}$ after summing over $d\leq x^{\rho_1}$ and over $\ell\leq x^{1-\varepsilon}$ weighted by $|g(\ell)|$.\\

We may write the first expression inside absolute values in \eqref{eq92} as
\begin{align}\label{eq73}
e\left(\frac{-(c_1-1)\xi}{QW}\right)\hspace{-0.1cm}\sum_{\substack{\ell p\sim QWx\\\ell p\equiv c_1-1\hspace{-0.1cm} \pmod{QW}\\\ell p\equiv -1\hspace{-0.1cm} \pmod d\\\ell\leq x^{1-\varepsilon}}}\hspace{-0.1cm}g(\ell)e\left(\frac{\xi}{QW}\ell p\right)+O(QW),
\end{align}
and the error $O(QW)$ is $\ll x^{\frac{1}{2}}$ after summing over $d\leq x^{\rho_1}$. We have ignored the conditions $(\ell,QW)=\delta, (\ell,d)=1$ above, since if either of them fails, $\ell p\equiv c_1-1 \hspace{-0.1cm} \pmod{QW}$, $\ell p\equiv -1\pmod d$ is impossible.\\

Crucially, our assumption is that the sequence $(g(\ell))_{\ell\geq 1}$ is of convolution type, so the sum in \eqref{eq73} can be rewritten as
\begin{align*}
\sum_{\substack{km p\sim QWx\\k m p\equiv c_1-1\hspace{-0.1cm} \pmod{QW}\\k m p\equiv -1\hspace{-0.1cm} \pmod d\\ k m\leq x^{1-\varepsilon}}}\alpha_k\beta_m e\left(\frac{\xi}{QW} k m p\right),
\end{align*}
where $(\alpha_k)$ is supported on $x^{\frac{1}{\sigma}}\leq k\leq (Qx)^{1-\frac{1}{\sigma}}$ for $\sigma=3+\varepsilon$. Putting
\begin{align*}
\beta^{*}_r=\sum_{r=mp}\beta_m
\end{align*}
and splitting the previous sum dyadically, it becomes $\ll \log x$ sums of the form
\begin{align*}
\sum_{\substack{k r\sim QWx\\k r\equiv c_1-1\hspace{-0.1cm} \pmod{QW}\\k r\equiv -1\hspace{-0.1cm} \pmod d\\k\asymp M}}\alpha_k\beta^{*}_r e\left(\frac{\xi}{QW} k r\right),
\end{align*}
where $x^{\frac{1}{\sigma}}\leq M\leq (Qx)^{1-\frac{1}{\sigma}}$, and by changing the roles of the variables, we may further assume that $(QWx)^{\frac{1}{2}}\leq M\leq Qx^{1-\frac{1}{\sigma}}$. Now our bilinear sums are exactly of the same form as in \eqref{eq100} (but with different $M$). Furthermore, we may assume that $|\xi-\frac{QWa}{q}|\leq \frac{1}{(QWq)^2}$ for the same reason as in Subsection \ref{sub: min sem}. If $M\geq x^{1-\rho_1-\varepsilon^2}$, denoting $D=\frac{x^{1-\varepsilon^2}}{M}\in [x^{\frac{1}{5}},x^{\rho_1}]$, we again see that $\Delta_1\Delta_2^2\leq Mx^{-\varepsilon^2}$ in \eqref{eq91} (with $\theta=0$). Therefore, we may apply the very same estimates as in the Cases 1 and 2 of Subsection \ref{sub: min sem}. If $M<x^{1-\rho_1-\varepsilon^2}$, we can apply precisely the same argument as in Case 3 of the previous subsection, since $x^{\rho_1}\leq \frac{x^{1-\varepsilon^2}}{M}$ and $x^{\rho_1}\leq Mx^{-\varepsilon^2}$. Summarizing, we have showed that the minor arcs for the linear sieve contribute $\ll x(\log x)^{-\frac{A}{100}}$, which is small enough for large $A$.\\

 We have now concluded the proof of Theorem \ref{theo_goldbach}, in view of Theorem \ref{t2} and Proposition \ref{prop2}.\qedd\\

\textbf{Proof of Theorem \ref{theo_sievebombieri}:} We take $Q=W=1$ and $L(n)=n$ in \eqref{eq41} and replace $L(n)\equiv 1 \pmod d$ by $L(n)\equiv b \pmod d$ (with $b\neq 0$ an arbitrary integer) there and note that the proof that \eqref{eq41} is $\ll_C x(\log x)^{-C}$ is verbatim the same as the minor arc argument for the semilinear sieve in this section, provided that $\xi$ is any real number with $|\xi-\frac{a}{q}|\leq \frac{1}{q^2}$ for some coprime $a$ and $q\in [(\log x)^{1000C},x(\log x)^{-1000C}]$. This proves Theorem \ref{theo_sievebombieri} in the case of lower bound sieve weights. The case of upper bound sieve weights follows very similarly by replacing $\lambda_{d}^{-,\textnormal{SEM}}$ with $\lambda_{d}^{+,\textnormal{SEM}}$ and making use of a remark after Lemma \ref{le1} (which is where the value $\rho_{+}=\frac{2}{5}-\varepsilon$ comes from).\qedd

\section[The distribution of fractional parts]{The distribution of $\xi p$ modulo $1$} \label{Sec: fractional parts}

We show that our considerations on primes $x^2+y^2+1$ in Bohr sets imply a result about the distribution of irrational multiples of such primes, in the form of Theorem \ref{theo_alphap}.\\

For proving Theorem \ref{theo_alphap}, it suffices to prove that, given an irrational $\xi>0$, there exist infinitely many integers $N\geq 1$ such that some prime $p\sim N$ of the form $x^2+y^2+1$ satisfies $\|\xi p+\kappa\|\leq \frac{N^{-\theta}}{2}$.  Let $\chi_0$ be a $1$-periodic function which is a lower bound for the characteristic function of $[-\frac{\eta}{2}, \frac{\eta}{2}]$ with $\eta=N^{-\theta}$. Specifically, as in \cite{matomaki-bombieri}, we choose $\chi_0$ so that
\begin{align*}
&0\leq \chi_0(t)\leq 1,\quad \chi_0(t)=0\quad \text{when}\quad t\not \in \left[-\frac{\eta}{2},\frac{\eta}{2}\right],\\
&\chi_0(t)=\frac{\eta}{2}+\sum_{|r|>0}c(r)e(rt)\quad \text{with}\,\,c(r)\ll \eta,\\
&\text{and}\,\,\sum_{|r|>R}|c(r)|\ll R^{-1}\quad \text{for} \quad R= \eta^{-1}(\log \eta^{-1})^C
\end{align*}
for some large constant $C$. This construction goes back to Vinogradov's work. What we want to show is that  
\begin{align}\label{eq25}
\sum_{\substack{p\sim N\\p\in \mathcal{S}+1}}\chi_0(\xi p+\kappa)\geq  \delta_0 \frac{\eta N}{(\log N)^{\frac{3}{2}}}
\end{align}
for some absolute constant $\delta_0>0$ and infinitely many $N$.  From now on, we choose a large integer $q$ satisfying $|\xi-\frac{a}{q}|\leq \frac{1}{q^2}$ for some $a$ coprime to $q$ (there are infinitely many such $q$) and take
\begin{align}\label{eq66}
N=q^2,\,\, R= \eta^{-1}(\log \eta^{-1})^C\asymp N^{\theta}(\log N^{\theta})^{C}.    
\end{align}
Concerning the term on the right-hand side of \eqref{eq25}, we note that
\begin{align*}
\sum_{n\sim N}\chi_0(\xi n+\kappa)-\frac{\eta}{2} N &\ll  \eta \sum_{0<|r|\leq R} \left|\sum_{n\sim N}e(\xi r n)\right| +\frac{N}{R}\\
&\ll \eta \sum_{0<|r|\leq R}\frac{1}{\|\xi r\|}+\eta N(\log N)^{-C}\\
&\ll \eta q\log{2q}+\eta N(\log N)^{-C}\\
&\ll \eta N(\log N)^{-C}
\end{align*}
for $2\varepsilon\leq \theta\leq \frac{1}{2}-\varepsilon$, so \eqref{eq25} takes the form 
\begin{align}\label{eq63}
\sum_{\substack{p\sim N\\p\in \mathcal{S}+1}}\chi_0(\xi p+\kappa)\geq \frac{\delta_1}{(\log N)^{\frac{3}{2}}}\sum_{n\sim N}\chi_0(\xi n +\kappa)    
\end{align}
for some absolute constant $\delta_1>0$. This is what we set out to prove.\\

\textbf{Proof of Theorem \ref{theo_alphap}.}  Pick any amenable linear polynomial, such as $L(n)=Kn+5$ with $K=6^4$. By applying Theorem \ref{t2} to $\omega_n=\chi_0(K\xi n+\kappa+5\xi)$ and $L(n)$, we see that \eqref{eq63} will follow (with $N$ replaced by $\frac{N}{K}$) once we establish Hypothesis \ref{h1} (with $\delta=(K,5-1)=4$) for this sequence $(\omega_n)$  and some parameters satisfying $\text{H}(\rho_1,\rho_2,\sigma)$ under the conditions \eqref{eq66}. Taking the definition of $\chi_0(\cdot)$ into account and making use of the classical Bombieri-Vinogradov theorem, it suffices to prove Hypothesis \ref{h1} for $\omega_n'=\sum_{0<|r|<R}c(r)e(K\xi r n)$ (with the choices \eqref{eq66}). Hence, what we must show is that
\begin{align*}
&\sum_{\substack{d\leq N^{\rho_2}\\(d,K)=1}}|\lambda_{d}^{-,\textnormal{SEM}}|\sum_{0<|r|<R}\bigg|\sum_{\substack{n\sim N\\Kn+5\in \mathbb{P}\\Kn+4\equiv 0\hspace{-0.1cm} \pmod d}}e(K\xi r n)-\frac{K}{\varphi(Kd)}\sum_{n\sim N}\frac{e(K\xi rn)}{\log(Kn)}\bigg|\quad \text{and}\\
&\sum_{\substack{d\leq N^{\rho_1}\\(d,K)=1}}|\lambda_{d}^{+,\textnormal{LIN}}|\sum_{0<|r|<R}\bigg|\sum_{\substack{\ell\leq N^{1-\varepsilon}\\(\ell,d)=1\\(\ell,K)=\delta}}g(\ell)\bigg(\sum_{\substack{n\sim N\\Kn+4=\ell p\\Kn+5\equiv 0\hspace{-0.1cm} \pmod d}}e(K\xi rn)-\frac{K}{\varphi(Kd)}\sum_{n\sim N}\frac{e(K\xi rn)}{\ell \log \frac{Kn}{\ell}}\bigg)\bigg|
\end{align*}
are $\ll \frac{N}{(\log N)^{100}}$, where $\lambda_d^{-,\textnormal{SEM}}$ has sifting parameter $z_2\ll N^{\frac{1}{\sigma}}$, while $\lambda_d^{+,\textnormal{LIN}}$ has sifting parameter $z_1\ll N^{\frac{1}{5}}$. We know that $|K\xi-\frac{a'}{q'}|\leq \frac{6^4}{q'^2}$ for some coprime $a'$ and $q'\asymp N^{\frac{1}{2}}$, so the minor arc arguments from Section \ref{Sec: hypotheses} allow replacing the previous Bombieri-Vinogradov sums (up to error $\ll N^{1-\varepsilon}$) with the sums
\begin{align}\begin{split} \label{eq99}
&\sum_{\substack{d\leq N^{\rho_2}\\(d,K)=1}}|\lambda_{d}^{-,\textnormal{SEM}}|\sum_{0<|r|<R}\bigg|\sum_{\substack{n\sim N\\Kn+5\in \mathbb{P}\\Kn+4\equiv 0\hspace{-0.1cm} \pmod d}}e(K\xi r n)\bigg|\quad \text{and}\\
&\sum_{\substack{d\leq N^{\rho_1}\\(d,K)=1}}|\lambda_{d}^{+,\textnormal{LIN}}|\sum_{0<|r|<R}\bigg|\sum_{\substack{\ell\leq N^{1-\varepsilon}\\(\ell,d)=1\\(\ell,K)=\delta}}g(\ell)\sum_{\substack{n\sim N\\Kn+4=\ell p\\Kn+5\equiv 0\hspace{-0.1cm} \pmod d}}e(K\xi rn)\bigg|.
\end{split}
\end{align}
Splitting the variables as in Subsection \ref{sub: splitting} and again employing the minor arc arguments from Section \ref{Sec: hypotheses}, the sums in \eqref{eq99} reduce to $\ll (\log N)^{10}$ sums of the same form as in Lemmas \ref{le8} and \ref{le9} with
\begin{align*}
R\leq N^{\theta}(\log N)^{C},\quad v=1,\quad q\asymp N^{\frac{1}{2}},\quad M\ll N^{\frac{1}{3}}
\end{align*} 
in the type I case, while
\begin{align*}
R\leq N^{\theta}(\log N)^{C},\quad v=1,\quad q\asymp N^{\frac{1}{2}},\quad M\in [N^{\frac{1}{2}},N^{\frac{2}{3}+\varepsilon^2}],\quad \Delta_1,\Delta_2 \quad \text{subject to} \quad \eqref{eq90}
\end{align*} 
(with $x$ replaced by $N$ in \eqref{eq90}) in the type II sums arising from the semilinear sieve weights and 
\begin{align*}
R\leq N^{\theta}(\log N)^{C},\quad v=1,\quad q\asymp N^{\frac{1}{2}},\quad M\in [N^{\frac{1}{2}},N^{\frac{3}{4}-\varepsilon}],\quad \Delta_1,\Delta_2 \quad \text{subject to} \quad \eqref{eq91}
\end{align*} 
(with $x$ replaced by $N$ in \eqref{eq91}) in the type II sums arising from the linear sieve weights.\\

From now on, we fix the values
\begin{align*}
\rho_1=\frac{1}{2}(1-4\theta)-\varepsilon,\quad \rho_2=\frac{3}{7}(1-4\theta)-\varepsilon,\quad \sigma=\frac{1}{\frac{1}{3}-2\theta}+\varepsilon.
\end{align*}
The bound offered by Lemma \ref{le8} for the type I sums we face is evidently $\ll N^{1-\varepsilon^2}$ for $\theta\leq \frac{1}{30}$. This takes care of the type I sums.\\

We turn to the type II sums that are of the same form as in Lemma \ref{le9}. Utilizing Lemma \ref{le9}, such Bombieri-Vinogradov sums are bounded by
\begin{align}\label{eq97}
\ll RN(\log N)^{1000}\left(\left(\frac{\Delta_1 M}{N}+\frac{\Delta_1\Delta_2^2}{M}\right)^{\frac{1}{2}}+\left(\frac{1}{\Delta_1}+\frac{1}{N^{\frac{1}{2}}}\right)^{\frac{1}{8}}\right)
\end{align}
when $\Delta_1\Delta_2\leq N^{\frac{1}{2}}$ and $\Delta_1\Delta_2^2\leq M$. For $R\leq N^{\theta}(\log N)^{C}$, the estimate \eqref{eq97} is $\ll N^{1-0.1\varepsilon^2}$, provided that
\begin{align}\label{eq98}
\Delta_1\leq \frac{N^{1-2\theta-\varepsilon^2}}{M}, \quad \Delta_1\Delta_2^2\leq MN^{-2\theta-\varepsilon^2},\quad \Delta_1\geq N^{0.1},\quad \theta\leq \frac{1}{80}-\varepsilon.
\end{align}
We deal with the type II sums in three cases. We will use $\rho$ to denote either $\rho_1$ or $\rho_2$.\\

\textbf{Case 1:} Suppose that $M\geq N^{1-\rho-2\theta-\varepsilon^2}, \Delta_1\geq N^{0.1}$. By taking $D=\frac{N^{1-2\theta-\varepsilon^2}}{M}$ in \eqref{eq90}-\eqref{eq91} and using the fact that $\frac{1}{\sigma}\leq \frac{1}{3}-2\theta-2\varepsilon^2$, we can indeed achieve \eqref{eq98} as long as $D\in [N^{\frac{1}{5}},N^{\rho}]$ in the case of the linear sieve and $D\in [N^{\frac{1}{3}-2\theta-2\varepsilon^2},N^{\rho}]$ in the case of the semilinear sieve. The inequality $D\leq N^{\rho}$ holds due to our lower bound on $M$. The inequality $D\geq N^{\frac{1}{5}}$ holds for $M\leq N^{\frac{3}{4}}$, which is true in the linear case. Similarly, the inequality $D\geq N^{\frac{1}{3}-2\theta-2\varepsilon^2}$ reduces to $M\leq N^{\frac{2}{3}+\varepsilon^2}$,  and this holds in the semilinear case. Therefore, in this case \eqref{eq98} is always valid, which means that our type II sums are $\ll N^{1-0.1\varepsilon^2}$, which is what we wanted.\\

\textbf{Case 2:} Suppose that $M\geq N^{1-\rho-2\theta-\varepsilon^2}, \Delta_1<N^{0.1}$. In this case we know that $\Delta_2=1$ from \eqref{eq90} and \eqref{eq91}. Now, choosing $F_2$ in Lemma \ref{le9}, we obtain for the type II Bombieri-Vinogradov sum the bound
\begin{align*}
\ll RN\Delta_1\left(\frac{1}{N^{\frac{1}{4}}}+\frac{1}{M^{\frac{1}{2}}}+\frac{M}{N}+\frac{N^{\frac{1}{4}}}{(RN)^{\frac{1}{2}}}\right)^{\frac{1}{2}}\ll RN\Delta_1N^{-\frac{1}{8}}\ll N^{0.999}
\end{align*}
when $\theta\leq \frac{1}{50}$.\\

\textbf{Case 3:} Suppose finally that $M< N^{1-\rho-2\theta-\varepsilon^2},\Delta_1\geq N^{0.1}$. Similarly as in Case 3 of Subsection \ref{sub: min sem}, we may take $\Delta_1=N^{\rho}$, $\Delta_2=1$. Again we require this choice to fulfill \eqref{eq98}. The first constraint in \eqref{eq98} follows directly from our upper bound on $M$. Since $M\geq N^{\frac{1}{2}}$, the second constraint in \eqref{eq98} holds for $\rho\leq \frac{1}{2}-2\theta-\varepsilon^2$, which certainly holds for our choices of $\rho_1$ and $\rho_2$. This means that also in Case 3 we get good enough bounds for the type II sums. Putting everything together, in each of the Cases 1-3 we get a good enough bound for the type II  sums.\\

Combining the analyses of the Cases 1-3, we see that Theorem \ref{theo_alphap} will follow with exponent $\theta$ if $\text{H}(\rho_1,\rho_2,\sigma)$ is true for $\sigma=\frac{1}{\frac{1}{3}-2\theta}+\varepsilon$, $\rho_1=\frac{1}{2}(1-4\theta)-\varepsilon$ and $\rho_2=\frac{3}{7}(1-4\theta)-\varepsilon$, provided that $\theta\leq \frac{1}{80}-\varepsilon$. By continuity, it suffices to check $\text{H}(\frac{1}{2}(1-4\theta),\frac{3}{7}(1-4\theta),\frac{1}{\frac{1}{3}-2\theta})$ for $\theta=\frac{1}{80}$, and this holds by a numerical computation (the difference between the left and right side of \eqref{eq96} is then $>10^{-3}$). This completes the proof of Theorem \ref{theo_alphap}.\qedd

\bibliography{refs_goldbach}{}
\bibliographystyle{plain}

\textsc{Department of Mathematics and statistics, University of Turku, 20014 Turku, Finland}\\
\textit{Email address:} \textup{\texttt{joni.p.teravainen@utu.fi}}

\end{document}